\input amstex
\documentstyle{amsppt}
\magnification=\magstep 1
\document
\vsize6.7in

\chardef\oldatsign=\catcode`\@
\catcode`\@=11
\newif\ifdraftmode			
\global\draftmodefalse



%
\font@\twelverm=cmr12 
\font@\twelvei=cmmi12 \skewchar\twelvei='177 
\font@\twelvesy=cmsy10 scaled\magstep1 \skewchar\twelvesy='060 
\font@\twelveex=cmex10 scaled\magstep1 
\font@\twelvemsa=msam10 scaled\magstep1 
\font@\twelvemsb=msbm10 scaled\magstep1 
\font@\twelvebf=cmbx12 
\font@\twelvett=cmtt12 
\font@\twelvesl=cmsl12 
\font@\twelveit=cmti12 
\font@\twelvesmc=cmcsc10 scaled\magstep1 
%
%
\font@\ninerm=cmr9 
\font@\ninei=cmmi9 \skewchar\ninei='177 
\font@\ninesy=cmsy9 \skewchar\ninesy='60 
\font@\ninemsa=msam9
\font@\ninemsb=msbm9
\font@\ninebf=cmbx9
%
%
%
\font@\ttlrm=cmbx12 scaled \magstep2 
\font@\ttlsy=cmsy10 scaled \magstep3 
\font@\tensmc=cmcsc10 
%
%
\def\normaltype{
	\def\pointsize@{12}%
	\abovedisplayskip18\p@ plus5\p@ minus9\p@
	\belowdisplayskip18\p@ plus5\p@ minus9\p@
	\abovedisplayshortskip1\p@ plus3\p@
	\belowdisplayshortskip9\p@ plus3\p@ minus4\p@
	\textonlyfont@\rm\twelverm
	\textonlyfont@\it\twelveit
	\textonlyfont@\sl\twelvesl
	\textonlyfont@\bf\twelvebf
	\textonlyfont@\smc\twelvesmc
	\ifsyntax@
		\def\big##1{{\hbox{$\left##1\right.$}}}%
	\else
		\let\big\twelvebig@
 \textfont0=\twelverm \scriptfont0=\ninerm \scriptscriptfont0=\sevenrm
 \textfont1=\twelvei  \scriptfont1=\ninei  \scriptscriptfont1=\seveni
 \textfont2=\twelvesy \scriptfont2=\ninesy \scriptscriptfont2=\sevensy
 \textfont3=\twelveex \scriptfont3=\twelveex  \scriptscriptfont3=\twelveex
 \textfont\itfam=\twelveit \def\it{\fam\itfam\twelveit}%
 \textfont\slfam=\twelvesl \def\sl{\fam\slfam\twelvesl}%
 \textfont\bffam=\twelvebf \def\bf{\fam\bffam\twelvebf}%
 \scriptfont\bffam=\ninebf \scriptscriptfont\bffam=\sevenbf
 \textfont\ttfam=\twelvett \def\tt{\fam\ttfam\twelvett}%
 \textfont\msafam=\twelvemsa \scriptfont\msafam=\ninemsa
 \scriptscriptfont\msafam=\sevenmsa
 \textfont\msbfam=\twelvemsb \scriptfont\msbfam=\ninemsb
 \scriptscriptfont\msbfam=\sevenmsb
	\fi
 \normalbaselineskip=\twelvebaselineskip
 \setbox\strutbox=\hbox{\vrule height12\p@ depth6\p@
      width0\p@}%
 \normalbaselines\rm \ex@=.2326ex%
}
%
%
%
\def\smalltype{
	\def\pointsize@{10}%
	\abovedisplayskip12\p@ plus3\p@ minus9\p@
	\belowdisplayskip12\p@ plus3\p@ minus9\p@
	\abovedisplayshortskip\z@ plus3\p@
	\belowdisplayshortskip7\p@ plus3\p@ minus4\p@
	\textonlyfont@\rm\tenrm
	\textonlyfont@\it\tenit
	\textonlyfont@\sl\tensl
	\textonlyfont@\bf\tenbf
	\textonlyfont@\smc\tensmc
	\ifsyntax@
		\def\big##1{{\hbox{$\left##1\right.$}}}%
	\else
		\let\big\tenbig@
	\textfont0=\tenrm \scriptfont0=\sevenrm \scriptscriptfont0=\fiverm 
	\textfont1=\teni  \scriptfont1=\seveni  \scriptscriptfont1=\fivei
	\textfont2=\tensy \scriptfont2=\sevensy \scriptscriptfont2=\fivesy 
	\textfont3=\tenex \scriptfont3=\tenex \scriptscriptfont3=\tenex
	\textfont\itfam=\tenit \def\it{\fam\itfam\tenit}%
	\textfont\slfam=\tensl \def\sl{\fam\slfam\tensl}%
	\textfont\bffam=\tenbf \def\bf{\fam\bffam\tenbf}%
	\scriptfont\bffam=\sevenbf \scriptscriptfont\bffam=\fivebf
	\textfont\msafam=\tenmsa
	\scriptfont\msafam=\sevenmsa
	\scriptscriptfont\msafam=\fivemsa
	\textfont\msbfam=\tenmsb
	\scriptfont\msbfam=\sevenmsb
	\scriptscriptfont\msbfam=\fivemsb
		\textfont\ttfam=\tentt \def\tt{\fam\ttfam\tentt}%
	\fi
 \normalbaselineskip 14\p@
 \setbox\strutbox=\hbox{\vrule height10\p@ depth4\p@ width0\p@}%
 \normalbaselines\rm \ex@=.2326ex%
}

\def\titletype{
	\def\pointsize@{17}%
	\textonlyfont@\rm\ttlrm
	\ifsyntax@
		\def\big##1{{\hbox{$\left##1\right.$}}}%
	\else
		\let\big\twelvebig@
		\textfont0=\ttlrm \scriptfont0=\twelverm
		\scriptscriptfont0=\tenrm
		\textfont2=\ttlsy \scriptfont2=\twelvesy
		\scriptscriptfont2=\tensy
	\fi
	\normalbaselineskip 25\p@
	\setbox\strutbox=\hbox{\vrule height17\p@ depth8\p@ width0\p@}%
	\normalbaselines
	\rm
	\ex@=.2326ex%
}

\def\tenbig@#1{
	{%
		\hbox{%
			$%
			\left
			#1%
			\vbox to8.5\p@{}%
			\right.%
			\n@space
			$%
		}%
	}%
}

\def\twelvebig@#1{%
	{%
		\hbox{%
			$%
			\left
			#1%
			\vbox to10.2\p@{}
			\right.%
			\n@space
			$%
		}%
	}%
}

%
%
%
%
%
\newif\ifl@beloutopen
\newwrite\l@belout
\newread\l@belin

\global\let\currentfile=\jobname

\def\getfile#1{%
	\immediate\closeout\l@belout
	\global\l@beloutopenfalse
	\gdef\currentfile{#1}%
	\input #1%
	\par
	\newpage
}

\def\getxrefs#1{%
	\bgroup
		\def\gobble##1{}
		\edef\list@{#1,}%
		\def\gr@boff##1,##2\end{
			\openin\l@belin=##1.xref
			\ifeof\l@belin
			\else
				\closein\l@belin
				\input ##1.xref
			\fi
			\def\list@{##2}%
			\ifx\list@\empty
				\let\next=\gobble
			\else
				\let\next=\gr@boff
			\fi
			\expandafter\next\list@\end
		}%
		\expandafter\gr@boff\list@\end
	\egroup
}

\def\testdefined#1#2#3{%
	\expandafter\ifx
	\csname #1\endcsname
	\relax
	#3%
	\else #2\fi
}

\def\document{%
	\minaw@11.11128\ex@ 
	\def\alloclist@{\empty}%
	\def\fontlist@{\empty}%
	\openin\l@belin=\jobname.xref	
	\ifeof\l@belin\else
		\closein\l@belin
		\input \jobname.xref
	\fi
}

\def\getst@te#1#2{%
	\edef\st@te{\csname #1s!#2\endcsname}%
	\expandafter\ifx\st@te\relax
		\def\st@te{0}%
	\fi
}

\def\setst@te#1#2#3{%
	\expandafter
	\gdef\csname #1s!#2\endcsname{#3}%
}

\outer\def\setupautolabel#1#2{%
	\def\newcount@{\global\alloc@0\count\countdef\insc@unt}	
	\def\newtoks@{\global\alloc@5\toks\toksdef\@cclvi}
	\expandafter\newcount@\csname #1Number\endcsname
	\expandafter\global\csname #1Number\endcsname=1%
	\expandafter\newtoks@\csname #1l@bel\endcsname
	\expandafter\global\csname #1l@bel\endcsname={#2}%
}

\def\reflabel#1#2{%
	\testdefined{#1l@bel}
	{
		\getst@te{#1}{#2}%
		\ifcase\st@te
			???
			\message{Unresolved forward reference to
				label #2. Use another pass.}%
		\or	
			\setst@te{#1}{#2}2
			\csname #1l!#2\endcsname 
		\or	
			\csname #1l!#2\endcsname 
		\or	
			\csname #1l!#2\endcsname 
		\fi
	}{
		{\escapechar=-1 
		\errmessage{You haven't done a
			\string\\setupautolabel\space for type #1!}%
		}%
	}%
}

{\catcode`\{=12 \catcode`\}=12
	\catcode`\[=1 \catcode`\]=2
	\xdef\Lbrace[{]
	\xdef\Rbrace[}]%
]%

\def\setlabel#1#2{%
	\testdefined{#1l@bel}
	{
		\edef\templ@bel@{\expandafter\the
			\csname #1l@bel\endcsname}%
		\def\@rgtwo{#2}%
		\ifx\@rgtwo\empty
		\else
			\ifl@beloutopen\else
				\immediate\openout\l@belout=\currentfile.xref
				\global\l@beloutopentrue
			\fi
			\getst@te{#1}{#2}%
			\ifcase\st@te
			\or	
			\or	
				\edef\oldnumber@{\csname #1l!#2\endcsname}%
				\edef\newnumber@{\templ@bel@}%
				\ifx\newnumber@\oldnumber@
				\else
					\message{A forward reference to label 
						#2 has been resolved
						incorrectly.  Use another
						pass.}%
				\fi
			\or	
				\errmessage{Same label #2 used in two
					\string\setlabel s!}%
			\fi
			\expandafter\xdef\csname #1l!#2\endcsname
				{\templ@bel@}
			\setst@te{#1}{#2}3%
			\immediate\write\l@belout 
				{\string\expandafter\string\gdef
				\string\csname\space #1l!#2%
				\string\endcsname
				\Lbrace\templ@bel@\Rbrace
				}%
			\immediate\write\l@belout 
				{\string\expandafter\string\gdef
				\string\csname\space #1s!#2%
				\string\endcsname
				\Lbrace 1\Rbrace
				}%
		\fi
		\templ@bel@	
		\expandafter\ifx\envir@end\endref 
			\gdef\marginalhook@{\marginal{#2}}%
		\else
			\marginal{#2}
		\fi
		\expandafter\global\expandafter\advance	
			\csname #1Number\endcsname
			by 1 %
	}{
		{\escapechar=-1
		\errmessage{You haven't done a \string\\setupautolabel\space
			for type #1!}%
		}%
	}%
}


\newcount\SectionNumber
\setupautolabel{t}{\number\SectionNumber.\number\tNumber}
\setupautolabel{r}{\number\rNumber}
\setupautolabel{T}{\number\TNumber}

\define\rref{\reflabel{r}}
\define\tref{\reflabel{t}}

\define\tnum{\setlabel{t}}
\define\rnum{\setlabel{r}}

%
\def\strutdepth{\dp\strutbox}%
\def\strutheight{\ht\strutbox}%

\newif\iftagmode
\tagmodefalse

\let\old@tagform@=\tagform@
\def\tagform@{\tagmodetrue\old@tagform@}

\def\marginal#1{%
	\ifvmode
	\else
		\strut
	\fi
	\ifdraftmode
		\ifmmode
			\ifinner
				\let\Vorvadjust=\Vadjust
			\else
				\let\Vorvadjust=\vadjust
			\fi
		\else
			\let\Vorvadjust=\Vadjust
		\fi
		\iftagmode	
			\llap{%
				\smalltype
				\vtop to 0pt{%
					\pretolerance=2000
					\tolerance=5000
					\raggedright
					\hsize=.72in
					\parindent=0pt
					\strut
					#1%
					\vss
				}%
				\kern.08in
				\iftagsleft@
				\else
					\kern\hsize
				\fi
			}%
		\else
			\Vorvadjust{%
				\kern-\strutdepth 
				{%
					\smalltype
					\kern-\strutheight 
					\llap{%
						\vtop to 0pt{%
							\kern0pt
							\pretolerance=2000
							\tolerance=5000
							\raggedright
							\hsize=.5in
							\parindent=0pt
							\strut
							#1%
							\vss
						}%
						\kern.08in
					}%
					\kern\strutheight
				}%
				\kern\strutdepth
			}
		\fi
	\fi
}


\newbox\Vadjustbox

\def\Vadjust#1{
	\global\setbox\Vadjustbox=\vbox{#1}%
	\ifmmode
		\ifinner
			\innerVadjust
		\fi		
	\else
		\innerVadjust
	\fi
}

\def\innerVadjust{%
	\def\nexti{\aftergroup\innerVadjust}%
	\def\nextii{%
		\ifvmode
			\hrule height 0pt 
			\box\Vadjustbox
		\else
			\vadjust{\box\Vadjustbox}%
		\fi
	}%
	\ifinner
		\let\next=\nexti
	\else
		\let\next=\nextii
	\fi
	\next
}%

\global\let\marginalhook@\empty

\def\endref{%
\setbox\tw@\box\thr@@
\makerefbox?\thr@@{\endgraf\egroup}%
  \endref@
  \endgraf
  \endgroup
  \keyhook@
  \marginalhook@
  \global\let\keyhook@\empty 
  \global\let\marginalhook@\empty 
}

\catcode`\@=\oldatsign

\nologo
\vsize6.7in


\topmatter

\title The resolution of the bracket powers of the maximal ideal in a diagonal hypersurface ring\endtitle
  \leftheadtext{Kustin, Rahmati, and Vraciu}
\rightheadtext{Bracket powers  in a diagonal hypersurface ring}
\author   
Andrew R. Kustin\footnote{Supported in part by the National Security Agency. \phantom{xxxxxxxxxxxxxxxxxxxxxxxxxxxxxxxxxxxxx}}, Hamid Rahmati, and Adela Vraciu\footnote{Supported in part by the National Security Agency.\phantom{xxxxxxxxxxxxxxxxxxxxxxxxxxxxxxxxxxxxx}}\endauthor

 \address
Mathematics Department,
University of South Carolina,
Columbia, SC 29208\endaddress
\email kustin\@math.sc.edu \endemail
\address Mathematics Department,
Syracuse University,
Syracuse, NY 13244\endaddress 
\email     hrahmati\@syr.edu \endemail
 \address
Mathematics Department,
University of South Carolina,
Columbia, SC 29208\endaddress
\email adela$\underline{\phantom{x}}$vraciu\@yahoo.com \endemail
\keywords Almost Complete Intersection, Enumeration of Plane Partitions, Frobenius Periodicity, Frobenius Power, Hilbert-Burch Matrix, Hilbert-Kunz Functions, Pfaffians, Socle Degrees, Syzygy, Syzygy Gap, Weak Lefschetz Property \endkeywords
\subjclass 13A35, 13D02, 13H10, 13C40, 13E10, 13D40\endsubjclass
\abstract Let \pmb k be a field. For each pair of positive integers $(n,N)$, we resolve $Q=R/(x^N,y^N,z^N)$
as a module over the ring $R=\pmb k[x,y,z]/(x^n+y^n+z^n)$. Write  $N$  in the form  $N=a n+r$  for integers $a$ and  $r$, with $r$ between $0$ and $n-1$.   If  $n$  does not divide  $N$  and  the characteristic of $\pmb k$ is fixed, then the value of $a$ determines whether $Q$ has finite or infinite projective dimension. If $Q$  has infinite projective dimension, then value of  $r$, together with the parity of $a$, determines the periodic part of the infinite resolution. When $Q$ has infinite projective dimension we give an explicit presentation for the module of first syzygies of $Q$. This presentation is quite complicated. We also give an explicit presentation the module of second syzygies for $Q$. This presentation is remarkably uncomplicated. We use linkage to find an explicit generating set for the grade three Gorenstein ideal $(x^N,y^N,z^N):(x^n+y^n+z^n)$ in the polynomial ring $\pmb k[x,y,z]$.

The question ``Does $Q$ have finite projective dimension?'' is intimately connected to the question ``Does $\pmb k[X,Y,Z]/(X^a,Y^a,Z^a)$ have the Weak Lefschetz Property?''. The second question is connected to the enumeration of plane partitions. 

When the field $\pmb k$ has positive characteristic, we investigate three questions about the Frobenius powers $F^t(Q)$ of $Q$. When does there exist a pair $(n,N)$ so that $Q$ has infinite projective dimension and $F(Q)$ has finite projective dimension?  Is the tail of the resolution of the Frobenius power $F^t(Q)$   eventually a periodic function of $t$, (up to shift)? 
In particular, we   exhibit  a situation where the tail of the resolution of $F^t(Q)$, after shifting, is periodic as a function of $t$, with an arbitrarily large period. 
Can one use socle degrees to predict that the tail of the resolution of $F^t(Q)$ is a shift of the tail of the resolution of  $Q$? 
\endabstract

\endtopmatter

\document

\centerline{Table of Contents}

\halign{
#\hfil&\quad#\hfil&\quad#\hfil&\quad#\hfil&\quad#\hfil\cr
0.&Introduction.\cr
1.&Terminology, notation, and elementary results. \cr
2.&The technique for resolving $Q_{\pmb k,n,N}$.\cr
3.&The resolutions of $Q_{\pmb k,n,N}$ when $\lfloor\frac Nn\rfloor$ is in $S_c$.\cr
4.&The proof of Theorem \tref{4}.\cr
5.&The projective dimension of $Q_{\pmb k,n,N}$ is finite when $\lfloor\frac Nn\rfloor$ is in $T_p$.\cr
6.&The set of non-negative integers may be partitioned as  $S_p\cup T_p$.\cr
7.&The proof of Lemma \tref{NandI}.\cr
8.&Applications to Frobenius powers.\cr
9.&Two variables.\cr
&References.\cr
}

\SectionNumber=0\tNumber=1
\heading Section \number\SectionNumber. \quad Introduction.
\endheading

Throughout this paper the diagonal hypersurface ring $R_{\pmb k,n}$ 
is the ring
 $$R_{\pmb k,n}=\pmb k[x,y,z]/(x^n+y^n+z^n),\tag\tnum{DHSR}$$ 
where  $\pmb k$ is a field and  $n$ is a positive integer.
For each positive integer $N$, we resolve the quotient ring  $$Q_{\pmb k,n,N}=R_{\pmb k,n}/(x^N,y^N,z^N),\tag\tnum{DHSR.1}$$ as a module over $R_{\pmb k,n}$. 
For each real number $\alpha$ which is not of the form $b+\frac12$ for some integer $b$, let $\{\alpha\}$ represent the integer which is closest to $\alpha$. In particular, if $b$ is an integer, then
$$\{{\tsize \frac b3}\}=\cases \frac b3&\text{if $b\equiv 0\mod 3$}\\
\frac {b-1}3&\text{if $b\equiv 1\mod 3$}\\\frac {b+1}3&\text{if $b\equiv 2\mod 3$}.\endcases\tag\tnum{rnd}$$
 
One of our main results is
\proclaim{Theorem  \tref{!!!}}Let $\pmb k$ be a field of characteristic $c$ and let $n$ and $N$ be positive integers. Then $\operatorname{pd}_{R_{\pmb k,n}} Q_{\pmb k,n,N}$ is finite if and only if at least one of the following conditions hold{\rm:}
\roster\item $n$ divides $N$, or
\item $c=2$ and $n\le N$, or
\item $c=p$ is an odd prime and there exist an odd integer $J$ and a power $q=p^e$ of $p$ with $e\ge 1$ and $|Jq-\frac Nn|< \{\frac q3\}$.\endroster
\flushpar In particular, if $c=0$, then $\operatorname{pd}_{R_{\pmb k,n}} Q_{\pmb k,n,N}$ is finite if and only if $n$ divides $N$. 
\endproclaim

 We use $\operatorname{pd}$ to mean {\it projective dimension}. When our meaning is clear, we write $R$ and $Q$ in place of $R_{\pmb k,n}$ and $Q_{\pmb k,n,N}$, respectively. Write $N$ in the form $N=\theta n+r$, for integers $\theta$ and $r$, with $0\le r\le n-1$. If $n$ does not divide $N$ and  the characteristic of $\pmb k$ is fixed, then the value of $\theta$ determines whether $Q_{\pmb k,n,N}$ has finite or infinite projective dimension. If $Q_{\pmb k,n,N}$ has infinite projective dimension, then value of $r$, together with the parity of $\theta$, determines the periodic part of the infinite resolution.

There are three steps in the proof of Theorem \tref{!!!}. In Definitions \tref{!D18.1}  and \tref{D5.1} we   define sets of non-negative integers $S_c$ and $T_c$  which depend   on the characteristic $c$ of $\pmb k$. (We use $c$ to be the characteristic of the field $\pmb k$; thus $c$ is either $0$ or a positive prime integer $p$.)
The first step is  Theorem \tref{I-MAIN}, where  we record an explicit infinite homogeneous minimal resolution of $Q$ when $\theta$ is in $S_c$ and $1\le r$. 
The second step is Theorem \tref{HW672}  where we identify a set of integers $T_p$ such    that if $\theta\in T_p$, then $\operatorname{pd}_RQ<\infty$. The proof in the finite projective dimension case is less explicit than the proof in the infinite projective dimension case. In the finite projective dimension  case we use a Theorem of Brenner and Kaid which identifies those integers $a$ for which the Hilbert-Burch matrix $\operatorname{HB}_a$ for $[X^a,Y^a,(X+Y)^a]$ in $\pmb k[X,Y]$ is unbalanced. (See \tref{A1} for the definitions.) The paper of Brenner and Kaid \cite{\rref{BK}} ties this unbalanced behavior of the Hilbert-Burch matrix to the absence of the  Weak Lefschetz Property (WLP) in the ring $\bar{ \pmb k}[X,Y,Z]/(X^a,Y^a,Z^a)$, where $\bar{\pmb k}$ is the algebraic closure of $\pmb k$. (We notice that Li and Zanello \cite{\rref{LZ}} have found ``a surprising, and still combinatorially obscure, connection'' between the monomial complete intersection ideals in three variables which satisfy the WLP, as a function of the characteristic of the base field, and the enumeration of plane partitions.)
The proof of the Brenner-Kaid result uses a Theorem of Han \cite{\rref{Han}}. An easier proof of Han's Theorem was given by Monsky \cite{\rref{M}}. The result in the present paper which connects unbalanced Hilbert-Burch matrices to $\operatorname{pd}_RQ<\infty$ is Theorem \tref{AVT}. The set $T_p$ contains and is slightly larger than     the set of indices $a$ where $\operatorname{HB}_a$ is unbalanced. 

Each set $S_p$ is defined recursively as a union of intervals. Each set $T_p$ is defined explicitly as a union of intervals. It is not immediately obvious that the union of $S_p$ and $T_p$ is the set of all non-negative integers. The proof of this fact, which is the third and final step in the proof of Theorem \tref{!!!},  is given in Theorem \tref{H72+}. The proof of Theorem \tref{H72+} makes use of a base-$p$ expansion of integers using unusual coefficients; see Notation \tref{N5.1}.

Three ingredients are used in the proof of Theorem \tref{I-MAIN}, which is the infinite projective dimension case. The linkage part of the argument, which is  similar to the work of Buchsbaum-Eisenbud (see for example \cite{\rref{BE}, Thm.~5.3}), is carried out in Section 2. In order to apply Section 2, one must find a matrix $\psi$ and a unit $u$ so that the quadratic equation that we have called (\tref{hyp20}) is satisfied. In Section 3 we propose a solution to the problem posed by (\tref{hyp20}). We prove that our proposed solution works over $\Bbb Z$ in Section 7. In Section 4 we prove that when $\theta\in S_p$, then the solution of (\tref{hyp20}) which works over $\Bbb Z$ can be made to work over $\pmb k$.

The critical parts of the  proof of Theorem \tref{I-MAIN} take place in the polynomial ring $P=\pmb k[x,y,z]$. In particular, we determine a generating set for the grade three Gorenstein ideal $$(x^N,y^N,z^N)\:\! (x^n+y^n+z^n)\tag\tnum{GknN}$$ in $P$. We prove (see Remark \tref{2.14}) that $Q_{\pmb k,n,N}$ has infinite projective dimension over $R_{\pmb k,n}$ if and only if  
the ideal of (\tref{GknN}) of $P$ has seven generators.

One consequence of our work (see Corollary \tref{I+}) is that the module of second syzygies of the $R_{\pmb k,n}$-module $Q_{\pmb k,n,N}$ is remarkably uncomplicated. Each such second syzygy module is either free or  is presented by a matrix of the form $$\varphi_{r,s}=\bmatrix 0&z^{r}&-y^{r}&x^{s} \\
-z^{r}&0&x^{r}&y^{s} \\
y^{r}&-x^{r}&0&z^{s} \\
-x^{s}&-y^{s}&-z^{s}&0 \endbmatrix, $$ where $r$ and $s$ are positive integers with $r+s=n$. This fact is somewhat surprizing because the presentation matrix for the module of first syzygies of $Q_{\pmb k,n,N}$ is very complicated. Its entries involve high degrees, many terms, and intricate coefficients given in terms of binomial coefficients.  The module of first syzygies of $Q_{\pmb k,n,N}$ is studied in a very recent preprint of Brenner and Kaid \cite{\rref{BK2}}. They use the syzygy bundle $\operatorname{Syz}(x^p,y^p,z^p)$ on $\operatorname{Proj}(R_{\pmb k,n})$, where $\pmb k$ is an algebraically closed field of characteristic $p$ to create vector bundles $\Cal E$ with $\Cal E$ isomorphic to $F^*(\Cal E)$, where $F\: R_{\pmb k,n}\to R_{\pmb k,n}$ is the Frobenius functor.

Recall that if 
$R$ is a ring of positive characteristic $p$, 
$J$ is an  ideal in $R$, and $q=p^e$, for some positive integer $e$, then the $e^{\text{th}}$ Frobenius power of $J$ is the ideal $J^{[q]}$ generated by all $j^q$ with $j\in J$.  Furthermore, the $e^{\text{th}}$ Frobenius power of the $R_{\pmb k,n}$-module $Q_{\pmb k,n,N}$ is  $Q_{\pmb k,n,qN}$.
Three questions pertaining to Frobenius powers are investigated in Section 8. 
The first question is 
``When does the Cohen-Macaulay local ring $(R,\frak m)$ have an $\frak m$-primary ideal $J$ so that $R/J$ has infinite projective dimension but the Frobenius power $R/J^{[p]}$ has finite projective dimension?'' It is shown that  if $R$ is not $F$-injective, then such an ideal exists. Furthermore, when $R$ is one of the rings $R_{\pmb k,n}$, then there exists an integer $N$ so that $J=(x^N,y^N,z^N)$ has the above property if and only if $R$ is not $F$-injective.
The second questions is ``Is the tail of the resolution of the Frobenius power $Q_{\pmb k,n,p^tN}$ (up to shift)  eventually a periodic function of $t$?'' The answer is yes. The third question is ``Can one use socle degrees to predict that the tail of the resolution of $Q_{\pmb k,n,p^tN}$ is a shift of the tail of the resolution of  $Q_{k,n,N}$?'' The third question was answered in \cite{\rref{KV}} in the case that both modules have finite projective dimension (hence the infinite tail of both resolutions  is zero). It is shown in \cite{\rref{KU}} how the socle degrees can be used to predict that the tail of the resolution of $Q_{\pmb k,n,p^tN}$ is a shift of the tail of the resolution of  $Q_{\pmb k,n,N}$, {\bf as a graded module}. We show that the condition of \cite{\rref{KU}}, when applied to $Q_{\pmb k,n,N}$, actually produces an isomorphism of complexes {\bf with differential}. 

In Section 9 we resolve  the module $\bar Q_{\pmb k,n,N}=\bar R_{\pmb k,n}/(x^N,y^N)$ over the two variable diagonal hypersurface   ring $\bar R=\pmb k[x,y]/(x^n+y^n)$ and we prove that the tail of the resolution of the Frobenius power $F^t(\bar Q_{\pmb k,n,N})$ is isomorphic to a shift of the tail of the resolution of $\bar Q_{\pmb k,n,N}$ as a graded module if and only if these objects are isomorphic as complexes if and only if the socle of $F^t(\bar Q_{\pmb k,n,N})$ is isomorphic to a shift of the socle of $\bar Q_{\pmb k,n,N}$. 
 We conclude by exhibiting a situation where the tail of the resolution of $F^t(\bar Q_{\pmb k,n,N})$, after shifting, is periodic as a function of $t$, with an arbitrarily large period.  
The results of Section 9 were announced in \cite{\rref{KU}}.

\SectionNumber=1\tNumber=1
\heading Section \number\SectionNumber. \quad Terminology, notation, and elementary results. 
\endheading 

In this section we gather the terminology notation, and elementary results which are used throughout the paper. We begin by defining the sets $S_c$ and $T_c$ which appear in the statement of Theorem \tref{???} and in the proof of Theorem \tref{!!!}. {\bf In this discussion $c$ is the characteristic of a field so $c$ is either zero or a positive prime integer $p$.}

\definition{Definition \tnum{!D18.1}} 
For each field $\pmb k$, we define a set of non-negative integers $S_c$, where $c$ is the characteristic of $\pmb k$. If $c$ is a positive prime integer $p$, with $p\ge 3$, then we also define a second  set of non-negative integers $D_p$.
\smallskip\flushpar{\bf(1)} The set $S_0$ is the set of all non-negative integers.
\smallskip\flushpar{\bf(2)} The set $S_2$ is $\{0\}$.
\smallskip\flushpar{\bf(3)} The sets $D_3$ and $S_3$ are defined recursively.
\smallskip  {\bf(a)} Define $D_3$ as follows. The number  $0$ is in $D_3$ and if $d\in D_3$, then $3d$ and $3d+2$ are also in $D_3$. 
\smallskip  {\bf(b)} Define $S_3$ as follows. The number  $0$ is in $S_3$ and if $a$ is an even element of $S_3$, then  $3a$, $3a+1$, $3a+4$, and $3a+5$ all are in $S_3$. 
\smallskip\flushpar{\bf(4)} Let 
 $p\ge5$ be prime. Define  $\pi_p=\pi$ to be the largest integer with 
 $\pi_p< \frac{p}3$. The sets $D_p$ and $S_p$ are defined recursively.
\smallskip  {\bf(a)} Define $D_p$ as follows. The closed interval of integers $[0,\pi]\subseteq D_p$; and if $d$ is a non-negative element of $D_p$, then $[pd-\pi,pd+\pi]\subseteq D_p$.
\smallskip {\bf(b)} Define $S_p$ as follows. The closed interval of integers $[0,2\pi]\subseteq S_p$; and if $\theta$ is a non-negative  even element of $S_p$, then $[p\theta-2\pi-1,p\theta+2\pi]\subseteq S_p$.\enddefinition
\remark{Remarks \tnum{!R18.1}}\smallskip\flushpar{\bf(1)} The sets $D_p$ and $S_p$ are related as follows:
 $$\cases D_p=\tsize\{\frac 12 \theta\mid \theta \text{ is an even element of $S_p$}\},&\text{for $3\le p$,}\\ S_3=\{2d\mid d\in D_3\}\cup \{2d+1\mid d\in D_3\},\\
 S_p=\{2d\mid d\in D_p\}\cup \{2d-1\mid d\in D_p\},&\text{for $5\le p$}.\endcases\tag\tnum{ae}$$One could also
define exactly one of the sets $D_p$ or $S_p$ in a recursive manner and then define the other set using (\tref{ae}).
\smallskip\flushpar{\bf(2)} If $a$ is an odd integer and $p\ge 3$ is  a prime integer, then
$$\cases a\in S_p\iff a-1\in S_p,&\text{if $p=3$, and}\\a\in S_p\iff a+1\in S_p,&\text{if $p\ge 5$}.\endcases$$
\endremark
\example{Example \tnum{!E18.2}}We see that
$$\matrix\format\l&\qquad \l\\
&D_3=\{0,2,6,8,18,20,24,26,54,56,60,62,72,74,78,80,\dots\},\\
&S_3={}\{0,1,4,5,12,13,16,17,36,37,40,41,48,49,52,53,\dots\},\\
 \pi_5=1,&D_5=[0,1]\cup[4,6]\cup[19,21]\cup[24,26]\cup\dots,\\&S_5=[0,2]\cup[7,12]\cup[37,42]\cup[47,52]\cup\dots,\\
\pi_7=2,&D_7=[0,2]\cup[5,9]\cup [12,16]\cup[33,37]\cup \dots,\\
&S_7=[0,4]\cup[9,18]\cup [23,32]\cup[65,74]\cup \dots,\\
\pi_{11}=3,&D_{11}=[0,3]\cup[8,14]\cup [19,25]\cup[30,36]\cup \dots,\\ 
&S_{11}=[0,6]\cup[15,28]\cup [37,50]\cup[59,72]\cup \dots,\\ 
\pi_{13}=4, &D_{13}=[0,4]\cup[9,17]\cup [22,30]\cup[35,43]\cup \dots, \ \text{and}\\
&S_{13}=[0,8]\cup[17,34]\cup [43,60]\cup[69,86]\cup \dots\ .\endmatrix $$
\endexample
In Remarks \tref{R1.1} and \tref{O-10-13} we offer an alternate, explicit, description of $S_p$ for $p\ge 3$. 

\definition{Notation \tnum{N5.1}}
Let $p\ge 5$ be a prime number.  We will write integers in base-$p$,   using even digits of the form 
$-2k,\ldots,-2, 0, 2, \ldots, 2k$, where $2k=v$ is the largest even integer such that $2k<2p/3$, and odd digits of the form $-u,\ldots,-1, 1, \ldots, u, $ where $u$ is the largest odd integer such that $u < p/3$. It is easy to see that every integer can be written uniquely in the form
$a_0 + a_1p + \ldots +a_tp^t$, with $a_0, a_1, \ldots, a_t$ digits as above.\enddefinition

\remark{Remark \tnum{R1.1}} Fix a prime number $p\ge 5$. Recall the set $S_p$ of Definition \tref{!D18.1}. It is not difficult to see that if $m$ is a non-negative even integer, then $m\in S_p$ if and only if the base-$p$ expansion of $m$, in the sense of Notation \tref{N5.1}, involves only even digits.\endremark 

\proclaim{Remark \tnum{O-10-13}}If $a$ is a non-negative integer, then $a\in S_3$ if and only if there are coefficients $\epsilon_i$ and $\epsilon$ in the set $\{0,1\}$ such that
$$a=\epsilon +4\sum_{i=0}^r \epsilon_i3^i.$$\endproclaim
\demo{Proof} ($\Rightarrow$) This direction proceeds by induction. It is clear that $0$ has the correct form. If the even element $a$ has the correct form $a= 4\sum_{i=0}^r \epsilon_i3^i$, then the integers
$$\matrix \format\l\\
3a= 4\sum_{i=0}^r \epsilon_i3^{i+1},\\
3a+1=1 +4\sum_{i=0}^r \epsilon_i3^{i+1},\\
3a+4= 4(\sum_{i=0}^r \epsilon_i3^{i+1}+1), \text{ and}\\
3a+5= 1+4(\sum_{i=0}^r \epsilon_i3^{i+1}+1)\endmatrix $$all also have the correct form. \medskip\flushpar ($\Leftarrow$) The recursive definition of $S_3$ shows that the integers $$4\epsilon_r,\quad  4(3\epsilon_r+\epsilon_{r-1}),\quad  4(3^2\epsilon_r+3\epsilon_{r-1}+\epsilon_{r-2}),\quad \dots$$ are all in $S_3$. \qed \enddemo

 The sets $S_p$ are constructed to ensure that certain formulas about the divisibility of binomial coefficients by $p$ hold for $\theta\in S_p$; see Theorem \tref{4}. We introduce the relevant ideas at this point.

\definition{Definition \tnum{Msubp}} If $M$ is an integer and $p$ a prime number, then we define $M_{\#p}$ as follows: 
$$M_{\# p}=\sup\{k\mid \text{ $p^k$ divides $M$}\}.$$ In particular, $0_{\#p}=\infty$. 
\enddefinition

\proclaim{Theorem \tnum{4}} Let $d\in D_p$ for some prime $p\ge 3$.
\roster \item If $p=3$, then 
\smallskip\flushpar {\bf(a)} $\binom{2d}d_{\#3}=\binom{3d}d_{\#3}$,
\smallskip\flushpar {\bf(b)} $\binom{2d+1}d_{\#3}=\binom{3d+2}d_{\#3}$,
\smallskip\flushpar {\bf(c)} $\left(\binom ad\binom{3d-a}d\right)_{\#3}\ge \binom{2d}d_{\#3}$,
\smallskip\flushpar {\bf(d)} $\left(\binom ad\binom{3d-1-a}d\right)_{\#3}\ge \binom{2d}d_{\#3}$, and
\smallskip\flushpar {\bf(e)} $\left(\binom ad\binom{3d+1-a}{d+1}\right)_{\#3}\ge \binom{2d+1}d_{\#3}$,
\smallskip\flushpar for all integers $a$ with $0\le a\le 2d$.

 \item If $p\ge 5$, then 
\smallskip\flushpar {\bf(a)} $\binom{2d}d_{\#p}=\binom{3d}d_{\#p}$,
\smallskip\flushpar {\bf(b)}
$\left(\binom{a}d\binom{3d-a}d\right)_{\#p}\ge \binom{2d}d_{\#p},$
\smallskip\flushpar {\bf(c)} $\left(\binom{a}d\binom{3d-1-a}d\right)_{\#p}\ge \binom{2d}d_{\#p}$, and 
\smallskip\flushpar {\bf(d)} $\left(\binom{a}{d-1}\binom{3d-2-a}d\right)_{\#p}\ge \binom{2d}d_{\#p}$
\smallskip\flushpar for all integers $a$ with $0\le a\le 2d$.
 \endroster
\endproclaim
The proof of Theorem \tref{4} is given in Section 4.

\definition{Definition \tnum{D5.1}} The set $T_0$ is empty and $T_2$ is the set of all positive integers. Assume that $p$ is an odd prime integer. The non-negative integer $a$ is in $T_p$ if there exists an odd integer $J$ and a power $q=p^e$ of $p$, with $e\ge 1$, such that
$$\cases Jq-\frac{q+1}3\le a\le Jq+\frac{q-2}3&\text{if $q\equiv 2\mod 3$}\\
Jq-\frac{q-1}3\le a\le Jq+\frac{q-4}3&\text{if $q\equiv 1\mod 3$}\\
Jq-\frac q3\le a\le Jq+\frac q3 -1&\text{if $q\equiv 0\mod 3$}.
\endcases\tag\tnum{.9}$$\enddefinition

\remark{Remark \tnum{R1.0}} In the language of (\tref{rnd}), the display (\tref{.9}) is equivalent to 
$$\tsize Jq-\{\frac q3\}\le a<Jq+\{\frac q3\}.$$\endremark

\remark{Remark \tnum{R1.1'}} The following statements are immediately clear.
\smallskip\flushpar{\bf(a)} If $a$ is odd, then $a\in T_3\iff |a-Jq|<\frac q3$ for some  odd integer $J$ and  power $q=p^e$ of $p=3$. \smallskip\flushpar{\bf(b)} If $a$ is odd, then $a\in T_3\iff a-1\in T_3$.\endremark 
A similar alternate definition of $T_p$ is available for each $p\ge 5$. 

\proclaim{Observation \tnum{also}}Fix a prime integer $p$ with $p\ge 5$. \roster \item A non-negative even integer $a$ is in $T_p$ if and only if there exists an odd integer $J$ and a power $q=p^e$ of $p$, with $e\ge 1$,  such that $|a-Jq|<\frac q3$.
\item An odd integer $a$ is in $T_p$ if and only if $a+1$ is in $T_p$.\endroster\endproclaim

\demo{Proof}We first prove (1).  Let $a$ and $J$ be integers with $a$ even and $J$ odd. Suppose   that $q\equiv 2\mod 3$. In this case, $Jq-\frac{q+1}3$ is odd; so, 
$$\tsize  Jq-\frac{q+1}3\le a\le Jq+\frac{q-2}3\iff  Jq-\frac{q+1}3+1\le a\le Jq+\frac{q-2}3\iff |Jq-a|\le \frac{q-2}3.$$ However, $|Jq-a|$ and $\frac{q-2}3$ are both integers and there are no   integers in the open interval $(\frac{q-2}3, \frac{q}3)$; so, 
$$\tsize Jq-\frac{q+1}3\le a\le Jq+\frac{q-2}3\iff  |Jq-a|< \frac{q}3.$$ Now suppose that $q\equiv 1\mod 3$. In this case, $Jq-\frac{q-1}3$ is odd; so, 
$$\tsize Jq-\frac{q-1}3\le a\le Jq+\frac{q-4}3\iff Jq-\frac{q-1}3+1\le a\le Jq+\frac{q-4}3\iff |Jq-a|\le \frac{q-4}3.$$ The integers $|Jq-a|$ and  $\frac{q-4}3$ are both odd. There are no odd integers in the open interval $(\frac{q-4}3,\frac{q}3)$; so, 
$$\tsize Jq-\frac{q-1}3\le a\le Jq+\frac{q-4}3\iff |Jq-a|< \frac q3.$$
Assertion (2) follows from the fact that both left hand end points $Jq-\frac{q+1}3$ (if $q\equiv 2$) and $Jq-\frac{q-1}3$  (if $q\equiv 1$) in the definition of $T_p$ are odd integers and both right hand end points are even integers. \qed
\enddemo

\example{Examples \tnum{E1.2}} We see that
$$\split T_3={}&[2,3]\cup [6,11]\cup [14,15]\cup [18,35]\cup[38,39]\cup[42,47]\cup[50,51]\cup[54,107]\\ &{}\cup[110,111]\cup[114,119]\cup [122,123]\cup[126,143]\cup[146,147]\cup\dots,\\ 
T_5={}&[3,6]\cup[13,36]\cup[43,46]\cup[53,56]\cup[63,186]\cup\dots,\ \text{and}\\
T_7={}&[5,8]\cup[19,22]\cup[33,64]\cup[75,78]\cup[89,92]\cup[103,106]\cup\dots\ 
.\endsplit$$
\endexample

\bigskip Now that the sets $S_c$ and $T_c$ have been thoroughly introduced, we turn our attention to other concepts which we use throughout the paper. The $s\times s$ matrix $\varphi=(\varphi_{i,j})$ is {\it alternating} if $\varphi_{j,i}=-\varphi_{i,j}$ and $\varphi_{i,i}=0$, for all $i$ and $j$. The {\it Pfaffian} of $\varphi$ is 
$$\operatorname{Pf}(\varphi)=\cases 0&\text{if $s$ is odd}\\1&\text{if $s=0$}\\
\varphi_{1,2}&\text{if $s=2$}\\
\sum_{j=2}^s(-1)^j\varphi_{1,j}\operatorname{Pf}\left(\matrix \format\l\\\text{$\varphi$ with rows and columns}\\\text{$1$ and $j$ deleted}\endmatrix\right)&\text{if $s\ge 4$ is even.}\endcases$$ If $s$ is odd, then, for each index $\ell$ with $1\le \ell\le s$, we define $\operatorname{Pf}_\ell(\varphi)$: 
$$\operatorname{Pf}_\ell(\varphi) =(-1)^{\ell+1} \operatorname{Pf} (\text{$\varphi$ with row  and column $\ell$ deleted}).\tag\tnum{pfell}$$
It is well-known that the classical adjoint of a square matrix $M$ satisfies $M\operatorname{Adj}M=(\det M)I$.   The corresponding statement for Pfaffians is recorded below.  

\definition{Definition \tnum{D24.1}} If $\varphi$ is an $s\times s$ alternating matrix for some positive even inter $s$, then define $\varphi\check{\phantom{.}}$ to be the alternating $s\times s$ matrix with
$$(\varphi\check{\phantom{.}})_{i,j}=\cases (-1)^{i+j} \operatorname{Pf}(\varphi \text{ with rows and columns $i$ and $j$ deleted})&\text{if $i<j$}\\0&\text{if $i=j$}\\(-1)^{i+j+1} \operatorname{Pf}(\varphi \text{ with rows and columns $i$ and $j$ deleted})&\text{if $j<i$}.\endcases$$\enddefinition
\proclaim{Observation \tnum{pf0}}If $\varphi$ is an $s\times s$ alternating matrix for some positive even integer $s$, then $\varphi\varphi\check{\phantom{.}}=\operatorname{Pf}(\varphi)\cdot I= \varphi\check{\phantom{.}}\varphi$.\endproclaim
\flushpar The proof is straightforward.  

\example{Examples \tnum{E24.1}} 1. If $\varphi=\left[\smallmatrix 0&f\\-f&0\endsmallmatrix\right]$, then $\varphi\check{\phantom{.}}=\left[\smallmatrix 0&-1\\1&0\endsmallmatrix\right]$.

\flushpar 2. If $$\varphi=\bmatrix 
0&\varphi_{1,2}&\varphi_{1,3}&\varphi_{1,4}\\
-\varphi_{1,2}&0&\varphi_{2,3}&\varphi_{2,4}\\
-\varphi_{1,3}&-\varphi_{2,3}&0&\varphi_{3,4}\\
-\varphi_{1,4}&-\varphi_{2,4}&-\varphi_{3,4}&0\endbmatrix$$ is an arbitrary $4\times 4$ alternating matrix, then 
$\varphi\check{\phantom{.}}$ is  the $4\times 4$ alternating matrix
$$\varphi\check{\phantom{.}}=\bmatrix 
0&-\varphi_{3,4}&\varphi_{2,4}&-\varphi_{2,3}\\
\varphi_{3,4}&0&-\varphi_{1,4}&\varphi_{1,3}\\
-\varphi_{2,4}&\varphi_{1,4}&0&-\varphi_{1,2}\\
\varphi_{2,3}&-\varphi_{1,3}&\varphi_{1,2}&0\endbmatrix.$$One easily sees that $\varphi\varphi\check{\phantom{.}}=\operatorname{Pf}(\varphi)I=\varphi\check{\phantom{.}}\varphi$.

\flushpar 3. 
We are particularly interested in the $4\times 4$ alternating matrices $\varphi_{r,n-r}$.
If $r$ and $s$ are non-negative integers, then $\varphi_{r,s}$ is the matrix
$$\varphi_{r,s}=\bmatrix 0&z^{r}&-y^{r}&x^{s} \\
-z^{r}&0&x^{r}&y^{s} \\
y^{r}&-x^{r}&0&z^{s} \\
-x^{s}&-y^{s}&-z^{s}&0 \endbmatrix  \tag\tnum{frn}$$ with entries in $\pmb k[x,y,z]$. We observe that, in the sense of Definition \tref{D24.1},  $\varphi_{r,s}\check{\phantom{.}}$ is equal to $-\varphi_{s,r}$. Thus, the matrices  $\varphi_{r,s}$ and $-\varphi_{s,r}$ form a matrix factorization of   $$(x^{r+s}+y^{r+s}+z^{r+s})=\operatorname{Pf}(\varphi_{r,s})=\operatorname{Pf}(\varphi_{s,r}).$$   In particular, the matrices $\varphi_{r,n-r}$ and $-\varphi_{n-r,r}$ form a matrix factorization of
the defining equation of the diagonal hypersurface ring $R_{\pmb k,n}$; that is, $$\varphi_{r,n}\varphi_{n-r,n}=-(x^n+y^n+z^n)I,\tag\tnum{ffhat}$$where $I$ is the $4\times 4$ identity matrix.\endexample

There is a curious identity which relates the maximal order Pfaffians of an odd-sized   alternating matrix to the Pfaffians of alternating matrices made from its constituent pieces. Remarkably enough, this identity plays a critical role in our calculations. The matrix $\varphi\check{\phantom{.}}$ is defined in Definition \tref{D24.1}. We write $\hat a$ to indicate that the entry $a$ has been deleted. The function $\operatorname{Pf}_\ell$ is defined in (\tref{pfell}). 

\proclaim{Observation \tnum{pf2}} Let $d_2$ be a $(m+3)\times (m+3)$ alternating matrix with $m$ even. Partition $d_2$ into submatrices 
$$d_2=\bmatrix \varphi&\psi^{\text{\rm T}}\\-\psi&\Phi\endbmatrix,$$
where $\varphi$ is a $m\times m$ alternating matrix, $\Phi$ is a $3\times 3$ alternating matrix, and $\psi$ is a $3\times m$ matrix. Then for each index $\ell$, with $1\le \ell\le 3$,
$$\operatorname{Pf}_{m+\ell}(d_2)=\operatorname{Pf}_{\ell}(\psi\varphi\check{\phantom{.}}\psi^{\text{\rm T}})+\operatorname{Pf}(\varphi)\cdot\operatorname{Pf}_{\ell}(\Phi).$$\endproclaim
 \demo{Proof}We prove the result for  $\ell=1$. Expand down column $m+3$ to see that $\operatorname{Pf}_{m+1}(d_2)$ is 
$$\sum_{j=1}^m(-1)^{j+1}\psi_{3,j}\operatorname{Pf}\pmatrix 
0&\dots&\widehat{\varphi_{1,j}}&\dots&\varphi_{1,m}&\psi_{2,1}\\
-\varphi_{1,2}&\dots&\widehat{\varphi_{2,j}}&\dots&\varphi_{2,m}&\psi_{2,2}\\
\vdots & &\vdots&&\vdots&\vdots\\
\widehat{-\varphi_{1,j}}&\dots&\widehat{0}&\dots&\widehat{\varphi_{2,j}}&\widehat{\psi_{2,j}}\\
\vdots & &\vdots&&\vdots&\vdots\\
-\varphi_{1,m}&\dots&\widehat{-\varphi_{j,m}}&\dots&0&\psi_{2,m}\\
-\psi_{2,1}&\dots&\widehat{-\psi_{2,j}}&\dots&-\psi_{2,m}&0\endpmatrix +\Phi_{2,3}\operatorname{Pf}(\varphi).$$
Expand down the column with entries $\{\psi_{2,*}\}$ to see that $\operatorname{Pf}_{m+1}(d_2)$ is $$=\sum_{j=1}^m(-1)^{j+1}\psi_{3,j}\sum_{i=0}^{m}\sigma_{i,j}\psi_{2,i}
\operatorname{Pf}\left(\matrix \format\l\\\text{$\varphi$ with rows and columns}\\\text{$i$ and $j$ deleted}\endmatrix\right)+\Phi_{2,3}\operatorname{Pf}(\varphi),$$where
$$\sigma_{i,j}=\cases (-1)^{i+1}&\text{if $i<j$}\\0&\text{if $i<j$}\\(-1)^{i}&\text{if $j<i$}.\endcases$$
Let $D$ be the $3\times 3$ alternating matrix 
$\psi\varphi\check{\phantom{.}}\psi^{\text{\rm T}}+\operatorname{Pf}(\varphi)\cdot \Phi$. Observe that $\operatorname{Pf}_{m+1}(d_2)$ is equal to the entry of $D$ in row 2 and column 3. On the other hand, $D_{2,3}=\operatorname{Pf}_1(D)$. Thus, $\operatorname{Pf}_{m+1}(d_2)=\operatorname{Pf}_1(D)$, and the proof is complete. \qed\enddemo

\definition{Definition \tnum{num?}} If $Q$ is a noetherian artinian graded $\pmb k$-algebra with unique homogeneous maximal ideal $Q_{+}$, then the socle of $Q$,
$$\operatorname{soc} Q=0\:\!   Q_+=\{q\in Q\mid q  Q_+=0\},$$ is a finite dimensional graded $\pmb k$-vector space:
$\operatorname{soc} Q=\bigoplus\limits_{i=1}^{\ell}\pmb k(-d_i)$. The numbers $d_1\le d_2\le \dots\le d_{\ell}$ are the {\it socle degrees of $Q$}. When recording socle degrees, we use the convention $d_j\!:\!r$ to indicate that $d_j=d_{j+1}=\dots=d_{j+r-1}$.\enddefinition

\remark{Remark \tnum{R?}} If the ring $Q$ of Definition \tref{num?} is the quotient of a polynomial ring $\pmb k[x_1,\dots,x_m]$, then one may read the socle degrees of $Q$ from a minimal homogeneous resolution of $Q$ by free $P$-modules. If 
$$0\to \bigoplus\limits_{i=1}^{\ell} P(-b_i) \to \dots\to P\to Q$$ is such a resolution (with $b_1\le \dots\le b_{\ell}$) and $\deg x_i=1$ for all $i$, then the socle degrees of $Q$ are $$b_1-m\le \dots \le b_{\ell}-m;$$ see, for example \cite{\rref{KV}, Cor.~1.7}.\endremark

If $M$ is an $R$-module and 
$$\Bbb F:\quad \dots @> d_{i+1}>> F_i@> d_{i}>>F_{i-1}  @> d_{i-1}>>\dots @> d_{2}>>F_1 @> d_{1}>>F_0 \to M\to 0$$ is a minimal resolution of $M$, then the $i^{\text{th}}$-syzygy of $M$ is $\operatorname{syz}_iM=\operatorname{im} d_i$. In particular, the truncation $$\Bbb F_{\ge i}:\quad \dots @>d_{i+2}>>F_{i+1}@> d_{i+1}>> F_i\tag\tnum{trun}$$ of $\Bbb F$ is a minimal resolution of $\operatorname{syz}_iM$.

Some of our calculations  involve  binomial coefficients. We recall that $\binom mi$ makes sense for any pair of integers $m$ and $i$; furthermore we  recall the standard properties of
these objects.

\definition{Definition \tnum{7binom}} For integers $i$ and $m$, the binomial
coefficient $\binom{m}{i}$ is defined to be
$$\binom{m}{i}=\left\{ \matrix\format \c&\quad\l\\
\dfrac{m(m-1)\cdots(m-i+1)}{i!}&\text{if $0<i$,}\\\vspace{5pt} 
1&\text{if $0=i$, and}\\\vspace{5pt} 0&\text{if
$i<0$.}\endmatrix \right.$$\enddefinition 

\proclaim{Facts \tnum{7Ob101}} 
\roster
\item"{(a)}" If $i$ and $m$ are integers with $0\le m<i$, then
$\binom{m}{i}=0$. 
\item"{(b)}" If $i$ and $m$ are integers, then 
$\binom{m}{i-1}+\binom{m}{i}=\binom{m+1}{i}$.
\item"{(c)}" If $i$ and $m$ are integers with $0\le m$, then $\binom{m}{i}=\binom{m}{m-i}.$
\item"{(d)}"If $i$ and $m$ are integers, then $(m-i)\binom mi=\binom m{i+1}(i+1)$.
\item"{(e)}"If $i$ and $m$ are integers, then $i\binom{m}{i}=\binom{m-1}{i-1}m.$
\item"{(f)}"If $0\le i\le m$ are integers, then $\binom mi=\frac {m!}{i!(m-i)!}$.
\endroster
\endproclaim

\SectionNumber=2\tNumber=1
\heading Section \number\SectionNumber. \quad The technique for resolving $Q_{\pmb k,n,N}$.\endheading

Lemma \tref{setup} describes our main technique for resolving $Q_{\pmb k,n,N}$ as a module over $R_{\pmb k,n}$. The Lemma is set in a slightly more general context. To apply the Lemma in the context of $Q_{\pmb k,n,N}$, one takes $P=\pmb k[x,y,z]$, $m=4$,  $$ \matrix\format\l&\quad\l\\    \pmb x=\bmatrix x^N&y^N&z^N\endbmatrix,&f=x^n+y^n+z^n,\\\vspace{5pt} X=\bmatrix 
0&z^N&-y^N\\
-z^N&0&x^N\\
y^N&-x^N&0\endbmatrix,\quad\text{and}& \varphi=\cases \varphi_{r,n-r}&\text{if $\theta$ is odd, or}\\
\varphi_{n-r,r}&\text{if $\theta$ is even,}\endcases\endmatrix \tag\tnum{ssetup}$$ where we have retained the notation $N=\theta n+r$, with $0\le r\le n-1$. Recall from (\tref{ffhat}) that
$\operatorname{Pf}(\varphi_{r,n-r})=\operatorname{Pf}(\varphi_{n-r,r})=f$ and $\varphi_{r,n-r}\varphi_{n-r,r}=-fI$. To apply Lemma \tref{setup}, one must find a matrix $\psi$ and a unit $u$ so that
$$\text{the entries of  $\psi\varphi\check{\phantom{.}}\psi^{\text{\rm T}}-uX$ are in the ideal $(f)P$.}\tag\tnum{hyp20}$$ (The matrix $\varphi\check{\phantom{.}}$ is defined in Definition \tref{D24.1}.) As soon as (\tref{hyp20}) is accomplished, then the $R_{\pmb k,n}$-resolution of $Q_{\pmb k,n,N}$ is given in conclusion (3) of Lemma \tref{setup}. Conclusions (1) and (2) are steps along the way to conclusion (3); these steps are of interest in their own right. Conclusion (1) lists the generators of $(x^N,y^N,z^N)\:\! f$ and conclusion (2) gives the $P$-resolution of $Q_{\pmb k,n,N}$. We can read the socle degrees of $Q_{\pmb k,n,N}$ and, indeed, the entire Hilbert function of $Q_{\pmb k,n,N}$, from the resolution of conclusion (2). On the other hand, it is difficult to find $\psi$ and $u$ so that (\tref{hyp20}) is satisfied. We show our solutions $\psi$ and $u$ to the problem posed in (\tref{hyp20}) in the proof of Theorem \tref{I-MAIN}. We prove that our solutions work in Section 7. 
\proclaim{Lemma \tnum{setup}}Let $P$ be a commutative noetherian ring, $\pmb x=[x_1,x_2,x_3]$ be a ${1\times 3}$ matrix whose entries generate a perfect grade three ideal in $P$, $m$ be a positive even integer, $\varphi$ be an $m\times m$ alternating matrix with entries in $P$, $\psi$ be a $3\times m$ matrix with entries in $P$, and $u$ be a unit of $P$. Define $f$ to be the Pfaffian of $\varphi$ and $X$ to be the $3\times 3$ alternating matrix with $\operatorname{Pf}_i(X)=x_i$ for $1\le i\le 3$. Assume that the entries of $\psi\varphi\check{\phantom{.}}\psi^{\text{\rm T}}-uX$ are in the ideal $(f)P$. Define $\Phi$ to be a  $3\times 3$ alternating matrix with
$$\psi \varphi\check{\phantom{.}}\psi^{\text{\rm T}}+f\Phi=uX,\tag\tnum{crit}$$  and define
$$d_2=\bmatrix \varphi&\psi^{\text{\rm T}}\\-\psi&\Phi\endbmatrix.\tag\tnum{d2}$$
Let $R$ be the ring $P/(f)$ and $Q$ be the ring $P/(f,I_1(\pmb x))$. Then the following statements hold.
\roster
\item The ideal $I_1(\pmb x):f$ of $P$ is generated by the maximal order Pfaffians of $d_2$.
\item The maps and modules 
$$\Bbb F:\quad 0\to F_3@> f_3>>F_2@> f_2>> F_1@> f_1>> F_0,$$
with
$$ F_3= P^m,\quad F_2= \matrix P^m\\\oplus\\P^3\endmatrix, \quad F_1= \matrix P^3\\\oplus\\P\endmatrix, \quad F_0=P,$$
$$f_3=\bmatrix \varphi\\-u\psi\endbmatrix,\quad f_2=\bmatrix u\psi\varphi\check{\phantom{.}}&fI\\-\pmb b&-\pmb x\endbmatrix,\quad f_1=\bmatrix \pmb x&f\endbmatrix$$ form a resolution of $Q$ by free $P$-modules, where $$\pmb b=\bmatrix \operatorname{Pf}_1(d_2)&\dots&\operatorname{Pf}_m(d_2)\endbmatrix.$$
\item  The maps and modules
$$\dots @>>> R^m@>\varphi>> R^m @> \varphi\check{\phantom{.}} >>  R^m@>\varphi>> R^m @> \psi\varphi\check{\phantom{.}} >>R^3@> \pmb x>> R$$
form a resolution of $Q$ by free $R$-modules.\endroster\endproclaim 

\demo{Proof}The ideal $I_1(\pmb x)$ is generated by a regular sequence; hence, the Koszul complex: 
$$0\to P@> \pmb x^{\text{\rm T}} >>P^3@>\pmb  X>> P^3@> \pmb x >>P\tag\tnum{7p/x}$$ is a resolution  of $P/I_1(\pmb x)$ by free $P$-modules.
Let $J$ be the ideal of $P$ generated by the maximal order Pfaffians of $d_2$. We will show that $J$ is a grade three Gorenstein ideal and we will compute that $\frac{I_1(\pmb x)\:\! J}{I_1(\pmb x)}$ is a cyclic module generated by the class of $f$. The theory of linkage then shows that $I_1(\pmb x)\:\! f=J$. 

In Observation \tref{pf2} we calculated that the Pfaffians
$$\bmatrix \operatorname{Pf}_{m+1}(d_2)&\operatorname{Pf}_{m+2}(d_2)&\operatorname{Pf}_{m+3}(d_2)\endbmatrix $$ of $d_2$ are precisely the same as the Pfaffians 
$$\bmatrix \operatorname{Pf}_1(\psi \varphi\check{\phantom{.}}\psi^{\text{\rm T}}+(\operatorname{Pf} \varphi)\Phi)&\operatorname{Pf}_2(\psi \varphi\check{\phantom{.}}\psi^{\text{\rm T}}+(\operatorname{Pf} \varphi)\Phi)&\operatorname{Pf}_3(\psi \varphi\check{\phantom{.}}\psi^{\text{\rm T}}+(\operatorname{Pf} \varphi)\Phi)\endbmatrix $$ of the alternating matrix $\psi \varphi\check{\phantom{.}}\psi^{\text{\rm T}}+(\operatorname{Pf} \varphi)\Phi$.
The hypothesis (\tref{crit}) ensures that $\psi \varphi\check{\phantom{.}}\psi^{\text{\rm T}}+(\operatorname{Pf} \varphi)\Phi=uX$; hence, 
$$\bmatrix \operatorname{Pf}_{m+1}(d_2)&\operatorname{Pf}_{m+2}(d_2)&\operatorname{Pf}_{m+3}(d_2)\endbmatrix =u \pmb x$$ and therefore $J$, which is the ideal of maximal order Pfaffians of $d_2$, has grade three and is a grade three Gorenstein ideal; see, for example \cite{\rref{BE}, Thm.~2.1}. Let $\pmb b$ be the matrix defined in the statement of (2).  The resolution of $P/J$ by free $P$-modules is 
 $$0\to P@> d_3 >> P^{m+3}@> d_2 >> P^{m+3} @> d_1 >> P,\tag\tnum{7p/j}$$with $$d_3={\bmatrix \pmb b^{\text{\rm T}}\\ u\pmb x^{\text{\rm T}}\endbmatrix},\quad \text{$d_2$ given in (\tref{d2})},\quad\text{and}\quad d_1={\bmatrix \pmb b& u\pmb x\endbmatrix} .$$ 

Let $\pi:P/I_1(\pmb x)\to P/J$ be the natural quotient map.  We claim that $$\CD
0@>>> P@> u\pmb x^{\text{\rm T}}>> P^3 @> u\pmb X>> P^3 @> u\pmb x >>P@>>>P/I_1(\pmb x)\\
@. @V\alpha_3 VV @VV \alpha_2V @VV \alpha_1 V @V = VV @V\pi VV\\
0@>>> P@> d_3>> P^{m+3} @>  d_2>> P^{m+3} @> d_1>> P@>>> P/J,\endCD\tag\tnum{CD}$$
with
$$\alpha_3=f ,\quad \alpha_2= {\bmatrix -\varphi\check{\phantom{.}}\psi^{\text{\rm T}}\\ fI\endbmatrix} , \quad \alpha_1= {\bmatrix 0_{4\times 3} \\I_{3\times 3}\endbmatrix},$$
is a map of complexes from a resolution of $P/I_1(\pmb x)$ isomorphic to  (\tref{7p/x}) to the resolution (\tref{7p/j}) of $P/J$. To see that $d_2\alpha_2=\alpha_1 X$, one uses the fact that $\varphi\varphi\check{\phantom{.}}=fI$ (see Observation \tref{pf0}) and the hypothesis (\tref{crit}). To see that $d_3\alpha_3=\alpha_2 (u\pmb x^{\text{\rm T}})$, observe first that $$d_2\alpha_2 (u\pmb x^{\text{\rm T}})=\alpha_1 (u\pmb X)(u\pmb x^{\text{\rm T}})=0.$$ Thus, $\alpha_2 (u\pmb x^{\text{\rm T}})$ is in the kernel of $d_2$. The complex (\tref{7p/j}) is exact; so,  $\alpha_2 (u\pmb x^{\text{\rm T}})=v d_3$    for some $v$ in $P$. On the other hand,
 $$ \alpha_2 (u\pmb x^{\text{\rm T}})=   {\bmatrix -\varphi\check{\phantom{.}}\psi^{\text{\rm T}}(u\pmb x^{\text{\rm T}})\\ fu\pmb x^{\text{\rm T}}\endbmatrix} \quad\text{and}\quad v d_3= v{\bmatrix \pmb b^{\text{\rm T}}\\ u\pmb x^{\text{\rm T}}\endbmatrix}.$$ It is obvious that $v$ must be $f$ and  (\tref{CD}) is indeed a map of complexes. 

Standard results from homological algebra now show that $I_1(\pmb x)\:\! f$ is equal to $J$.
Indeed, the long exact sequence of homology which is obtained by applying $\operatorname{Hom}_P(\underline{\phantom{X}},P)$ to the short exact sequence
$$0\to J/I_1(\pmb x)\to P/I_1(\pmb x)@> \pi>> P/J\to 0,$$ yields that 
$$\pi^*\: \operatorname{Ext}^3_P(P/J,P)\to \operatorname{Ext}^3_P(P/I_1(\pmb x),P)\tag\tnum{pi*}$$ is an injection. On the other hand, one may apply $\operatorname{Hom}_P(\underline{\phantom{X}},P)$ to the comparison map of resolutions which is given in (\tref{CD}) to see  that the injection (\tref{pi*})  is also equal to $$P/J@> f >> P/I_1(\pmb x).   $$In other words, $I_1(\pmb x)\:\! f$ is equal to $J$. This establishes (1). Take the mapping cone of the dual of (\tref{CD}) to establish (2). 

We prove (3). Let $\overline{\phantom{X}}$ represent the functor $\underline{\phantom{X}}\otimes_P R$. Take the resolution $\Bbb F$ of $Q$ from (2). We see that $\overline {\Bbb F}$ is a complex of free $R$-modules with
$$\operatorname{H}(\overline {\Bbb F})=\operatorname{Tor}_i^P(R,Q)=\cases Q&\text{if $i=0$ or $1$}\\0&\text{otherwise}.\endcases$$ Furthermore, the cycle 
$$\xi= \bmatrix 0\\0\\0\\1\endbmatrix\tag\tnum{xi}$$ in $\overline{F}_1$ represents a generator of $\operatorname{H}_1(\overline{\Bbb F})$. We kill the homology in $\overline{\Bbb F}$. Define $R$-module homomorphisms $\beta_i\: \overline{F}_i\to \overline{F}_{i+1}$ by 
$$\beta_2=\bmatrix \varphi\check{\phantom{.}}&0\endbmatrix,\quad \beta_1=\bmatrix 0&0\\-I&0\endbmatrix,\quad \beta_0=\bmatrix 0\\1\endbmatrix.$$ A straightforward calculation shows that
$$\CD 0@>>> \overline{F}_3 @>\overline{f}_3>> \overline{F}_2 @> \overline{f}_2 >> \overline{F}_1
@> \overline{f}_1 >> \overline{F}_0\\ @. @. @V \beta_2 VV @V \beta_1 VV @V \beta_0 VV\\
@. 0@>>> \overline{F}_3 @>\overline{f}_3>> \overline{F}_2 @> \overline{f}_2 >> \overline{F}_1
@> \overline{f}_1 >> \overline{F}_0\endCD\tag\tnum{ll}$$is a map of complexes. It is clear that $\beta_0$ induces an isomorphism from $\operatorname{H}_0$ of the top line of (\tref{ll}) to $\operatorname{H}_1$ of the bottom  line of (\tref{ll}). Let $\Bbb M$ be the total complex of (\tref{ll}). We have shown that the homology of $\Bbb M$ is concentrated in positions $0$ and $3$ and the $\xi$ from (\tref{xi}) of the summand $\overline{F}_1$ in $\Bbb M_3=\overline{F}_1\oplus \overline{F}_3$ represents a generator of $\operatorname{H}_3(\Bbb M)$. Iterate this process to see that the $Q$ is resolved by the total complex $\Bbb T$ of the following infinite double complex:
$$\CD 
@.  \vdots  @. \vdots   @. \vdots  \\
@. @V \beta_2  VV @V  \beta_1 VV @V \beta_0  VV\\
0@>>> \overline{F}_3 @>\overline{f}_3>> \overline{F}_2 @> \overline{f}_2 >> \overline{F}_1
@> \overline{f}_1 >> \overline{F}_0\\ @. @. @V \beta_2 VV @V \beta_1 VV @V \beta_0 VV\\
@.0@>>> \overline{F}_3 @>\overline{f}_3>> \overline{F}_2 @> \overline{f}_2 >> \overline{F}_1
@> \overline{f}_1 >> \overline{F}_0\\@. @. @. @V \beta_2 VV @V \beta_1 VV @V \beta_0 VV\\
@.@. 0@>>> \overline{F}_3 @>\overline{f}_3>> \overline{F}_2 @> \overline{f}_2 >> \overline{F}_1
@> \overline{f}_1 >> \overline{F}_0.\endCD$$The complex of (3) is a summand of $\Bbb T$. \qed
\enddemo

\proclaim{Corollary \tnum{5gen}} If the ring $Q_{\pmb k,n,N}$ of {\rm (\tref{DHSR.1})} has finite projective dimension over $R_{\pmb k,n}$, then the ideal $(x^N,y^N,z^N)\:\! (x^n+y^n+z^n)$ in $P=\pmb k[x,y,z]$ can be generated by $5$ elements. \endproclaim

\demo{Proof} Let $R=R_{\pmb k,n}$ and $Q=Q_{\pmb k,n,N}$. Assume $\operatorname{pd}_RQ<\infty$. The rings $R$  and $ Q$ have depth $2$ and $0$, respectively. The Auslander-Buchsbaum Theorem guarantees that $\operatorname{pd}_RQ=2$ and therefore, the Hilbert-Burch Theorem ensures that there exists  a $3\times 2$ matrix $\bar \psi$, over $R$, whose signed maximal order minors are $x^N,y^N,z^N$. (See (\tref{smom}), if necessary.) Lift $\bar \psi$ back to a matrix $\psi$ in $P$. Let $\varphi$ be the matrix $\left[\smallmatrix 0&f\\-f&0\endsmallmatrix\right]$ over $P$, where $f=x^n+y^n+z^n$. Recall from Example \tref{E24.1} that $\varphi\check{\phantom{.}}=\left[\smallmatrix 0&-1\\1&0\endsmallmatrix\right]$. Observe that the non-zero entries of $\psi \varphi\check{\phantom{.}}\psi^{\text{\rm T}}$ are the   maximal order minors of $\psi$, up to sign. Indeed, if $u=-1$ and
$$X=\bmatrix 0&z^N&-y^N\\-z^N&0&x^N\\y^N&-x^N&0\endbmatrix,$$ then entries of of the alternating matrix $\psi \varphi\check{\phantom{.}}\psi^{\text{\rm T}}-uX$ are in the ideal $(f)P$. The critical hypothesis (\tref{hyp20}) of Lemma \tref{setup} is satisfied; hence, the ideal $(x^N,y^N,z^N)\:\! f$ of $P$ if generated by the maximal order Pfaffians of the $5\times 5$ alternating matrix which is produced in Lemma \tref{setup}. \qed \enddemo

\remark{\bf Remark \tnum{2.14}} The ideal $(x^N,y^N,z^N)\:\! f$ was introduced in (\tref{GknN}). It is shown in Theorem \tref{???} that $\operatorname{pd}_RQ=\infty$ if and only if the hypotheses of Theorem \tref{I-MAIN} hold. When the hypotheses of Theorem \tref{I-MAIN} hold, then the ideal (\tref{GknN}) is $7$-generated. When Theorems \tref{???} and \tref{I-MAIN} are combined with Corollary \tref{5gen}, we see that 
$\operatorname{pd}_RQ$ is infinite  if and only if the ideal of (\tref{GknN}) is $7$-generated.\endremark

\SectionNumber=3\tNumber=1
\heading Section \number\SectionNumber. \quad The resolutions of $Q_{\pmb k,n,N}$ when $\lfloor\frac Nn\rfloor$ is in $S_c$.  \endheading

The following notation and  hypotheses are  in effect throughout this section. Recall the set $S_c$ from Definition \tref{!D18.1}. \definition {Data \tnum{Data2}} Let $\pmb k$ be a field of characteristic $c$, $n$ and $N$ be positive integers with $n$ not dividing $N$. Write $N=\theta n+r$ for integers $\theta$ and $r$ with $0< r<n$. Let $P$ be the polynomial ring $P=\pmb k[x,y,z]$, $R$ be the diagonal hypersurface ring $R_{\pmb k,n}$ of {\rm(\tref{DHSR})}, and $Q$ be the quotient ring $Q_{\pmb k,n,N}$ of {\rm(\tref{DHSR.1})}. {\bf Assume that $\theta$ is in the set $S_c$.} \enddefinition

In Theorem \tref{I-MAIN} we give the minimal homogeneous resolution of $Q$ by free $P$-modules and also the minimal homogeneous resolution of $Q$ by free $R$-modules.  
Our proof  is based on the method of Lemma \tref{setup}. One must find a matrix $\psi$ and a unit $u$ so that (\tref{hyp20}) is satisfied. We use the polynomials $\operatorname{Poly}_{d,a,b}$ of Definition \tref{D21.1} to build our candidate for $(\psi,u)$. The key calculation which shows that our candidate for $(\psi,u)$ satisfies (\tref{hyp20}) is called Lemma \tref{NandI}. The proof of Lemma \tref{NandI} is given in Section 7.

\definition{Definition \tnum{D21.1}} Given non-negative integers $d$, $a$, $b$, we define
$\operatorname{Poly}_{d,a,b}(A,B)$ to be the following polynomial in the polynomial ring $\Bbb Z[A,B]$:
$$\operatorname{Poly}_{d,a,b}(A,B)=\sum\limits_{i=0}^{d}(-1)^i \binom{a+d-i}{a}\binom{b+i}{b}A^{d-i}B^i.$$\enddefinition

\remark{Remark}We see that $\operatorname{Poly}_{d,b,a}(B,A)=(-1)^d\operatorname{Poly}_{d,a,b}(A,B)$.\endremark

\definition{Definition \tnum{D21.1.5}}For each positive integer $\delta$ define the polynomials
$\frak P_{2\delta-1}$ and $\frak P_{2\delta}$ in $\Bbb Z[A,B,C]$ by
$$\frak P_{2\delta-1}(A,B,C)=\left\{\matrix \format\l\\
(-1)^{\delta}A\operatorname{Poly}_{\delta-1,\delta,\delta-1}(A,B)\operatorname{Poly}_{\delta-1,\delta,\delta-1}(A,C)\\
+B\operatorname{Poly}_{\delta-1,\delta,\delta-1}(B,A)\operatorname{Poly}_{\delta-1,\delta,\delta-1}(A,C)\\
+C\operatorname{Poly}_{\delta-1,\delta,\delta-1}(A,B)\operatorname{Poly}_{\delta-1,\delta,\delta-1}(C,A)\endmatrix\right. $$
and 
$$\frak P_{2\delta}(A,B,C)=\left\{\matrix \format\l\\
(-1)^{\delta+1}\operatorname{Poly}_{\delta,\delta,\delta}(A,B)\operatorname{Poly}_{\delta,\delta,\delta}(A,C)\\
+B\operatorname{Poly}_{\delta-1,\delta,\delta}(B,A)\operatorname{Poly}_{\delta,\delta,\delta}(A,C)\\
+C\operatorname{Poly}_{\delta,\delta,\delta}(A,B)\operatorname{Poly}_{\delta-1,\delta,\delta}(C,A).\\
\endmatrix \right.$$ 
\enddefinition
\remark{Remark} We see that $$\frak P_{2\delta-1}(A,B,C)=\frak P_{2\delta-1}(A,C,B)\quad\text{and}\quad \frak P_{2\delta}(A,B,C)=\frak P_{2\delta}(A,C,B).$$ \endremark 
\proclaim{Lemma \tnum{NandI}}For each positive integer $\delta$, the polynomials
$$ P_{2\delta-1}(A,B,C)=\frak P_{2\delta-1}(A,B,C)
+ (-1)^{\delta+1}\binom{2\delta}{\delta}\binom{3\delta-1}{\delta-1}A^{2\delta-1}
$$
and $$P_{2\delta}(A,B,C)=\frak P_{2\delta}(A,B,C)
+(-1)^{\delta}\binom{2\delta}\delta\binom{3\delta}\delta A^{2\delta}$$ are in the ideal $(A+B+C)\Bbb Z[A,B,C]$.\endproclaim

\proclaim{Theorem \tnum{I-MAIN}} Assume that   the notation and hypotheses of Data {\rm\tref{Data2}} are in effect. Then the following statements hold. 
\roster\item If $\theta=2\delta-1$ is odd, then the minimal homogeneous resolution of $Q$ by free $P$-modules has the form 
$$0\to \matrix P(-3\delta n+n-2r)^3\\\oplus\\P(-3\delta n)\endmatrix \to \matrix P(-3\delta n+n-r)^3\\\oplus\\
P(-3\delta n+2n-3r)\\\oplus\\ P(-2\delta n-r)^3\endmatrix@>>>\matrix P(-N)^3\\\oplus\\P(-n)\endmatrix @>>> P $$
and the minimal homogeneous resolution of $Q$ by free $R$-modules has the form
$$\multline \shoveleft\cdots @> \varphi_{r,n-r}>> \matrix R(-3\delta n-r)^3\\\oplus \\ R(-3\delta n+n-3r)\endmatrix @> \varphi_{n-r,r}>>
\matrix R(-3\delta n +n- 2r)^3\\\oplus \\ R(-3\delta n)\endmatrix \hfill\\\vspace{8pt}\hfill \shoveright {  @> \varphi_{r,n-r}>>\matrix R(-3\delta n +n- r)^3\\\oplus \\ R(-3\delta  n+2n-3r)\endmatrix @>>> R(-N)^3\to R. }\endmultline
 $$In particular, $\operatorname{pd}_RQ=\infty$ and the socle degrees of $Q$ are  
$$(3\delta n-n+2r-3)\!:\!3,\quad 3\delta n-3.$$  
\item If $\theta=2\delta$ is even, then the minimal homogeneous resolution of $Q$ by free $P$-modules has the form 
$$0\to \matrix P(-3\delta n-n-r)^3\\\oplus\\P(-3\delta n-3r)\endmatrix\to \matrix P(-3\delta n-2r)^3\\\oplus\\ P(-3\delta n-n)\\\oplus\\P(-2\delta n-n-r)^3\endmatrix \to \matrix P(-N)^3\\\oplus\\P(-n)\endmatrix \to P $$
and the minimal  homogeneous resolution of $Q$ by free $R$-modules has the form
$$\multline \shoveleft \dots @>_{n-r,r}>> \matrix R(-3\delta n-n-2r)^3\\\oplus\\R(-3\delta n-2n)\endmatrix @> \varphi_{r,n-r}>> \matrix R(-3\delta n-n-r)^3\\\oplus\\R(-3\delta n-3r)\endmatrix \hfill\\\vspace{8pt}\hfill \shoveright {@> \varphi_{n-r,r}>> \matrix R(-3\delta n-2r)^3\\\oplus\\R(-3\delta n-n) \endmatrix \to R(-N)^3\to  R.}\endmultline$$
In particular, $\operatorname{pd}_RQ=\infty$ and the socle degrees of $Q$ are $$(3\delta n+n+r-3)\!:\!3,\quad 3\delta n+3r-3.$$
\endroster \endproclaim

\remark{Remark} All of the differentials in the resolutions of Theorem \tref{I-MAIN} are explicitly described in Lemma \tref{setup}.\endremark

Before we prove Theorem \tref{I-MAIN}, we observe that this Theorem contains one direction of Theorem \tref{???}.
\proclaim{Corollary \tnum{COR}} Let $\pmb k$ be a field of characteristic $c$ and let $n$ and $N$ be positive integers. If  $N=\theta n+r$, with $\theta\in S_c$ and $1\le r\le n-1$, then $\operatorname{pd}_RQ=\infty$. \endproclaim

We also offer the following reformulation of part of Theorem \tref{I-MAIN}
\proclaim{Corollary \tnum{I+}} If the hypotheses of Data \tref{Data2} are in effect, then there exist an exact sequence of $R$-modules{\rm:}
$$ 0\to M_{\pmb k,n,r}(-3\delta n+2n-2r)\to R(-N)^3\to R\to Q\to 0,\quad\text{if $N=(2\delta-1)n+r$},$$ or
$$0\to M_{\pmb k,n,n-r}(-3\delta n-r)\to R(-N)^3\to R\to Q\to 0,\quad\text{if $N=2\delta n+r$},$$where the $R_{\pmb k,n}$-module
 $M_{\pmb k,n,r}$ is defined to be the following cokernel{\rm:}
$$\matrix R(-n)^3\\\oplus\\R(2r-2n)\endmatrix @> \varphi_{r,n-r}>> 
\matrix R(r-n)^3\\\oplus\\R(-r)\endmatrix@> >>M_{\pmb k,n,r}\to 0,$$ and $\varphi_{r,n-r}$ is given in {\rm(\tref{frn})}. 
\endproclaim

\demo{Proof of Theorem {\rm\tref{I-MAIN}}} We apply Lemma \tref{setup} with $\pmb x$, $X$, $f$, and $\varphi$ as given in (\tref{ssetup}) and  we consider two cases. In the first case $\operatorname{char} \pmb k\neq 3$ and in the second case $\operatorname{char} \pmb k=3$.  

We begin with   $\operatorname{char} \pmb k\neq 3$.
Define the integer $\gamma$  by 
$$\gamma=\cases 1&\text{if $\operatorname{char} \pmb k$ equals $0$ or $2$}\\
p^{\ell}&\text{if $\operatorname{char} \pmb k=p\ge 5$ and $\binom{2 \delta}\delta_{\#p}=\ell$}.\endcases$$
We know from the hypothesis  that $\theta\in S_c$. We next show that  $$\tsize \text{$\binom{2\delta}{\delta}\binom{3\delta}{\delta}/\gamma^2$ is a unit in $\pmb k$.}\tag\tnum{unit}$$ There is nothing to show if  $\operatorname{char} \pmb k=0$. If $\operatorname{char} \pmb k=2$, then $\theta\in S_2=\{0\}$; thus, $\delta=0$ and again there is nothing to show. If $\operatorname{char} \pmb k=p\ge 5$, then Remark \tref{!R18.1} ensures that $\delta\in D_p$ and Theorem \tref{4} yields that $\binom{3\delta}{\delta}_{\#p}=\ell$; and therefore (\tref{unit}) holds. We next observe that 
$\gamma$ divides  every coefficient of each of the polynomials
$$\operatorname{Poly}_{\delta,\delta,\delta}(A,B),\quad \operatorname{Poly}_{\delta-1,\delta,\delta}(A,B),\quad\text{and}\quad \operatorname{Poly}_{\delta-1,\delta,\delta-1}(A,B).$$(Again, the assertion is obvious if $\operatorname{char} \pmb k$ is equal to $0$ or $2$ and the assertion follows from Theorem \tref{4} if $\operatorname{char} \pmb k=p\ge 5$.)  Let
$$\matrix \format \r &\ \l\ &\l \quad &\r &\ \l\ &\l\\
\operatorname{Poly}'_{\delta-1,\delta,\delta-1}&=&\operatorname{Poly}_{\delta-1,\delta,\delta-1}/\gamma& \operatorname{Poly}'_{\delta,\delta,\delta}&=&\operatorname{Poly}_{\delta,\delta,\delta}/\gamma,\\ 
\operatorname{Poly}'_{\delta-1,\delta,\delta}&=&\operatorname{Poly}_{\delta-1,\delta,\delta}/\gamma,&\frak P_{2\delta-1}'&=&\frak P_{2\delta-1}/\gamma^2,\ \text{and}\\  \frak P_{2\delta}'&=&\frak P_{2\delta}/\gamma^2\endmatrix $$ in $\Bbb Z[A,B,C]$ and 
$$u_{2\delta-1}=  (-1)^{\delta+1}\binom{2\delta}{\delta}\binom{3\delta-1}{\delta-1}/\gamma^2\quad\text{and}\quad
u_{2\delta}=(-1)^{\delta}\binom{2\delta}{\delta}\binom{3\delta}{\delta}/\gamma^2.$$
in $\Bbb Z$. Recall that $\binom{3\delta-1}{\delta-1}=\frac 13 \binom{3\delta}{\delta}$ and $3$ is a unit in $\pmb k$; so
$u_{2\delta-1}$ and $u_{2\delta}$ both are units in $\pmb k$. 
We have $$ \matrix\format\l&\l\\ \gamma^2(\frak P'_{2\delta-1}(A,B,C)+u_{2\delta-1}A^{2\delta-1})&{}={}\frak P_{2\delta-1}(A,B,C)+(-1)^{\delta+1}\binom{2\delta}{\delta}\binom{3\delta-1}{\delta-1}A^{2\delta-1}\\&{}={}P_{2\delta-1}(A,B,C)\ \text{and}\endmatrix$$
$$\matrix  \gamma^2(\frak P'_{2\delta}(A,B,C)+u_{2\delta}A^{2\delta})= \frak P_{2\delta}(A,B,C)+(-1)^{\delta}\binom{2\delta}{\delta}\binom{3\delta}{\delta}A^{2\delta}=P_{2\delta}(A,B,C).\endmatrix $$Therefore, according to Lemma \tref{NandI},
$$\gamma^2(\frak P'_{2\delta-1}(A,B,C)+u_{2\delta-1}A^{2\delta-1}) \quad\text{and}\quad \gamma^2(\frak P'_{2\delta}(A,B,C)+uA^{2\delta})$$ are both in the prime ideal $(A+B+C)\Bbb Z[A,B,C]$  of $\Bbb Z[A,B,C]$. It follows that 
$$ \frak P'_{2\delta-1}(A,B,C)+u_{2\delta-1}A^{2\delta-1}  \quad\text{and}\quad  \frak P'_{2\delta}(A,B,C)+u_{2\delta}A^{2\delta} $$ 
are in $(A+B+C)\Bbb Z[A,B,C]$. Define
$\Cal R_{2\delta-1}(A,B,C)$ and $\Cal R_{2\delta}(A,B,C)$ to be the polynomials in $\Bbb Z[A,B,C]$ which satisfy
$$\matrix\format\l&\c&\l\\ \frak P'_{2\delta-1}(A,B,C)+ u_{2\delta-1}A^{2\delta-1}&{}={}&(A+B+C)\Cal R_{2\delta-1}(A,B,C)
\ \text{and}\\ \frak P'_{2\delta}(A,B,C)+ u_{2\delta}A^{2\delta}&{}={}&(A+B+C)\Cal R_{2\delta}(A,B,C)
.\endmatrix  \tag\tnum{zzz}$$

We now focus on assertion (1); that is, we have  $\theta=2\delta-1$ is odd. Recall, from (\tref{ssetup}), that  $\varphi$ is   the matrix $\varphi_{r,n-r}$ of {\rm(\tref{frn})}. It follows that $\varphi\check{\phantom{.}}$ is the matrix $-\varphi_{n-r,r}$.  We see from (\tref{ffhat}) that $\varphi\varphi\check{\phantom{.}}=\varphi\check{\phantom{.}}\varphi=fI$. 
Define $\psi$ to be the matrix $$\eightpoint \left[\smallmatrix 
0&(-1)^{\delta-1}z^r\operatorname{Poly}'_{\delta-1,\delta,\delta-1}(z^n,y^n)&y^r\operatorname{Poly}'_{\delta-1,\delta,\delta-1}(y^n,z^n)&0  \\
z^r\operatorname{Poly}'_{\delta-1,\delta,\delta-1}(z^n,x^n)&0&(-1)^{\delta-1}x^r\operatorname{Poly}'_{\delta-1,\delta,\delta-1}(x^n,z^n)&0 \\
(-1)^{\delta-1}y^r\operatorname{Poly}'_{\delta-1,\delta,\delta-1}(y^n,x^n)&x^r\operatorname{Poly}'_{\delta-1,\delta,\delta-1}(x^n,y^n)&0&0\endsmallmatrix\right] $$ and $$\Phi= 
\bmatrix 
0& z^r \Cal R_{2\delta-1}(z^n,x^n,y^n)&-y^r \Cal R_{2\delta-1}(y^n,x^n,z^n)\\
-z^r\Cal R_{2\delta-1}(z^n,x^n,y^n)&0& x^r\Cal R_{2\delta-1}(x^n,y^n,z^n)\\
 y^r \Cal R_{2\delta-1}(y^n,x^n,z^n)&-x^r\Cal R_{2\delta-1}(x^n,y^n,z^n)&0\endbmatrix.$$
One may easily calculate that $\psi\varphi\check{\phantom{.}}\psi^{\text{\rm T}}$  is equal to 
$$\bmatrix 0&-z^r\frak P'_{2\delta-1}(z^n,x^n,y^n)&y^r\frak P'_{2\delta-1}(y^n,x^n,z^n)\\z^r\frak P'_{2\delta-1}(z^n,x^n,y^n)&0&-x^r\frak P'_{2\delta-1}(x^n,y^n,z^n)\\-y^r\frak P'_{2\delta-1}(y^n,x^n,z^n)&x^r\frak P'_{2\delta-1}(x^n,y^n,z^n)&0\endbmatrix.$$
Therefore, the definition (\tref{zzz}) of $\Cal R_{2\delta-1}$ tells us that 
$$\psi\varphi\check{\phantom{.}}\psi^{\text{\rm T}}+f\Phi= u_{2\delta-1}\pmb X.$$ 
The condition (\tref{hyp20}) is satisfied. Lemma \tref{setup} applies. All of the entries in all of the matrices are homogeneous of positive degree. Thus, an explicit minimal generating set for   the ideal $(x^N,y^N,z^N)\:\! f$ of $P$; an explicit minimal homogeneous resolution of $Q$ by free $P$-modules, and an explicit minimal homogeneous resolution of $Q$ by free $R$-modules are all given in Lemma \tref{setup}. (If $r$ had been $0$ or $n$, then Lemma \tref{setup} would still produce a generating set for $(x^N,y^N,z^N)\:\! f$ and free resolutions of $Q$; however these objects would not be minimal.) To complete the proof of (1), we need only record the degrees of the entries of the matrices in these resolutions. We know that $\deg f=n$, each entry of $\pmb x$ has degree $N$,
$$\deg b_1=\deg b_2=\deg b_3= 3\delta n -2n+r,\quad \deg b_4=3\delta n-3n+3r,$$ each entry in columns $1$, $2$, and $3$ of $\psi \varphi\check{\phantom{.}}$ has degree $\delta n$, each entry of column 4 of $\psi \varphi\check{\phantom{.}}$ has degree $\delta n -n+2r$, each entry of columns $1$, $2$, and $3$ of $\psi$ has degree $\delta n-n+r$, column four of $\psi$ is a zero matrix, and 
$$\varphi=\bmatrix\format\l&\qquad\l\\ \matrix\format\l\\ \text{$3\times3$ matrix with entries}\\\text{of degree $r$}\endmatrix& \matrix\format\l\\\text{$3\times1$ matrix with entries}\\\text{of degree $n-r$}\endmatrix\\\vspace{5pt}\matrix\format\l\\\text{$1\times3$ matrix with entries}\\\text{of degree $n-r$}\endmatrix&
 \text{$1\times1$ zero matrix } 
\endbmatrix.$$ The socle degrees of $Q$ may be read from the $P$-resolution of $Q$; see Remark \tref{R?}.

We now focus on assertion (2); that is, we have  $\theta=2\delta$ is even. Recall, from (\tref{ssetup}), that  $\varphi$ is   the matrix $\varphi_{n-r,r}$ of {\rm(\tref{frn})}. It follows that $\varphi\check{\phantom{.}}$ is the matrix $-\varphi_{r,n-r}$.  We still have $\varphi\varphi\check{\phantom{.}}=\varphi\check{\phantom{.}}\varphi=fI$. 
Define $\psi$ to be the matrix $$\eightpoint  \left[\matrix\format \c&\ \c&\ \c&\ \r\\ 
\operatorname{Poly}'_{\delta,\delta,\delta}(y^n,z^n)&0&0&y^rz^r\operatorname{Poly}'_{\delta-1,\delta,\delta}(y^n,z^n)\\
0&(-1)^{\delta}\operatorname{Poly}'_{\delta,\delta,\delta}(x^n,z^n)&0&(-1)^{\delta+1}x^rz^r\operatorname{Poly}'_{\delta-1,\delta,\delta}(x^n,z^n)\\
0&0&\operatorname{Poly}'_{\delta,\delta,\delta}(x^n,y^n)&x^ry^r\operatorname{Poly}'_{\delta-1,\delta,\delta}(x^n,y^n)\endmatrix\right] $$ and $$\Phi= 
\bmatrix 
0& z^r \Cal R_{2\delta}(z^n,x^n,y^n)&-y^r \Cal R_{2\delta}(y^n,x^n,z^n)\\
-z^r\Cal R_{2\delta}(z^n,x^n,y^n)&0& x^r\Cal R_{2\delta}(x^n,y^n,z^n)\\
 y^r \Cal R_{2\delta}(y^n,x^n,z^n)&-x^r\Cal R_{2\delta}(x^n,y^n,z^n)&0\endbmatrix.$$
One may easily calculate that $\psi\varphi\check{\phantom{.}}\psi^{\text{\rm T}}$ is equal to 
$$-\bmatrix 0&-z^r\frak P'_{2\delta}(z^n,x^n,y^n)&y^r\frak P'_{2\delta}(y^n,x^n,z^n)\\z^r\frak P'_{2\delta}(z^n,x^n,y^n)&0&-x^r\frak P'_{2\delta}(x^n,y^n,z^n)\\-y^r\frak P'_{2\delta}(y^n,x^n,z^n)&x^r\frak P'_{2\delta}(x^n,y^n,z^n)&0\endbmatrix.$$
Therefore, the definition (\tref{zzz}) of $\Cal R_{2\delta}$ tells us that 
$$-\psi\varphi\check{\phantom{.}}\psi^{\text{\rm T}}+f\Phi= u_{2\delta}\pmb X.$$ 
The condition (\tref{hyp20}) is satisfied. Lemma \tref{setup} applies. The proof of (2) is completed just like the proof of (1).   We know that $\deg f=n$, each entry of $\pmb x$ has degree $N$,
$$\deg b_1=\deg b_2=\deg b_3= 3\delta n -n+2r,\quad \deg b_4=3\delta n,$$ each entry in columns $1$, $2$, and $3$ of $\psi \varphi\check{\phantom{.}}$ has degree $\delta n+r$, each entry of column 4 of $\psi \varphi\check{\phantom{.}}$ has degree $\delta n +n-r$, each entry of columns $1$, $2$, and $3$ of $\psi$ has degree $\delta n$, each entry in column four of $\psi$ has degree $\delta n -n+2r$, and 
$$\varphi=\bmatrix\format\l&\qquad\l\\ \matrix\format\l\\ \text{$3\times3$ matrix with entries}\\\text{of degree $n-r$}\endmatrix& \matrix\format\l\\\text{$3\times1$ matrix with entries}\\\text{of degree $r$}\endmatrix\\\vspace{5pt}\matrix\format\l\\\text{$1\times3$ matrix with entries}\\\text{of degree $r$}\endmatrix&
 \text{$1\times1$ zero matrix } 
\endbmatrix. $$  

The proof is complete when $\operatorname{char} \pmb k\neq 3$. Henceforth, we consider $\operatorname{char} \pmb k=3$. The idea of the proof is exactly the same as in the first case; but all of the details change a little bit. Let 
$$d=\cases \delta&\text{if $\theta$ is even}\\\delta-1&\text{if $\theta$ is odd}.\endcases$$It follows that $$\theta=\cases 2d &\text{if $\theta$ is even}\\2d+1&\text{if $\theta$ is odd}.\endcases$$The parameter $\theta$ is in $S_3$ by hypothesis; hence, Remark \tref{!R18.1} shows that $d\in D_3$ and therefore all five statements of part (1) of Theorem \tref{4} hold for $d$. 

Define the integers $\gamma$ and $\Gamma$ by  
$\gamma=3^{\ell}$ and $\Gamma=3^L$  where $\binom{2 d}d_{\#3}=\ell$ and $\binom{2d+1}d_{\#3}=L$. Theorem \tref{4} shows that
$\binom{3d}{d}_{\#3}=\binom{2d}{d}_{\#3}$; so $\gamma^2$ divides $\binom{2d}{d}\binom{3d}{d}$ in $\Bbb Z$ and 
$$u_{2d}=(-1)^{d}\binom{2d}{d}\binom{3d}{d}/\gamma^2$$ is a unit in $\pmb k$. Theorem \tref{4} also  guarantees that
$\binom{2d+1}d_{\#3}=\binom{3d+2}{d}_{\#3}$. Recall that $\binom{2d+2}{d+1}=2\binom{2d+1}d$; so, $\binom{2d+2}{d+1}_{\#3}=\binom{2d+1}d_{\#3}$. It follows that $\Gamma^2$ divides 
$\binom{2d+2}{d+1}\binom{3d+2}{d}$
 in $\Bbb Z$ and 
$$u_{2d+1}=(-1)^d\binom{2d+2}{d+1}\binom{3d+2}{ d }/\Gamma^2$$ is a unit in $\pmb k$. According to (1c) and (1d) from Theorem \tref{4}, $\gamma$ divides   every coefficient of each of the polynomials
$$\operatorname{Poly}_{d,d,d}(A,B),\quad \operatorname{Poly}_{d-1,d,d}(A,B);$$ furthermore,
(1e) from Theorem \tref{4} shows that $\Gamma$ divides   every coefficient of   the polynomial
$$\operatorname{Poly}_{d,d+1,d}(A,B).$$
  Let
$$\matrix \format \r &\ \l\ &\l \quad &\r &\ \l\ &\l\\
\operatorname{Poly}'_{d,d,d}&=&\operatorname{Poly}_{d,d,d}/\gamma&\operatorname{Poly}'_{d,d+1,d}&=&\operatorname{Poly}_{d,d+1,d}/\Gamma,\\ 
\operatorname{Poly}'_{d-1,d,d}&=&\operatorname{Poly}_{d-1,d,d}/\gamma,&\frak P_{2d+1}'&=&\frak P_{2d+1}/\Gamma^2,\ \text{and}\\  \frak P_{2d}'&=&\frak P_{2d}/\gamma^2\endmatrix $$ in $\Bbb Z[A,B,C]$. 
We have $$ \matrix\format\l&\l\\ \Gamma^2(\frak P'_{2d+1}(A,B,C)+u_{2d+1}A^{2d+1})&{}={}\frak P_{2d+1}(A,B,C)+(-1)^{d}\binom{2d+2}{d+1}\binom{3d+2}{d+1}A^{2d+1}\\&{}={}P_{2d+1}(A,B,C)\ \text{and}\endmatrix$$
$$\matrix  \gamma^2(\frak P'_{2d}(A,B,C)+u_{2d}A^{2d})= \frak P_{2d}(A,B,C)+(-1)^{d}\binom{2d}{d}\binom{3d}{d}A^{2d}=P_{2d}(A,B,C).\endmatrix $$Therefore, according to Lemma \tref{NandI},
$$\Gamma^2(\frak P'_{2d+1}(A,B,C)+u_{2d+1}A^{2d+1}) \quad\text{and}\quad \gamma^2(\frak P'_{2d}(A,B,C)+u_{2d}A^{2d})$$ are both in the prime ideal $(A+B+C)\Bbb Z[A,B,C]$  of $\Bbb Z[A,B,C]$. It follows that 
$$ \frak P'_{2d+1}(A,B,C)+u_{2d+1}A^{2d+1}  \quad\text{and}\quad  \frak P'_{2d}(A,B,C)+u_{2d}A^{2d} $$ 
are in $(A+B+C)\Bbb Z[A,B,C]$. Define
$\Cal R_{2d+1}(A,B,C)$ and $\Cal R_{2d}(A,B,C)$ to be the polynomials in $\Bbb Z[A,B,C]$ which satisfy
$$\matrix\format\l&\c&\l\\ \frak P'_{2d+1}(A,B,C)+ u_{2d+1}A^{2d+1}&{}={}&(A+B+C)\Cal R_{2d+1}(A,B,C)
\ \text{and}\\ \frak P'_{2d}(A,B,C)+ u_{2d}A^{2d}&{}={}&(A+B+C)\Cal R_{2d}(A,B,C)
.\endmatrix  \tag\tnum{zzz'}$$

We now focus on assertion (1); that is, we have  $\theta=2d+1$ is odd. Recall, from (\tref{ssetup}), that  $\varphi$ is   the matrix $\varphi_{r,n-r}$ of {\rm(\tref{frn})}. It follows that $\varphi\check{\phantom{.}}$ is the matrix $-\varphi_{n-r,r}$.  We see from (\tref{ffhat}) that $\varphi\varphi\check{\phantom{.}}=\varphi\check{\phantom{.}}\varphi=fI$. 
Define $\psi$ to be the matrix $$\eightpoint \left[\smallmatrix 
0&(-1)^{d}z^r\operatorname{Poly}'_{d,d+1,d}(z^n,y^n)&y^r\operatorname{Poly}'_{d,d+1,d}(y^n,z^n)&0  \\
z^r\operatorname{Poly}'_{d,d+1,d}(z^n,x^n)&0&(-1)^{d}x^r\operatorname{Poly}'_{d,d+1,d}(x^n,z^n)&0 \\
(-1)^{d}y^r\operatorname{Poly}'_{d,d+1,d}(y^n,x^n)&x^r\operatorname{Poly}'_{d,d+1,d}(x^n,y^n)&0&0\endsmallmatrix\right] $$ and $$\Phi= 
\bmatrix 
0& z^r \Cal R_{2d+1}(z^n,x^n,y^n)&-y^r \Cal R_{2d+1}(y^n,x^n,z^n)\\
-z^r\Cal R_{2d+1}(z^n,x^n,y^n)&0& x^r\Cal R_{2d+1}(x^n,y^n,z^n)\\
 y^r \Cal R_{2d+1}(y^n,x^n,z^n)&-x^r\Cal R_{2d+1}(x^n,y^n,z^n)&0\endbmatrix.$$
One may easily calculate that $\psi\varphi\check{\phantom{.}}\psi^{\text{\rm T}}$ is equal to 
$$\bmatrix 0&-z^r\frak P'_{2d+1}(z^n,x^n,y^n)&y^r\frak P'_{2d+1}(y^n,x^n,z^n)\\z^r\frak P'_{2d+1}(z^n,x^n,y^n)&0&-x^r\frak P'_{2d-1}(x^n,y^n,z^n)\\-y^r\frak P'_{2d+1}(y^n,x^n,z^n)&x^r\frak P'_{2d+1}(x^n,y^n,z^n)&0\endbmatrix.$$
Therefore, the definition (\tref{zzz'}) of $\Cal R_{2d+1}$ tells us that 
$$\psi\varphi\check{\phantom{.}}\psi^{\text{\rm T}}+f\Phi= u_{2d+1}\pmb X.$$ 
The condition (\tref{hyp20}) is satisfied. Lemma \tref{setup} applies. One completes the proof of (1) when $\operatorname{char}\pmb k=3$ just like one completes the proof of (1) when $\operatorname{char}\pmb k\neq 3$. Once one has the resolution of $Q$ written in terms of the parameter $d$, then the equation $d=\delta-1$ may be used to write the resolution in the form which is recorded in (1).

Now that all of the words have been defined, one completes the proof of $\theta=2d$ and $\operatorname{char}\pmb k=3$ using the exact same words as were used in the proof of $\theta=2d$ and $\operatorname{char}\pmb k\neq 3$. $\qed$
\enddemo

\SectionNumber=4\tNumber=1
\heading Section \number\SectionNumber. \quad The proof of Theorem \tref{4}.
\endheading

In this section we prove Theorem \tref{4}: Proposition \tref{!L18.2} is the second assertion   and Proposition \tref{P8-13} is the first assertion.  The proof is carried out by induction. The key step in this induction  is Lemmas \tref{Aug21} and \tref{L8} where we show how the numbers of the form $\binom AD_{\#p}$ involved in Theorem \tref{4}, with $A=ap+r$ and $D=dp+\epsilon$, for small values of $r$ and $\epsilon$,  are related to numbers of the form $\binom ad_{\#p}$.

\proclaim{Lemma \tnum{Aug21}} Let $p$, $d$, and $\epsilon$ be integers with $p\ge 3$ and prime, $d\ge 1$,  and $-\frac p3< \epsilon< \frac p3$.
Let $a$,  $r$, $s$, and $t$ be non-negative integers with $ r\le p-1$, $ s\le 1$, and $ t\le 2$. If $A=pa+r$ and $D=pd+\epsilon$, then 
\flushpar{\bf(a)} 
$$\binom{2D}D_{\#p}=\cases 
\left(pd\binom{2d}{d}\right)_{\#p}&\text{if $\epsilon<0$}\\
\binom{2d}d_{\#p}&\text{if $0\le\epsilon$},\endcases$$
\flushpar{\bf(b)} $$\binom{3D}D_{\#p}=\cases 
\left(\frac{pd}3\binom{3d}{d}\right)_{\#p}&\text{if $\epsilon<0$}\\
\binom{3d}d_{\#p}&\text{if $0\le\epsilon$},\endcases$$
\flushpar{\bf(c)}
$$\binom{A}{D-s}_{\#p}=\cases 
\left(p(a-d)\binom{a}d\right)_{\#p}&\text{if $r+s+1-\epsilon\le 0$}\\
\binom ad_{\#p}&\text{if $s\le \epsilon$ and $1\le r+s+1-\epsilon\le p$}\\
\left(pd\binom ad\right)_{\#p}&\text{if $\epsilon\le s-1$ and $1\le r+s+1-\epsilon\le p$}\\
\binom a{d-1}_{\#p}&\text{if $\epsilon\le s-1$ and  $p+1\le r+s+1-\epsilon$},\\
\endcases$$
\flushpar{\bf(d)}
if $0\le \epsilon$, then  
$$\binom{3D-A-t}D_{\#p}\hskip-.8pt=\cases
\binom{3d-a-2}d_{\#p}&\text{if $3\epsilon+p+1\le r+t$}   
\\
\left(p(2d-a-1)\binom{3d-a-1}d\right)_{\#p}&\text{if  
$2\epsilon+p+1\le r+t\le 3\epsilon+p$} 
\\
\binom{3d-a-1}d_{\#p} &\text{if $3\epsilon+1\le r+t\le p+2\epsilon$} 
\\
\left(p(2d-a) \binom{3d-a}{d}\right)_{\#p}&\text{if $2\epsilon+1\le r+t\le 3\epsilon$} 
\\
\binom{3d-a}d_{\#p}&\text{if  $r+t\le 2\epsilon$, and} 
\endcases
$$
\flushpar{\bf(e)}
if $\epsilon\le -1$, then  
$$\binom{3D-A-t}D_{\#p}= \cases
\left(pd\binom{3d-a-2}d\right)_{\#p}&\text{if  
$2\epsilon+p+1\le r+t$} 
\\
\binom{3d-a-2}{d-1}_{\#p}&  \text{if $3\epsilon+p+1\le r+t\le 2\epsilon+p$  } 
\\ 
\left(pd\binom{3d-a-1}d\right)_{\#p}&\text{if $r+t\le p+3\epsilon$}.  
\endcases
$$
\endproclaim

\demo{Proof} The key to all of these calculations is the observation that if $L$ is a set of integers, then
$$\left(\prod\limits_{\{\ell\in L\}}\ell\right)_{\#p}=\left(\prod\limits_{\{\ell\in L\mid \text{ $p$ divides $\ell$}\}}\ell\right)_{\#p}.$$In particular, for example, if $0\le \epsilon<\frac p3$, then
$$\left((pd+\epsilon)!\right)_{\#p}=\left([p(d)]\cdot [p(d-1)]\cdots [p(2)] \cdot [p(1)]\right)_{\#p}=\left(p^dd!\right)_{\#p}.$$

 Assertion (a) is easy to verify if $d=1$. We assume $2\le d$. We have 
$$\binom{2D}D_{\#p}=\left(\frac{(2D)\cdots (D+1)}{D!}\right)_{\#p}=\left(\frac{[p(2d)+2\epsilon]\cdots [p(d)+\epsilon+1]}{(pd+\epsilon)!}\right)_{\#p}$$
$$=
\cases 
\left(\frac{[p(2d-1)]\cdots [pd]}{[p(d-1)]\cdots [p(1)]}\right)_{\#p}=\left(pd\binom{2d-1}{d-1}\right)_{\#p}&\text{if $\epsilon<0$}\\
\left(\frac{[p(2d)]\cdots [p(d+1)]}{[p(d)]\cdots [p(1)]}\right)_{\#p}=\binom{2d}d_{\#p}&\text{if $0\le\epsilon$}.\endcases
$$ Notice that $ \binom{2d-1}{d-1}=\frac 12 \binom{2d}d$ and $\left(\frac 12 \binom{2d}d\right)_{\#p}=\binom{2d}d_{\#p}$, since $p\neq 2$. 

The  technique that was used for (a) shows that
$$\binom{3D}D_{\#p}=\cases \left(pd\binom{3d-1}{d-1}\right)_{\#p}&\text{if $\epsilon<0$}\\
\binom{3d}d_{\#p}&\text{if $0\le\epsilon$}.\endcases
$$Furthermore, we know that  $\binom{3d-1}{d-1}=\frac 13\binom{3d}{d}$.

 We now prove (c).  We know that  $ \binom{A}{D-s}_{\#p}$ is equal to
$$\left(\frac{A\cdots (A-D+s+1)}{(D-s)!}\right)_{\#p}
=\left(\frac{ [pa+r]\cdots [p(a-d)+r+s+1-\epsilon]}{(pd+\epsilon-s)!}\right)_{\#p}.$$
Notice that that $-p<\epsilon-s<p$ and that $-p<r+s+1-\epsilon<2p$. It follows that 
$$\left((pd+\epsilon-s)!\right)_{\#p}=\cases 
\left([p(d)]\cdots [p(1)]\right)_{\#p}&\text{if $s\le \epsilon$}\\
\left([p(d-1)]\cdots [p(1)]\right)_{\#p}&\text{if $\epsilon\le s-1$}\\
\endcases$$and $\left([pa+r]\cdots [p(a-d)+r+s+1-\epsilon]\right)_{\#p}$
$$=\cases 
\left([pa]\cdots [p(a-d)]\right)_{\#p}&\text{if $r+s+1-\epsilon\le 0$}\\
\left([pa]\cdots [p(a-d+1)]\right)_{\#p}&\text{if $1\le r+s+1-\epsilon\le p$}\\
\left([pa]\cdots [p(a-d+2)]\right)_{\#p}&\text{if $p+1\le r+s+1-\epsilon$.}\endcases$$
Observe that $s\le \epsilon\implies r+s+1-\epsilon\le p$ and $\epsilon\le s-1\implies 1\le r+s+1-\epsilon$. There are a total of four options:
$ \binom{A}{D-s}_{\#p}$ is equal to 
$$ \cases 
\left(\frac{ [pa]\cdots [p(a-d)]}{[p(d)]\cdots [p(1)]}\right)_{\#p}=\left(p(a-d)\binom{a}d\right)_{\#p}&\text{if $s\le \epsilon$ and $r+s+1-\epsilon\le 0$}\\
\left(\frac{ [pa]\cdots [p(a-d+1)]}{[p(d)]\cdots [p(1)]}\right)_{\#p}=\binom ad_{\#p}&\text{if $s\le \epsilon$ and $1\le r+s+1-\epsilon\le p$}\\
\left(\frac{ [pa]\cdots [p(a-d+1)]}{[p(d-1)]\cdots [p(1)]}\right)_{\#p}=\left(pd\binom ad\right)_{\#p}&\text{if $\epsilon\le s-1$ and $1\le r+s+1-\epsilon\le p$}\\
\left(\frac{ [pa]\cdots [p(a-d+2)]}{[p(d-1)]\cdots [p(1)]}\right)_{\#p}=\binom a{d-1}_{\#p}&\text{if $\epsilon\le s-1$ and  $p+1\le r+s+1-\epsilon$}.
\endcases$$
In the top case, $r+s+1-\epsilon\le 0\implies 1\le r+1\le \epsilon-s$, so the constraint $s\le \epsilon$ is redundant.
This calculation is valid for $1\le d$; however, special care is needed  when $d=1$. In this case,
$p$ does not divide $(pd+\epsilon-s)!$ when $\epsilon-s\le -1$ and $p$ does not divide $[pa+r]\cdots [p(a-d)+r+s+1-\epsilon]$ when $p+1\le r+s+1-\epsilon$.

We now prove (d) and (e). The argument we give is valid for all $d\ge 1$. We have 
$$ \binom{3D-A-t}D_{\#p}=\left(\frac{(3D-A-t)\cdots (2D-A-t+1)}{D!}\right)_{\#p}$$
$$ =\left(\frac{[p(3d-a)+3\epsilon-r-t]\cdots [p(2d-a)+2\epsilon-r-t+1]}{(pd+\epsilon)!}\right)_{\#p}.$$
We know that 
$$((pd+\epsilon)!)_{\#p}= \cases
\left({[p(d)]\cdots [p(1)]}\right)_{\#p}&\text{if $0\le \epsilon$}\\
\left(p^{d-1}(d-1)!\right)_{\#p}&\text{if $\epsilon\le -1$}.\endcases
$$
We see that   $$0\le \epsilon\implies -2p< 3\epsilon-r-t<p\quad\text{and}\quad \epsilon\le -1\implies -2p\le 3\epsilon-r-t<0.$$ Thus, $\binom{3D-A-t}D_{\#p}$ is equal to 
$$\cases
\left(\frac{[p(3d-a-2)]\cdots [p(2d-a)+2\epsilon-r-t+1]}{[p(d)]\cdots [p(1)]}\right)_{\#p}&\text{if $0\le \epsilon$ and
$3\epsilon-r-t\le -p-1$ }\\
\left(\frac{[p(3d-a-1)]\cdots [p(2d-a)+2\epsilon-r-t+1]}{[p(d)]\cdots [p(1)]}\right)_{\#p}&\text{if $0\le \epsilon$ and 
$-p\le 3\epsilon-r-t\le -1$}\\
\left(\frac{[p(3d-a)]\cdots [p(2d-a)+2\epsilon-r-t+1]}{[p(d)]\cdots [p(1)]}\right)_{\#p}&\text{if $0\le \epsilon$ and 
$0\le 3\epsilon-r-t$}\\
\left(\frac{[p(3d-a-2)]\cdots [p(2d-a)+2\epsilon-r-t+1]}{ p^{d-1}(d-1)!}\right)_{\#p}&\text{if $\epsilon\le -1$ and  $3\epsilon-r-t\le -p-1$}\\
\left(\frac{[p(3d-a-1)]\cdots [p(2d-a)+2\epsilon-r-t+1]}{ p^{d-1}(d-1)!}\right)_{\#p}&\text{if $\epsilon\le -1$ and  $-p\le 3\epsilon-r-t$.}
\\
\endcases
$$We see that $$-2p<-\frac{5p}3=-\frac{2p}3-(p-1)-2+1\le 2\epsilon-r-t+1< \frac{2p}3+1\le p;$$ and therefore, the smallest multiple of $p$ which is at least $[p(2d-a)+2\epsilon-r-t+1]$ is 
$$\cases p(2d-a-1)&\text{if $2\epsilon-r-t+1\le -p$}\\
p(2d-a)&\text{if $-p+1\le 2\epsilon-r-t+1\le 0$}\\
p(2d-a+1)&\text{if $1\le 2\epsilon-r-t+1$.}\endcases$$
If $0\le \epsilon$, then  $2\epsilon\le 3\epsilon\le p+2\epsilon\le p+3\epsilon$ and  the constraints 
$$\left\{\matrix\format\l\\  3\epsilon-r-t\le -p-1\\
-p\le 3\epsilon-r-t\le -1\\
0\le 3\epsilon-r-t\endmatrix\right.\quad\text{and}\quad 
\left\{\matrix\format\l\\  2\epsilon-r-t+1\le -p \\
 -p+1\le 2\epsilon-r-t+1\le 0\\
1\le 2\epsilon-r-t+1 \endmatrix\right.$$ may be merged to produce
$$\left\{\matrix\format\l\\  3\epsilon+p+1\le r+t\\   
2\epsilon+p+1\le r+t\le 3\epsilon+p\\
3\epsilon+1\le r+t\le p+2\epsilon\\
2\epsilon+1\le r+t\le 3\epsilon\\
r+t\le 2\epsilon.\endmatrix\right.$$ Thus, if $0\le \epsilon$, then $\binom{3D-A-t}D_{\#p}$ is equal to 
$$\eightpoint \cases
\left(\frac{[p(3d-a-2)]\cdots [p(2d-a-1)]}{[p(d)]\cdots [p(1)]}\right)_{\#p}=\binom{3d-a-2}d_{\#p}&\text{if  $3\epsilon+p+1\le r+t$}  
\\
\left(\frac{[p(3d-a-1)]\cdots [p(2d-a-1)]}{[p(d)]\cdots [p(1)]}\right)_{\#p}=\left(p(2d-a-1)\binom{3d-a-1}d\right)_{\#p}&\text{if  
$2\epsilon+p+1\le r+t\le 3\epsilon+p$}
\\
\left(\frac{[p(3d-a-1)]\cdots [p(2d-a)]}{[p(d)]\cdots [p(1)]}\right)_{\#p}=\binom{3d-a-1}d_{\#p} &\text{if   $3\epsilon+1\le r+t\le p+2\epsilon$} 
\\
\left(\frac{[p(3d-a)]\cdots [p(2d-a)]}{[p(d)]\cdots [p(1)]}\right)_{\#p}=\left(p(2d-a) \binom{3d-a}{d}\right)_{\#p}&\text{if $2\epsilon+1\le r+t\le 3\epsilon$}  
\\
\left(\frac{[p(3d-a)]\cdots [p(2d-a+1)]}{[p(d)]\cdots [p(1)]}\right)_{\#p}=\binom{3d-a}d_{\#p}&\text{if  $r+t\le 2\epsilon$.}  
\endcases
$$In a similar manner, if $\epsilon\le -1$, then  
$$2\epsilon<0<3\epsilon+p<2\epsilon+p$$ and
the constraints 
$$\left\{\matrix\format\l\\  3\epsilon-r-t\le -p-1\\
-p\le 3\epsilon-r-t\endmatrix\right.\quad\text{and}\quad 
\left\{\matrix\format\l\\   2\epsilon-r-t+1\le -p \\
 -p+1\le 2\epsilon-r-t+1\le 0\\
1\le 2\epsilon-r-t+1 \endmatrix\right.$$ may be merged to produce
$$\left\{\matrix\format\l\\ 2\epsilon+p+1\le r+t\\3\epsilon+p+1\le r+t\le 2\epsilon+p\\r+t\le p+3\epsilon,\endmatrix\right.$$since $r+t$ is automatically non-negative.
Therefore, if $\epsilon\le -1$, then $\binom{3D-A-t}D_{\#p}$ is equal to 
$$\cases
\left(\frac{[p(3d-a-2)]\cdots [p(2d-a-1)]}{ p^{d-1}(d-1)!}\right)_{\#p}=\left(pd\binom{3d-a-2}d\right)_{\#p}&\text{if  
$2\epsilon+p+1\le r+t$}
\\
\left(\frac{[p(3d-a-2)]\cdots [p(2d-a)]}{ p^{d-1}(d-1)!}\right)_{\#p}=\binom{3d-a-2}{d-1}_{\#p}&\text{if $3\epsilon+p+1\le r+t\le 2\epsilon+p$} 
\\ 
\left(\frac{[p(3d-a-1)]\cdots [p(2d-a)]}{ p^{d-1}(d-1)!}\right)_{\#p}=\left(pd\binom{3d-a-1}d\right)_{\#p}&\text{if $r+t\le p+3\epsilon$}. \qed
\endcases
$$\enddemo

\proclaim{Proposition \tnum{!L18.2}} Assertion {\rm(2)} of Theorem {\rm \tref{4}} holds.\endproclaim

\demo{Proof} Fix a prime integer $p$ with $p\ge 5$.  All four assertions are  clear for $0\le d< \frac p3$ because $p$ does not divide either $\binom{2d}d$ or $\binom{3d}d$.  
The proof proceeds by induction. Let $D$ be an element of $D_p$ which is greater than $\frac p3$. The recursive definition of $D_p$ guarantees that $D=pd+\epsilon$ for some $d\in D_p$ and some integer $\epsilon$ with $-\frac p3< \epsilon< \frac p3$. We see that $d<D$; so the induction hypothesis guarantees that  $\binom{2d}d_{\#p}=\binom{3d}d_{\#p}$. According to Lemma \tref{Aug21} we have
$$\cases\binom{2D}D_{\#p}=\binom{2d}d_{\#p}=\binom{3d}d_{\#p}=\binom{3D}D_{\#p}&\text{if $0\le \epsilon$ and}\\
\binom{2D}D_{\#p}=\left(pd\binom{2d}d\right)_{\#p}=\left(pd\binom{3d}d\right)_{\#p}=\left(\frac{pd}3\binom{3d}d\right)_{\#p}=\binom{3D}D_{\#p}&\text{if $\epsilon\le -1$.}\endcases$$ So $\binom{2D}D_{\#p}=\binom{3D}D_{\#p}$ in all cases and assertion (a) is established. 

The induction hypothesis also guarantees that   all of  the inequalities (b), (c), and (d)  hold when $0\le a\le 2d$. Notice also that the induction hypothesis ensures that
$$ \left(\binom{3d-a-2}{d-1}\binom ad\right)_{\#p}\ge \binom{2d}d_{\#p}\tag\tnum{!-also}$$ also holds for $0\le a\le 2d$ because this is a different form of (d).

Fix $A=pa+r$  with $$0\le a,\quad 0\le r\le p-1,  \quad\text{and}\quad 0\le A\le 2D.$$ 
We prove that all three inequalities (b)--(d) hold at the pair $(A,D)$.   

Observe first that $0\le a\le 2d$; so all four inequalities (b)--(d) and (\tref{!-also}) hold at $(a,d)$ by induction. 
 
We first establish (b) and (c) at $(A,D)$, when $0\le \epsilon$. Fix $t$ with $0\le t\le 1$ and $0\le \epsilon$. We study 
$$\left(\binom AD\binom {3D-A-t}A\right)_{\#p}.$$Lemma \tref{Aug21} shows that  
$$\binom AD_{\#p}=\cases 
\left(p(a-d)\binom{a}d\right)_{\#p}&\text{if $r\le \epsilon-1$}\\
\binom ad_{\#p}&\text{if $\epsilon\le r$}.
\endcases\tag\tnum{!-17.4}$$
The hypotheses $0\le \epsilon$ and $t\le 1$ ensure that $r+t\le 3\epsilon+p$; so, Lemma \tref{Aug21} also gives
$$\eightpoint \binom {3D-A-t}A_{\#p}=\cases
\left(p(2d-a-1)\binom{3d-a-1}d\right)_{\#p}&\text{if  
$2\epsilon+p+1\le r+t$} 
\\
\binom{3d-a-1}d_{\#p} &\text{if $3\epsilon+1\le r+t\le p+2\epsilon$} 
\\
\left(p(2d-a) \binom{3d-a}{d}\right)_{\#p}&\text{if $2\epsilon+1\le r+t\le 3\epsilon$} 
\\
\binom{3d-a}d_{\#p}&\text{if  $r+t\le 2\epsilon$.} 
\endcases\tag\tnum{!-17.5}$$We combine (\tref{!-17.4}) and (\tref{!-17.5}). Observe that the condition $r\le \epsilon-1$ forces $r+t\le 2\epsilon$. It follows that  $\left(\binom AD\binom {3D-A-t}A\right)_{\#p}$ is equal to 
$$\cases
\left(\binom adp(2d-a-1)\binom{3d-a-1}d\right)_{\#p}&\text{if  
$2\epsilon+p+1\le r+t$} 
\\
\left(\binom ad\binom{3d-a-1}d\right)_{\#p} &\text{if $3\epsilon+1\le r+t\le p+2\epsilon$} 
\\
\left(\binom adp(2d-a) \binom{3d-a}{d}\right)_{\#p}&\text{if $2\epsilon+1\le r+t\le 3\epsilon$} 
\\
\left(\binom ad\binom{3d-a}d\right)_{\#p}&\text{if  $r+t\le 2\epsilon$ and $\epsilon\le r$} \\
\left(p(a-d)\binom{a}d\binom{3d-a}d\right)_{\#p}&\text{if  $r+t\le 2\epsilon$ and $r\le \epsilon-1$}.
\endcases$$ The hypothesis that all (b)--(d) hold at $(a,d)$   guarantees that
$$\left(\binom AD\binom {3D-A-t}A\right)_{\#p}\ge \binom{2d}d_{\#p}=\binom{2D}D_{\#p},$$ as claimed. 

We next establish (b) and (c) at $(A,D)$   simultaneously, when $\epsilon\le -1$. Fix $t$ with $0\le t\le 1$ and $\epsilon\le -1$. We study 
$$\left(\binom AD\binom {3D-A-t}A\right)_{\#p}.$$
Apply part (c) of Lemma \tref{Aug21} with $s=0$ and $\epsilon\le -1$. The cases $r+s+1-\epsilon\le 0$ and $s\le \epsilon$ do not occur in this situation and the condition $\epsilon\le r$ holds automatically; so, 
we have   $$\binom{A}{D}_{\#p}=\cases 
\left(pd\binom ad\right)_{\#p}&\text{$r\le p-1+\epsilon$}\\
\binom a{d-1}_{\#p}&\text{$p+\epsilon\le r$}.\\
\endcases\tag\tnum{!-17.8}$$Part (e) of Lemma \tref{Aug21} gives
$$\binom{3D-A-t}D_{\#p}= \cases
\left(pd\binom{3d-a-2}d\right)_{\#p}&\text{if  
$2\epsilon+p+1\le r+t$} 
\\
\binom{3d-a-2}{d-1}_{\#p}&  \text{if $3\epsilon+p+1\le r+t\le 2\epsilon+p$  } 
\\ 
\left(pd\binom{3d-a-1}d\right)_{\#p}&\text{if $r+t\le p+3\epsilon$}.  
\endcases\tag\tnum{!-17.9}
$$We notice that $p+\epsilon\le r$ implies
$$p+2\epsilon+1=p+\epsilon+(\epsilon+1)\le r+0\le r+t;$$ furthermore, $p+3\epsilon<p+2\epsilon<p+2\epsilon+1$. Combine (\tref{!-17.8}) and (\tref{!-17.9}) to see that 
$\left(\binom{A}{D} \binom{3D-A-t}D\right)_{\#p}$ is equal to 
$$\cases
\left(\binom a{d-1}pd\binom{3d-a-2}d\right)_{\#p}&\text{if  
$2\epsilon+p+1\le r+t$ and $p+\epsilon\le r$} 
\\
\left(pd\binom adpd\binom{3d-a-2}d\right)_{\#p}&\text{if  
$2\epsilon+p+1\le r+t$ and $r\le p-1+\epsilon$} 
\\
\left(pd\binom ad\binom{3d-a-2}{d-1}\right)_{\#p}&  \text{if $3\epsilon+p+1\le r+t\le 2\epsilon+p$  } 
\\ 
\left(pd\binom adpd\binom{3d-a-1}d\right)_{\#p}&\text{if $r+t\le p+3\epsilon$}.  
\endcases
$$
Notice that $d\binom ad=\binom a{d-1}(a-d+1)$; so the second option in the above display is
$\left(p(a-d+1)\binom a{d-1}pd\binom{3d-a-2}d\right)_{\#p}$. For each case, apply one of the inequalities (b)--(d) or (\tref{!-also}) at $(a,d)$. We conclude that
  $$\left(\binom{A}{D} \binom{3D-A-t}D\right)_{\#p}\ge \left(pd\binom {2d}d\right)_{\#p}=\binom {2D}D_{\#p}.$$

We next establish $(d)$ when $0\le \epsilon$. Lemma \tref{Aug21} gives
$$\binom{A}{D-1}_{\#p}=\cases 
\left(p(a-d)\binom{a}d\right)_{\#p}&\text{if $r\le \epsilon-2$}\\
\binom ad_{\#p}&\text{if $1\le \epsilon$ and $\epsilon-1\le r\le p-2+\epsilon$}\\
\left(pd\binom ad\right)_{\#p}&\text{if $\epsilon= 0$ and $ r\le p-2$}\\
\binom a{d-1}_{\#p}&\text{if $\epsilon= 0$ and  $p-1= r$}\\
\endcases$$
and  $\binom{3D-A-2}D_{\#p}$ is equal to 
$$\cases
\binom{3d-a-2}d_{\#p}&\text{if $3\epsilon+p-1\le r$}   
\\
\left(p(2d-a-1)\binom{3d-a-1}d\right)_{\#p}&\text{if  
$2\epsilon+p-1\le r\le 3\epsilon+p-2$} 
\\
\binom{3d-a-1}d_{\#p} &\text{if $3\epsilon-1\le r\le p-2+2\epsilon$} 
\\
\left(p(2d-a) \binom{3d-a}{d}\right)_{\#p}&\text{if $2\epsilon-1\le r\le 3\epsilon-2$} 
\\
\binom{3d-a}d_{\#p}&\text{if  $r\le 2\epsilon-2$.} 
\endcases
$$

First consider $\epsilon=0$. In this case, $$\binom{A}{D-1}_{\#p}=\cases 
\left(pd\binom ad\right)_{\#p}&\text{if  $ r\le p-2$}\\
\binom a{d-1}_{\#p}&\text{if  $p-1= r$}\\
\endcases$$ and $$\binom{3D-A-2}D_{\#p}=\cases
\binom{3d-a-2}d_{\#p}&\text{if $p-1\le r$}   
\\
\binom{3d-a-1}d_{\#p} &\text{if $r\le p-2$}; 
\endcases
$$so,
$$\left(\binom{A}{D-1}_{\#p}\binom{3D-A-2}D\right)_{\#p}=\cases
\binom a{d-1}_{\#p}\binom{3d-a-2}d_{\#p}&\text{if $p-1\le r$}   
\\
\left(pd\binom ad\right)_{\#p}\binom{3d-a-1}d_{\#p} &\text{if $r\le p-2$}, 
\endcases
$$which is at least $\binom{2d}d_{\#p}=\binom{2D}D_{\#p}$. Now 
consider $1\le \epsilon$. We have $$\binom{A}{D-1}_{\#p}=\cases 
\left(p(a-d)\binom{a}d\right)_{\#p}&\text{if $r\le \epsilon-2$}\\
\binom ad_{\#p}& \text{if $\epsilon-1\le r$}\\
\endcases$$ and (since $3\epsilon+p-1\le r\le p-1$ is impossible) $\binom{3D-A-2}D_{\#p}$ is equal to 
$$\cases
\left(p(2d-a-1)\binom{3d-a-1}d\right)_{\#p}&\text{if  
$2\epsilon+p-1\le r\le 3\epsilon+p-2$} 
\\
\binom{3d-a-1}d_{\#p} &\text{if $3\epsilon-1\le r\le p-2+2\epsilon$} 
\\
\left(p(2d-a) \binom{3d-a}{d}\right)_{\#p}&\text{if $2\epsilon-1\le r\le 3\epsilon-2$} 
\\
\binom{3d-a}d_{\#p}&\text{if  $r\le 2\epsilon-2$.} 
\endcases
$$Of course, if $r\le \epsilon-2$, then $r\le 2\epsilon-2$; so, $\left(\binom{A}{D-1}\binom{3D-A-2}D\right)_{\#p}$
is equal to 
$$\eightpoint \cases
\left(\binom adp(2d-a-1)\binom{3d-a-1}d\right)_{\#p}&\text{if  
$2\epsilon+p-1\le r\le 3\epsilon+p-2$} 
\\
\left(\binom ad\binom{3d-a-1}d\right)_{\#p} &\text{if $3\epsilon-1\le r\le p-2+2\epsilon$} 
\\
\left(\binom adp(2d-a) \binom{3d-a}{d}\right)_{\#p}&\text{if $2\epsilon-1\le r\le 3\epsilon-2$} 
\\
\left(\binom ad\binom{3d-a}d\right)_{\#p}&\text{if  $\epsilon-1\le r\le 2\epsilon-2$} 
\\
\left(p(a-d)\binom{a}d\binom{3d-a}d\right)_{\#p}&\text{if  $r\le \epsilon-2$,} 
\endcases
$$which is at least $\binom {2d}d_{\#p}=\binom {2D}D_{\#p}$.

Finally, we study (d) with $\epsilon\le -1$. We have 
$$\binom{A}{D-1}_{\#p}=\cases 
\left(pd\binom ad\right)_{\#p}&\text{if  $r\le p-2+\epsilon$}\\
\binom a{d-1}_{\#p}&\text{if   $p-1+\epsilon\le r$}\\
\endcases$$and 
$$\binom{3D-A-2}D_{\#p}= \cases
\left(pd\binom{3d-a-2}d\right)_{\#p}&\text{if  
$2\epsilon+p-1\le r$} 
\\
\binom{3d-a-2}{d-1}_{\#p}&  \text{if $3\epsilon+p-1\le r\le 2\epsilon+p-2$  } 
\\ 
\left(pd\binom{3d-a-1}d\right)_{\#p}&\text{if $r\le p+3\epsilon-2$}.  
\endcases
$$Of course, $p-1+2\epsilon<p-1+\epsilon$ (because $\epsilon<0$). 
Thus, $\left(\binom{A}{D-1}\binom{3D-A-2}D\right)_{\#p}$ is equal to 
$$\eightpoint \cases
\left(\binom a{d-1}pd\binom{3d-a-2}d\right)_{\#p}&\text{$p-1+\epsilon\le r$} \\
\left(pd\binom adpd\binom{3d-a-2}d\right)_{\#p}&\text{if  
$2\epsilon+p-1\le r\le p-2+\epsilon$} 
\\
\left(pd\binom ad\binom{3d-a-2}{d-1}\right)_{\#p}&  \text{if $3\epsilon+p-1\le r\le 2\epsilon+p-2$  } 
\\ 
\left(pd\binom adpd\binom{3d-a-1}d\right)_{\#p}&\text{if $r\le p+3\epsilon-2$}.  
\endcases$$
In the second line, we use $d\binom{a}d=\binom a{d-1}(a-d+1)$. In the third line we use 
(\tref{!-also}). We conclude that  $$\left(\binom{A}{D-1}\binom{3D-A-2}D\right)_{\#p} \ge \left(pd\binom{2d}d\right)_{\#p}=\binom{2D}D_{\#p}. \qed$$ \enddemo

\proclaim{Lemma \tnum{L8}}If $A=3a+r$ with $0\le r\le 2$ and $D=3d+\epsilon$ with $\epsilon$ equal to $0$ or $2$, then
\smallskip\flushpar{\bf(a)} $\binom AD_{\#3}=\cases \binom ad_{\#3}&\text{if $\epsilon=0$ or $\epsilon=r=2$}\\
\left(3(d+1)\binom{a}{d+1}\right)_{\#3}&\text{if $\epsilon=2$ and  $0\le r\le 1$},\endcases$

\smallskip\flushpar{\bf(b)} $\binom {2D}D_{\#3}=\cases \binom {2d}d_{\#3}&\text{if $\epsilon=0$}\\
\left(3(d+1)\binom{2d+1}{d}\right)_{\#3}&\text{if $\epsilon=2$},\endcases$

\smallskip\flushpar{\bf(c)} $\binom {2D+1}D_{\#3}=\cases \binom {2d}d_{\#3}&\text{if $\epsilon=0$}\\
\binom{2d+1}{d}_{\#3}&\text{if $\epsilon=2$},\endcases$

\smallskip\flushpar{\bf(d)} $\binom {3D}D_{\#3}=\cases \binom {3d}d_{\#3}&\text{if $\epsilon=0$}\\
\left(3(d+1)\binom{3d+2}{d}\right)_{\#3}&\text{if $\epsilon=2$,}\endcases$

\smallskip\flushpar{\bf(e)} $\binom {3D+2}D_{\#3}=\cases \binom {3d}d_{\#3}&\text{if $\epsilon=0$}\\
\binom{3d+2}{d}_{\#3}&\text{if $\epsilon=2$,}\endcases$ 

\smallskip\flushpar{\bf(f)} $\binom {3D-1-A}D_{\#3}=\cases 
\binom {3d-1-a}d_{\#3}&\text{if $\epsilon=0$}\\
\binom {3d+1-a}d_{\#3}&\text{if $\epsilon=2$ and $r=0$}\\
\left(3(2d-a+1)\binom{3d+1-a}{d}\right)_{\#3}&\text{if $\epsilon=2$ and $1\le r\le 2$,}\endcases$ 

\smallskip\flushpar{\bf(g)} $\binom {3D-A}D_{\#3}=\cases 
\binom {3d-a}d_{\#3}&\text{if $\epsilon=r=0$}\\
\binom {3d-1-a}d_{\#3}&\text{if $\epsilon=0$ and $1\le r\le 2$}\\
\left(3(3d-a+2)\binom{3d+1-a}{d}\right)_{\#3}&\text{if $\epsilon=2$ and $r=0$}\\
\binom {3d+1-a}d_{\#3}&\text{if $\epsilon=2$ and $r=1$}\\
\left(3(d+1)\binom{3d+1-a}{d+1}\right)_{\#3}&\text{if $\epsilon=r=2$, and}
\endcases$

\smallskip\flushpar{\bf(h)} $\binom {3D+1-A}{D+1}_{\#3}=\cases 
\binom {3d-a}d_{\#3}&\text{if $\epsilon=r=0$}\\
\left(3(2d-a)\binom {3d-a}d\right)_{\#3}&\text{if $\epsilon=0$ and $r=1$}\\
\binom {3d-1-a}d_{\#3}&\text{if $\epsilon=0$ and $r= 2$}\\
\binom {3d+2-a}{d+1}_{\#3}&\text{if $\epsilon=2$ and $0\le r\le 1$}\\
\binom{3d+1-a}{d+1}_{\#3}&\text{if $\epsilon=r=2$}.\endcases$ 
 
\endproclaim

\demo{Proof} We compute
$$ \binom AD_{\#3}=\left(\frac{A\cdots(A-D+1)}{D!}\right)_{\#3}=\left(\frac{[3a+r]\cdots[3(a-d)+(r+1-\epsilon)]}{(3d+\epsilon)!}\right)_{\#3}$$$$=\cases \left(\frac{[3a]\cdots[3(a-d)]}{[3d]\cdots [3(1)]}\right)_{\#3}=\left(3(d+1)\binom a{d+1}\right)_{\#3}&\text{if $r+1\le \epsilon$}\\
\left(\frac{[3a]\cdots[3(a-d+1)]}{[3d]\cdots [3(1)]}\right)_{\#3}= \binom a{d}_{\#3}&\text{if $\epsilon\le r$,}\endcases$$ and (a) is established. One may obtain (b) -- (g) from (a) or by direct calculation. We prove (h). Observe that
$$\binom {3D+1-A}{D+1}_{\#3}=\left(\frac{[3(3d+\epsilon-a)+1-r]\cdots [3(2d-a)+(2\epsilon-r+1)]}{(3d+\epsilon+1)!}\right)_{\#3}.$$
We see that  the largest multiple of $3$ in the numerator is 
$$\cases 3(3d+\epsilon-a-1)&\text{if $r=2$}\\3(3d+\epsilon-a)&\text{if $0\le r\le 1$.}\endcases$$We see that
$-1\le 2\epsilon-r+1\le 5$ and the smallest multiple of $3$ in the numerator is 
 $$\cases 3(2d-a)&\text{if $2\epsilon-r+1\le 0$}\\3(2d-a+1)&\text{if $1\le 2\epsilon-r+1\le 3$}\\
3(2d-a+2)&\text{if $4\le 2\epsilon-r+1$.}
\endcases$$The largest multiple of $3$ in the denominator is 
$$\cases 3(d+1)&\text{if $\epsilon=2$}\\3d&\text{if $\epsilon=0$}.\endcases$$It is not difficult to put the pieces together. \qed\enddemo

\proclaim{Proposition \tnum{P8-13}} Assertion {\rm(1)} of Theorem {\rm\tref{4}} holds.\endproclaim
\demo{Proof} It is clear that all five statements hold for $d=0$. The proof proceeds by induction. Let $D>0$ be an element of $D_3$. So $D=3d+\epsilon$ with $d\in D_3$, $d<D$, and $\epsilon$ equal to $0$ or $2$. The induction hypothesis ensures that all five statements hold for $d$. Lemma   \tref{L8}, together with the induction hypothesis, yields   $$\eightpoint \binom {2D}D_{\#3}=\left.\cases \binom {2d}d_{\#3}&\text{if $\epsilon=0$}\\
\left(3(d+1)\binom{2d+1}{d}\right)_{\#3}&\text{if $\epsilon=2$}\endcases \right\}\hskip-1.15pt=
\hskip-1.15pt\left.\cases \binom {3d}d_{\#3}&\text{if $\epsilon=0$}\\
\left(3(d+1)\binom{3d+2}{d}\right)_{\#3}&\text{if $\epsilon=2$}\endcases\right\}=\binom {3D}D_{\#3}$$ 
and$$\eightpoint \left.\binom {2D+1}D_{\#3}=\cases \binom {2d}d_{\#3}&\text{if $\epsilon=0$}\\
\binom{2d+1}{d}_{\#3}&\text{if $\epsilon=2$}\endcases\right\}=\left.\cases \binom {3d}d_{\#3}&\text{if $\epsilon=0$}\\
\binom{3d+2}{d}_{\#3}&\text{if $\epsilon=2$}\endcases\right\}= \binom {3D+2}D_{\#3}.$$ Statements (a) and (b) hold at $D$. For (c), use Lemma \tref{L8} and the induction hypothesis to see that 
$$\eightpoint\split  \left(\binom AD\binom{3D-A}D\right)_{\#3}&{}\ge  \cases 
\left(\binom ad\binom{3d-a}d\right)_{\#3}&\text{if $\epsilon=r=0$}\\
\left(\binom ad\binom{3d-1-a}d\right)_{\#3}&\text{if $\epsilon=0$ and $1\le r$}\\
\left(3(d+1)\binom a{d+1}\binom{3d+1-a}d\right)_{\#3}&\text{if $\epsilon=2$ and $r\le 1$}\\
\left(3(d+1)\binom ad\binom{3d+1-a}{d+1}\right)_{\#3}&\text{if $\epsilon=r=2$}\endcases \\&{}\ge 
\left.\cases \binom{2d}d_{\#3}&\text{if $\epsilon=0$}\\\left(3(d+1)\binom{2d+1}d\right)_{\#3}&\text{if $\epsilon=2$}\endcases\right\}=\binom{2D}D_{\#3}.\endsplit $$ Use Lemma \tref{L8}, the fact that
$$(2d-a+1)\binom{3d+1-a}d=(d+1)\binom{3d+1-a}{d+1},$$ and the induction hypothesis to see that
$$\eightpoint\split \left(\binom AD\binom{3D-1-A}D\right)_{\#3}&{}\ge  
\cases\left(\binom ad\binom {3d-1-a}d\right)_{\#3}&\text{if $\epsilon=0$}\\
\left(3(d+1)\binom a{d+1}\binom {3d+1-a}d\right)_{\#3}&\text{if $\epsilon=2$ and $r\le 1$}\\
\left(\binom a{d}3(d+1)\binom {3d+1-a}{d+1}\right)_{\#3}&\text{if $\epsilon=r=2$}\endcases\\&{}\ge 
\left.\cases \binom{2d}d_{\#3}&\text{if $\epsilon=0$}\\\left(3(d+1)\binom{2d+1}d\right)_{\#3}&\text{if $\epsilon=2$}\endcases\right\}=\binom{2D}D_{\#3};\endsplit $$thereby verifying (d) at $D$. Finally, we consider (e). In addition to Lemma \tref{L8} and the induction hypothesis, we use
$$\eightpoint(d+1)\binom a{d+1}=(a-d)\binom ad\quad\text{and}\quad (d+1)\binom{3d+2-a}{d+1}=(3d+2-a)\binom{3d+1-a}d,$$ to see that
$$\eightpoint\split \left(\binom AD\binom{3D+1-A}{D+1}\right)_{\#3}&{}\ge  
\cases\left(\binom ad\binom {3d-a}d\right)_{\#3}&\text{if $\epsilon=0$ and $r\le 1$}\\
\left(3(d+1)\binom a{d}\binom {3d-1-a}d\right)_{\#3}&\text{if $\epsilon=0$ and $r=2$}\\
\left(\binom a{d+1}\binom {3d+1-a}{d}\right)_{\#3}&\text{if $\epsilon=2$ and $r\le 1$}\\
\left(\binom a{d}\binom {3d+1-a}{d+1}\right)_{\#3}&\text{if $\epsilon=r=2$}\\
\endcases\\&{}\ge 
\left.\cases \binom{2d}d_{\#3}&\text{if $\epsilon=0$}\\\binom{2d+1}d_{\#3}&\text{if $\epsilon=2$}\endcases\right\}=\binom{2D+1}D_{\#3}. \qed\endsplit $$
\enddemo

\SectionNumber=5\tNumber=1
\heading Section \number\SectionNumber. \quad The projective dimension of $Q_{\pmb k,n,N}$ is finite when $\lfloor\frac Nn\rfloor$ is in $T_p$.
\endheading

Throughout this section, $\pmb k$ is a field, $n$ is a positive integer, $R$ is the diagonal hypersurface ring $R_{\pmb k,n}$ of {\rm(\tref{DHSR})}, and $B=\pmb k[X,Y]$ is a polynomial ring. 

\definition{Notation \tnum{A1}} For each positive integer $a$,
we let $\operatorname{HB}_a$ denote a homogeneous  Hilbert-Burch matrix for the row vector  $\bmatrix X^a& Y^a& (X+Y)^a\endbmatrix$ in $B=\pmb k[X, Y]$. In particular, 
$$0\to B(-a-d_1)\oplus B(-a-d_2)@> \operatorname{HB}_a >> B(-a)^3@> \bmatrix X^a& Y^a& (X+Y)^a\endbmatrix >>B$$ is an exact sequence of homogeneous $B$-module homomorphisms and if 
$$\operatorname{HB}_a=\bmatrix H_{1,1}&H_{1,2}\\H_{2,1}&H_{2,2}\\H_{3,1}&H_{3,2}\endbmatrix,$$ then
$H_{i,j}$ is a homogeneous form in $B$ of degree $d_j$; furthermore the signed maximal order minors of $\operatorname{HB}_a$ are
$$\matrix\format \r&\ \l \\X^a= \det \bmatrix H_{2,1}&H_{2,2}\\H_{3,1}&H_{3,2}\endbmatrix, 
&Y^a= -\det \bmatrix H_{1,1}&H_{1,2}\\H_{3,1}&H_{3,2}\endbmatrix,
\ \text{and}\\  (X+Y)^a= \det \bmatrix H_{1,1}&H_{1,2}\\H_{2,1}&H_{2,2}\endbmatrix.\endmatrix\tag\tnum{smom}$$ Each of the homogeneous forms $H_{i,j}$ has degree $d_j$. We say that  column $j$ of $\operatorname{HB}_a$ is a relation on $\bmatrix X^a& Y^a& (X+Y)^a\endbmatrix$ of degree $d_j$. A quick look at (\tref{smom}) shows that the sum of the degrees $d_1+d_2$ must equal $a$. If $|d_1-d_2|\le 1$, then the Hilbert-Burch matrix $\operatorname{HB}_a$ is called {\it balanced}; otherwise, $\operatorname{HB}_a$ is called {\it unbalanced}. 

We notice that the Hilbert-Burch matrix $\operatorname{HB}_a$ continues to make sense even if $a=p^e$, where $p>0$ is the characteristic of $\pmb k$ and $e\ge 1$ is an integer. In this case, 
$$\operatorname{HB}_a=\bmatrix 1&-Y^a\\1&X^a\\-1&0\endbmatrix,$$ the formulas of (\tref{smom}) hold, and 
$\operatorname{HB}_a$ is unbalanced. 

The Hilbert-Burch matrix $\operatorname{HB}_a$ gives information which is relevant to the $R_{\pmb k,n}$ modules $Q_{\pmb k,n,N}$ by way of the $\pmb k$-algebra homomorphism $\alpha\:\pmb k[X,Y]\to R_{\pmb k,n}$. The map $\alpha$ is  defined by $\alpha(X)=x^n$ and $\alpha(Y)=y^n$. We see that 
$$\alpha(X+Y)=x^n+y^n=-z^n.$$\enddefinition

\proclaim{Observation \tnum{O5.1}} Let $\pmb k$ be a field and $n$ and $N$ be positive integers. If $n$ divides $N$, then $Q_{\pmb k,n,N}$ has finite projective dimension as an $R_{\pmb k,n}$-module. \endproclaim

\demo{Proof} Let $N=an$ and  let $\operatorname{HB}_a$ be a Hilbert-Burch matrix for $[X^a,Y^a,(X+Y)^a]$ over $\pmb k$. If $\alpha\:\pmb k[X,Y]\to R_{\pmb k,n}$ is the $\pmb k$-algebra homomorphism    defined by $\alpha(X)=x^n$ and $\alpha(Y)=y^n$, then the ideal generated by the maximal order minors of $\alpha(\operatorname{HB}_a)$ is $(X^N,Y^N,Z^N)$ and 
$$0\to R_{\pmb k,n}^2@> \alpha(\operatorname{HB}_a)>> R_{\pmb k,n}^3 @>>> R_{\pmb k,n} $$ is a finite free resolution of $Q_{\pmb k,n,N}$. \qed \enddemo

\proclaim{Lemma \tnum{L-8-5}} Fix a positive integer $a$ and a field $\pmb k$. Suppose that some minimal, non-zero, homogeneous relation on $[X^a,Y^a,(X+Y)^a]$ in $\pmb k[X,Y]$ has the form 
$$\xi=\bmatrix X\xi_1'\\Y\xi_2'\\(X+Y)\xi_3'\endbmatrix,$$ 
where the $\xi_i'$ are homogeneous elements of $\pmb k[X,Y]$. If $N$ and $n$ are positive integers with $N=an+r$ for some $r$ with $0\le r\le n$, then $\operatorname{pd}_{R_{\pmb k,n}}Q_{\pmb k,n,N}$ is finite.\endproclaim
\demo{Proof} 
We complete $\xi$ to form a homogeneous Hilbert-Burch matrix $$\operatorname{HB}_a=\left[\matrix X\xi_1'&\eta_1\\ Y\xi_2'&\eta_2\\(X+Y)\xi_3'&\eta_3\endmatrix\right]$$ for $[X^a,Y^a,(X+Y)^a]$ in $\pmb k[X,Y]$. 
The signed maximal order minors of $\operatorname{HB}_a$ in the sense of (\tref{smom}) are $X^a$, $Y^a$, and $(X+Y)^a$. Apply the $\pmb k$-algebra homomorphism $\alpha\:\pmb k[X,Y]\to R_{\pmb k,n}$ with $\alpha(X)=x^n$ and $\alpha(Y)=y^n$. We see that the ideal generated by   the maximal order minors of
$$\operatorname{hb}=\bmatrix x^{n-r}\alpha(\xi_1')&(yz)^r\alpha(\eta_1)\\
y^{n-r}\alpha(\xi_2')&(xz)^r\alpha(\eta_2)\\
-z^{n-r}\alpha(\xi_3')&(xy)^r\alpha(\eta_3)\endbmatrix$$in $R_{\pmb k,n}$ is generated by $(x^N,y^N,z^N)$; and therefore
$$0\to R_{\pmb k,n}^2@> \operatorname{hb} >> R_{\pmb k,n}^3\to R_{\pmb k,n}$$ is a finite free resolution of $Q_{\pmb k,n,N}$. The entries of $\operatorname{hb}$ are in $R_{\pmb k,n}$. The only constraint on $r$ is that $r$ and $n-r$ both must be non-negative integers. \qed
\enddemo

\proclaim{Lemma \tnum{s-f}} Let $\pmb k$ be a field, $n$ and $a$ be  positive integers,  and   $R$ be the diagonal hypersurface ring $R_{\pmb k,n}$ of {\rm(\tref{DHSR})}. 
\flushpar{\bf (1)} If some Hilbert-Burch matrix $\operatorname{HB}_a$ has the form
$$
\operatorname{HB}_a=\left[\matrix\format \c\qquad &\c\\ F & XI \\ YG & J \\ (X+Y)H & K\\ \endmatrix\right],$$ for   homogeneous forms $F,G,H,I,J,K$ in $B=\pmb k[X,Y]$, then the $R$-module $Q_{\pmb k,n,N}$ has finite  projective dimension for all $N=na+r$, with $0 \le r \le n$.

\smallskip\flushpar{\bf (2)} If some Hilbert-Burch matrix $\operatorname{HB}_a$ has the form
$$\operatorname{HB}_a=\left[ \matrix\format \c\qquad &\c\\  Y(X+Y) F & I \\ (X+Y)G & XJ \\ YH & XK \\ \endmatrix\right],
$$for   homogeneous forms $F,G,H,I,J,K$ in $B=\pmb k[X,Y]$,  then the $R$-module $Q_{\pmb k,n,N}$ has finite projective dimension for all $N=na-r$, with $0 \le r \le n$.
\endproclaim

\demo{Proof}
For the first assertion,   the maximal order minors of the matrix 
$$
\operatorname{hb}=\left[ \matrix\format \r\qquad &\r\\ \alpha(F) & y^rz^rx^{n-r}\alpha(I) \\ x^ry^{n-r}\alpha(G) & z^r\alpha(J) \\ -x^rz^{n-r}\alpha(H) & y^r\alpha(K)\\ \endmatrix\right]
$$ generate the ideal $(x^{na+r}, y^{na+r}, z^{na+r})$ of $R$. 
For the second assertion, the maximal order minors of the matrix
$$
\operatorname{hb}=\left[\matrix\format \r\qquad &\r\\ x^ry^{n-r}z^{n-r}\alpha(F) & \alpha(I) \\
z^{n-r}\alpha(G) & x^{n-r}y^r\alpha(J) \\ y^{n-r}\alpha(H) & -x^{n-r}z^r\alpha(K) \\ \endmatrix\right]$$
generate the ideal   $(x^{na-r}, y^{na-r}, z^{na-r})$ of $R$. In each case
$$0\to R^3@> \operatorname{hb}>> R^2@>   >> R$$ is a resolution of $Q_{\pmb k,n,N}$ by free $R$-modules.
\qed \enddemo

\proclaim {Theorem \tnum{AVT}} Consider the data $(\pmb k,n,N)$ with $N=an+r$ and $0\le r\le n$. If $\operatorname{HB}_a$ is unbalanced over $\pmb k$, then  the $R_{\pmb k,n}$-module $Q_{\pmb k,n,N}$ has finite projective dimension.
\endproclaim
\demo{Proof}Let 
$$
\operatorname{HB}_a=\left[ \matrix\format \c\qquad &\c\\ F & I \\G & J \\ H & K \\ \endmatrix\right],
\qquad
\operatorname{HB}_{a+1}=\left[\matrix\format \c\qquad &\c\\ M & Q \\N& R \\ P & S \\ \endmatrix\right],$$
$d_a=\deg(F)$, $D_a=\deg(I)$, $d_{a+1}=\deg(M)$, and  $D_{a+1}=\deg(Q)$. The hypothesis that $\operatorname{HB}_a$ is unbalanced guarantees that $|D_a-d_a|\ge 2$.
Note that $$[XM, YN, (X+Y)P]^{\text{\rm T}}$$ is a relation on $[X^a, Y^a, (X+Y)^a]$, and thus we have
$$ 
\left[\matrix\format \c  \\ XM \\ YN \\ (X+Y)P \\ \endmatrix\right]=f_1\left[\matrix\format \c  \\ F \\ G \\ H \\ \endmatrix\right] + g_1 \left[ \matrix\format \c \\ I \\ J \\ K \\ \endmatrix\right]
\tag\tnum{eq1}$$
for some $f_1, g_1 \in \pmb k[X, Y]$, which can be taken to be homogeneous.
Similarly we have
$$ 
\left[ \matrix\format \c \\ XQ \\ YR \\ (X+Y)S\\ \endmatrix\right] = f_2 \left[ \matrix\format \c \\ F \\ G \\ H \endmatrix\right] +g_2\left[ \matrix\format \c \\ I \\ J \\ K \\ \endmatrix\right],
\tag\tnum{eq2}$$for some homogeneous $f_2, g_2 \in \pmb k[X, Y]$. 
Using equations (\tref{eq1}) and (\tref{eq2}) to calculate the $2\times 2$ minors of the matrix $\operatorname{HB}_{a+1}$ yields
$$ 
f_1g_2-g_1f_2=XY(X+Y).
\tag\tnum{eq3}$$
Comparing degrees in (\tref{eq1}) and (\tref{eq2}) yields
$$\deg(f_1)-\deg(g_1)=\deg(f_2)-\deg(g_2)=D_a-d_a.\tag\tnum{refo}$$

We first consider the case where all four polynomials   $f_1, g_1, f_2, g_2$ are nonozero. Then (\tref{eq3}) implies $\deg(f_1)+\deg(g_2)=\deg(g_1)+\deg(f_2)=3$.
There are two possibilities: either one of $f_1, f_2, g_1, g_2$ is a unit, in which case one can modify $\operatorname{HB}_a$ so that it has a column of the form $[XM, YN, (X+Y)P]^{\text{\rm T}}$ (which gives the conclusion by Lemma \tref{L-8-5}), or $\deg(f_1), \deg(g_1), \deg(f_2), \deg(g_2) \in \{1, 2\}$, which implies that $|\deg(F)-\deg(I)|=0$ or $1$, contradicting the hypothesis.

It remains to consider the case when one of $f_1, g_1, f_2, g_2$ is zero. Without loss of generality, assume $g_1=0$, and none  of the polynomials $f_1$, $f_2$, or  $g_2$ is a unit.
Then (\tref{eq3}) becomes $f_1g_2=XY(X+Y)$. Since $X$, $Y$, and $X+Y$ are irreducible in the Unique Factorization Domain $\pmb k[X, Y]$, we may assume without loss of generality that one of the following holds:
\roster
\item case 1: $f_1=X, g_2=Y(X+Y)$
\item case 2: $f_1=Y(X+Y), g_2=X$
\endroster
The hypothesis $|\deg(F)-\deg(I)|\ge 2$, together with the equation (\tref{refo}), ensures that $|\deg(f_2)-\deg(g_2)|\ge 2$. We have $\deg g_2\in\{1,2\}$ and $f_2$ is not a unit; so either $f_2$ is zero or $\deg f_2\ge 3$. Write $f_2=Xu+cY^s$ for some homogeneous $u\in \pmb k[X, Y]$ and $c\in \pmb k$, with $\deg f_2=\deg u+1=s$. We have seen that if $c\neq 0$, then $s\ge 3$.

Consider case 1.
Equations (\tref{eq1}) and (\tref{eq2}) give $Y|G$ and $(X+Y)|H$ in $\pmb k[X,Y]$. We will apply part (1) of  Lemma \tref{s-f} to some appropriate modification of $\operatorname{HB}_a$.  Equation (\tref{eq2}) tells us that
$$XQ=(Xu+cY^s)F+Y(X+Y)I.$$Thus,
$$X(Q-uF-YI)=Y^2(cY^{s-2}F+I)$$ and $X$ divides $cY^{s-2}F+I$ in $\pmb k[X,Y]$. (Recall that the expression $cY^{s-2}$ is meaningful in $\pmb k[X,Y]$ because either $c=0$ or $s\ge 3$.) We replace the column $[I, J, K]^{\text{\rm T}}$ in $\operatorname{HB}_a$ by 
$[I, J, K]^{\text{\rm T}}+cY^{s-2}[F, G, H]^{\text{\rm T}}$. Note that this new matrix can be used in  place of $\operatorname{HB}_a$, and it has the form specified in the first part of Lemma \tref{s-f}. Thus, the conclusion follows from Lemma~\tref{s-f}.

Now consider case 2.
Equation (\tref{eq1}) implies that $M$ is divisible by $Y(X+Y)$, $N$ is divisible by $X+Y$, $P$ is divisible by $Y$. We will apply part (2) of  Lemma \tref{s-f} to some appropriate modification of $\operatorname{HB}_{a+1}$. Recall that (\tref{eq2}) gives
$$\bmatrix XQ\\YR\\(X+Y)S\endbmatrix =(cY^s+uX)\bmatrix F\\G\\H\endbmatrix +X\bmatrix I\\J\\K\endbmatrix $$ and (\tref{eq1}) gives
$$\bmatrix XM\\YN\\(X+Y)P\endbmatrix =(Y^2+XY)\bmatrix F\\G\\H\endbmatrix.$$ 
We see that 
(\tref{eq2}) minus $cY^{s-2}$ times (\tref{eq1}) is
$$\eightpoint \bmatrix X&0&0\\0&Y&0\\0&0&X+Y\endbmatrix\left(\bmatrix Q\\R\\S\endbmatrix -cY^{s-2}\bmatrix M\\N\\P\endbmatrix\right)=X\left(u\bmatrix F\\G\\H\endbmatrix +\bmatrix I\\J\\K\endbmatrix -cY^{s-1}\bmatrix F\\G\\H\endbmatrix\right).$$ Thus, the bottom two entries of 
$$\bmatrix Q\\R\\S\endbmatrix -cY^{s-2}\bmatrix M\\N\\P\endbmatrix$$ are divisible by $X$.
We may replace the column $[Q, R, S]^{\text{\rm T}}$ in $\operatorname{HB}_{a+1}$ by $[Q,R, S]^{\text{\rm T}}-cY^{s-2}[M, N, P]^{\text{\rm T}}$. This new matrix can be used in the role of $\operatorname{HB}_{a+1}$, and it has the form specified in the second part of Lemma \tref{s-f}. Thus, the conclusion follows from Lemma \tref{s-f}.
\qed \enddemo

\proclaim{Theorem \tnum{T-BK} (Brenner and Kaid)}Let $\pmb k$ be a field of positive characteristic $p$ and let $a$ be a positive integer.
\roster
\item If $a$ is odd, then $\operatorname{HB}_a$ is unbalanced over $\pmb k$ if and only if there exists an odd integer J and a power $q=p^e$ of $p$ such that $|a-Jq|<\frac{q-1}3$.
\item If $a$ is even, then $\operatorname{HB}_a$ is unbalanced over $\pmb k$ if and only if there exists an odd integer J and a power $q=p^e$ of $p$ such that $|a-Jq|<\frac{q}3$.\endroster\endproclaim
\demo{Proof}Let $\bar{ \pmb k}$ be the algebraic closure of $\pmb k$ and let $A$ be the $\bar{ \pmb k}$-algebra $$A=\bar{ \pmb k}[X,Y,Z]/(X^a,Y^a,Z^a).$$ It is shown in \cite{\rref{BK}, Cor.~2.2}  that $A$ has the weak Lefschetz property (WLP) if and only if $\operatorname{HB}_a$ over $\bar{\pmb k}$ is balanced. In \cite{\rref{BK}, Thm.~2.6} numerical conditions, which characterize the set of $a$ for which $A$ does not have the WLP, are given. We have reformulated these numerical conditions. Notice that
$$\frac {3a+\epsilon}{6k+4}< q< \frac{3a-\epsilon}{6k+2}\iff (2k+1)q-\frac{q-\epsilon}3<a<(2k+1)q+\frac{q-\epsilon}3.$$
We notice that $\bar{\pmb k}[X,Y]$ is a free $\pmb k[X,Y]$-module. If $\Bbb F$ is a resolution of $$\pmb k[X,Y]/(X^a,Y^a,(X+Y)^a)$$ by free $\pmb k[X,Y]$-modules, then $\Bbb F\otimes_{\pmb k[X,Y]}\bar{\pmb k}[X,Y]$ is a resolution of $$\bar{\pmb k}[X,Y]/(X^a,Y^a,(X+Y)^a)$$ by free $\bar{\pmb k}[X,Y]$-modules. In particular, if $\operatorname{HB}_a$ is a Hilbert-Burch matrix for the matrix $[X^a,Y^a,(X+Y)^a]$ over $\pmb k$, then $\operatorname{HB}_a$ is also a Hilbert-Burch matrix for $[X^a,Y^a,(X+Y)^a]$ over $\bar{\pmb k}$. In other words, $\operatorname{HB}_a$ is balanced over $\pmb k$ if and only if $\operatorname{HB}_a$ is balanced over $\bar{\pmb k}$. \qed\enddemo

\proclaim{Lemma \tnum{L-8-5.2}}Fix a non-negative integer $a$ and a field $\pmb k$. Let $n$ and $N$ be positive integers with $N=an+r$ and  $0\le r\le n$. Suppose that 
$$\xi=\bmatrix X^j\xi_1'\\Y^j\xi_2'\\(X+Y)^j\xi_3'\endbmatrix$$ is a  non-zero, homogeneous relation on $[X^a,Y^a,(X+Y)^a]$ in $\pmb k[X,Y]$ for some homogeneous polynomials $\xi_i'\in \pmb k[X,Y]$.
\roster
\item If $j=1$ and $\deg \xi\le \lfloor \frac a2\rfloor$, then  $\operatorname{pd}_{R_{\pmb k,n}}Q_{\pmb k,n,N}$ is finite.
\item If $j=2$ and $\deg \xi\le \lceil \frac a2\rceil$,  then  $\operatorname{pd}_{R_{\pmb k,n}}Q_{\pmb k,n,N}$ is finite.\endroster
\endproclaim
\demo{Proof} (1)  If $\xi$ is a minimal relation, then use Lemma \tref{L-8-5}. If $\xi$ is not a minimal relation, then $\operatorname{HB}_a$ is unbalanced; so use Theorem \tref{AVT}.

(2) If $\xi$ is a minimal relation on $\rho_a=[X^a,Y^a,(X+Y)^a]$, then the conclusion follows from Lemma \tref{L-8-5}. Henceforth, we assume that $\xi$ is not a minimal relation. The  Hilbert-Burch Theorem guarantees that some minimal   relation on $\rho_a$ has degree at least $\lceil \frac a2\rceil$. Thus, $\xi$ is a multiple of some minimal relation on $\rho_a$.
In other words, $\xi=A\chi$ for some homogeneous $A$ in $\pmb k[X,Y]$ and some minimal relation $\chi$. If degree $A=1$, then the conclusion follows from Lemma \tref{L-8-5} applied to $\chi$. If $\deg A\ge 2$, then $\operatorname{HB}_a$ is unbalanced and the conclusion follows from Theorem \tref{AVT}. 
\qed 
\enddemo

\proclaim{Observation \tnum{A-5}} Let $J$ be an odd integer, $\pmb k$ be a field of characteristic $p\ge 5$, $q=p^e$ for some integer $e\ge 1$, and $N=an+r$ for $0\le r\le n$.

\roster\item If $q\equiv 2\mod 3$ and   $a=Jq-\frac{q+1}3$, then $\operatorname{pd}_{R_{\pmb k,n}}Q_{\pmb k,n,N}$ is finite. 
\item If $q\equiv 1\mod 3$ and   $a=Jq-\frac{q-1}3$, then $\operatorname{pd}_{R_{\pmb k,n}}Q_{\pmb k,n,N}$ is finite. \endroster
\endproclaim
\demo{Proof} (1) Write $J$ as $2k+1$. Select a relation $\bmatrix A,B,C\endbmatrix^{\text{\rm T}}$ on $[X^J,Y^J,(X+Y)^J]$ of degree $\le k$. Observe that $$\xi= \bmatrix X^{\frac{q+1}3}A^q\\Y^{\frac{q+1}3}B^q\\(X+Y)^{\frac{q+1}3}C^q\endbmatrix$$ is a relation on $[X^a,Y^a,(X+Y)^a]$ of degree at most $\frac{q+1}3+qk$. 
We see that $$a=qJ-\frac{q+1}3=q(2k+1)-\frac{q+1}3=2qk+\frac{2q-1}3=2\left(kq+\frac{q-2}3\right)+1.$$
So, $\deg \xi\le \lceil \frac a2\rceil$. Apply Lemma \tref{L-8-5.2}. Notice that $\frac{q+1}3\ge 2$ since $q\ge 5$. 

(2) Observe that  $$\xi= \bmatrix X^{\frac{q-1}3}A^q\\Y^{\frac{q-1}3}B^q\\(X+Y)^{\frac{q-1}3}C^q\endbmatrix$$ is a relation on $[X^a,Y^a,(X+Y)^a]$ of degree at most $\frac{q-1}3+qk=\lfloor \frac a2\rfloor$. 
\qed \enddemo

\proclaim{Theorem \tnum{HW672}} Let $\pmb k$ be a field of characteristic $c$ and let $n$ and $N$ be positive integers. Write   $N=an+r$ with $0\le r\le n$. 
If   $a\in T_c$,  then  the $R_{\pmb k,n}$-module $Q_{\pmb k,n,N}$ has finite projective dimension.\endproclaim

\demo{Proof} The set $T_0$ is empty; so $c=p$ is positive.  Lemma \tref{T24.5} establishes the result for $p=2$. Throughout the rest of the proof we take $p\ge 3$.

Identify an odd integer $J$ and a power $q=p^e$ of $p$ so that $a$ is close to $Jq$ in the sense of Definition \tref{D5.1}.
First take $p=3$.
If   $Jq-\frac q3+1\le a\le Jq+\frac q3 -1$, then $|Jq-a|\le \frac q3-1<\frac{q-1}3$; so,
Theorem \tref{T-BK} yields that $\operatorname{HB}_a$ is unbalanced and the conclusion follows from Theorem \tref{AVT}. We now assume that $a=Jq-\frac{q}3$. Write $J=2k+1$. Let $[A,B,C]^{\text{\rm T}}$ be a relation on $[X^J,Y^J,(X+Y)^J]$ of degree at most $k$. Observe that 
$$\xi=\bmatrix A^qX^{q/3}\\B^qY^{q/3}\\ C^q(X+Y)^{q/3}\endbmatrix$$is a relation on $[X^a,Y^a,(X+Y)^a]$ of degree at most $kq+\frac q3=\frac a2$. If $\xi$ is a minimal relation on $[X^a,Y^a,(X+Y)^a]$, then $1\le \frac q3$ and Lemma \tref{L-8-5} yields the conclusion. Otherwise, $\operatorname{HB}_a$ is unbalanced at we may apply Theorem \tref{AVT}.

Henceforth, we consider $p\ge 5$.  First we assume that $q\equiv 2\mod 3$. If 
$$Jq-\frac{q+1}3+1\le a\le Jq+\frac{q-2}3,$$ then $|Jq-a|\le \frac{q-2}3$ and Theorem \tref{T-BK} shows that $\operatorname{HB}_a$ is unbalanced over $\pmb k$ and Theorem \tref{AVT} shows that $Q_{\pmb k,n,N}$ has finite projective dimension. If $a=Jq-\frac{q+1}3$, then apply Observation \tref{A-5}. Now we assume that $q\equiv 1\mod 3$. If $$Jq-\frac{q-1}3+1\le a\le Jq+\frac{q-4}3,$$ then $|Jq-a|\le \frac{q-4}3$ and the conclusion follows from Theorems \tref{T-BK} and \tref{AVT}. If $a=Jq-\frac{q-1}3$, then we apply Observation \tref{A-5}. 
\qed\enddemo

\proclaim{Lemma \tnum{T24.5}} Let $\pmb k$ be a field of characteristic $2$ and $n$ and $N$ be positive integers. If $n\le N$, then $Q_{\pmb k,n,N}$ has finite projective dimension as a module over $R_{\pmb k,n}$  \endproclaim
\demo{Proof} Write $N$ in the form  $N=qn+r$ for integers $q$ and $r$ with $q=2^e$ and $0\le r\le qn$. Consider the matrix $$\operatorname{hb}= \bmatrix y^rz^r&x^{N-2r}\\x^rz^r&y^{N-2r}\\x^ry^r&z^{N-2r} \endbmatrix.$$
Observe that the ideal generated by the maximal order minors of $\operatorname{hb}$ is $(x^N,y^N,z^N)$. It follows that
$$0\to R_{\pmb k,n}^2@> \operatorname{hb}>> R_{\pmb k,n}^3\to R_{\pmb k,n}$$ is a finite resolution of $Q_{\pmb k,n}$ by free $R_{\pmb k,n}$-modules.\qed \enddemo

\SectionNumber=6\tNumber=1
\heading Section \number\SectionNumber. \quad The set of non-negative integers may be partitioned as  $S_p\cup T_p$.
\endheading

\proclaim{Theorem \tnum{H72+}}  If $c$ is the characteristic of a field, then the set of non-negative integers is the disjoint union of $S_c$ and $T_c$.  \endproclaim

\demo{Proof} The assertion is obvious if $c$ is equal to $0$ or $2$. Henceforth, we take $c$ to be a prime integer $p\ge 3$.

  Let $a$ be a non-negative integer and $\pmb k$ be a field of characteristic $p$. Pick $n$ and $N$ so that $N=an+r$ with $1\le r\le n-1$. The projective dimension of $Q_{\pmb k,n,N}$ over $R_{\pmb k,n}$ is either finite or infinite. If $a\in S_p$, then Theorem \tref{I-MAIN} shows that $\operatorname{pd} Q_{\pmb k,n,N}$ is infinite. If $a\in T_p$, then Theorem \tref{HW672} shows that $\operatorname{pd} Q_{\pmb k,n,N}$ is finite. Thus the sets $S_p$ and $T_p$ are disjoint. 

Suppose now that $a$ is odd. We saw in Remark \tref{!R18.1}, Remark \tref{R1.1'}, and Observation \tref{also} 
that if $p=3$, then 
$$a\in S_p\iff a-1\in S_p\quad\text{and}\quad a\in T_p \iff a-1\in T_p;$$ and if $p\ge5$, then
$$a\in T_p\iff a+1\in T_p\quad\text{and}\quad a\in T_p \iff a+1\in T_p.$$ 
Thus, it suffices that prove that every non-negative even integer is in $S_p\cup T_p$ for each prime $p\ge 3$. 
 
Let $p\ge 5$ be a prime integer and let $m$ be  an even integer with $m \notin S_p$. We prove that $m\in T_p$.  
Write $m$ as $a_0 + a_1p + \ldots + a_t p^t$ with $a_0, \ldots, a_t$ digits in the sense of Notation \tref{N5.1}.  Since $m$ is even and we are assuming $m \notin S_p$, it follows that at least two of the digits, say $a_k$ and $a_\ell$ are odd and less that $p/3$. More precisely, assume that $k$ is the smallest index such that $a_k$ is odd, and there are an odd number of indices among $a_{k+1}, \ldots, a_t$ for which the corresponding digits are odd.

We can write
$$
m=p^{k+1}J + a_kp^k+a_{k-1}p^{k-1} + \ldots + a_0
$$
where $J=a_{k+1} + \ldots + a_t p^{t-k-1}$ is an odd integer. Let $q=p^{k+1}$. We obtain an upper bound for $|m-Jq|$. Let $u$ denote the largest odd integer which is less than $p/3$, and let $v$ denote the largest even integer which is less than $2p/3$.
Since $a_k, a_{k-1}, \ldots, a_0$ are digits, and $a_k$ is odd, we have
$$ \split 
|m-Jq|&{}=|a_kp^k + a_{k-1} p^{k-1} + \ldots + a_0|\\
&{}< u p^k + v(p^{k-1} + \ldots + 1)\\
&{}=
up^k + v\frac{p^k-1}{p-1}< \left( u+\frac{v}{p-1}\right)p^k.\endsplit
$$
 If $p\equiv 1\mod 3$, then $u=(p-4)/3$ and $v=(2p-2)/3$. If $p\equiv 2 \mod 3$, then $u=(p-2)/3$ and $v=(2p-4)/3$. A straightforward calculation establishes $$
u + \frac{v}{p-1} < \frac{p}{3} 
$$ in each case. We conclude that $|m-Jq|<(\frac{p}{3})p^k=\frac q3$; and therefore $m\in T_p$ by Observation \tref{also}.

  Now take  $p=3$. We prove  that every non-negative even integer is in $S_3\cup T_3$. We know that $0\in S_3$ and $2\in T_3$. The proof proceeds by induction. Suppose that $m\ge 4$ is an even integer with $m\notin T_3$. Suppose further that every non-negative even integer $a$ with $a<m$ is in $S_3\cup T_3$. 
We observe first that $m$ is not congruent to $2\mod 3$. Indeed, if $m\equiv 2 \mod 3$. Then $m+1$ is equal to $3J$ for some odd integer $J$; thus, $m+1\in T_3$ by the definition of $T_3$ and $m\in T_3$ by Remark \tref{R1.1'} and this is a contradiction. Thus, $m\equiv 0$ or $m\equiv 1\mod 3$. We treat the two cases separately. 

We first suppose that  $m\equiv 0 \mod 3$. In this case, we consider the even integer $a=\frac m3$. We know that $0\le a<m$; so the induction hypothesis shows that $a\in S_3\cup T_3$. On the other hand, it is not possible for $a  \in T_3$. Indeed if $a  \in T_3$, then there exists an odd integer $J$ and a power $q=3^e$ with $Jq-\frac q3\le a \le Jq +\frac q3 -1$. Multiply by $3$ to see that 
$$\tsize J(3q)-\frac{3q}3\le m\le J(3q)+\frac{3q}3-3<J(3q)+\frac{3q}3-1;$$ and therefore $m\in T_3$, which is a contradiction. We have  $a  \in S_3\cup T_3$ and $a\notin T_3$. It follows that $a$ is an even integer in $S_3$. The definition of $S_3$ shows that $m=3a\in S_3$.

Finally, we suppose that  $m\equiv 1 \mod 3$. In this case, we consider the even integer $a=\frac {m-4}3$. We know that $0\le a<m$; so the induction hypothesis shows that $a$ is in $S_3\cup T_3$. On the other hand, it is not possible for $a  \in T_3$. Indeed if $a  \in T_3$, then there exists an odd integer $J$ and a power $q=3^e$ with $Jq-\frac q3\le a \le Jq +\frac q3 -1$. Multiply by $3$ and add $4$ to see that 
$$\tsize J(3q)-\frac {3q}3<J(3q)-\frac {3q}3+4\le m\le J(3q)+\frac {3q}3+1.$$
The integer $m$ is even and $\equiv 1 \mod 3$; whereas $J(3q)+\frac {3q}3+1$ is odd and $\equiv 1 \mod 3$.
It follows that $m\le J(3q)+\frac {3q}3-2$. At any rate, we would have
$$\tsize J(3q)-\frac {3q}3\le m \le J(3q)+\frac {3q}3-2< J(3q)+\frac {3q}3-1;$$ and therefore $m\in T_3$, which is a contradiction. We have  $a  \in S_3\cup T_3$ and $a\notin T_3$. It follows that $a$ is an even integer in $S_3$. The definition of $S_3$ shows that $m$, which is equal to $3a+4$, is in $S_3$.
\qed \enddemo
 
\proclaim{Theorem \tnum{???}}Let $\pmb k$ be a field of characteristic $c$ and let $n$ and $N$ be positive integers.
Write $N$ in the form $N=\theta n+r$, for integers $\theta$ and $r$, with $0\le r\le n-1$. 
Then 
$$\operatorname{pd}_{R_{\pmb k,n}} Q_{\pmb k,n,N}=\infty \iff \theta\in S_c\ \text{and}\ 1\le r.$$
\endproclaim
\demo{Proof} Corollary \tref{COR} establishes ``$(\Leftarrow)$''.  Theorem \tref{H72+} shows that if $\theta\notin S_c$, then $\theta$ is in $T_c$ and Theorem \tref{HW672} shows that if $\theta\in T_c$, then $\operatorname{pd}_{R_{\pmb k,n}} Q_{\pmb k,n,N}$ is finite.
Furthermore, Observation \tref{O5.1} shows that if  $r=0$, then  $\operatorname{pd}_{R_{\pmb k,n}} Q_{\pmb k,n,N}$ is finite. \qed
\enddemo
We rephrase Theorem \tref{???} without appealing directly to the sets $S_c$ and $T_c$. Recall the meaning of the operation $\{\ \}$ from (\tref{rnd}).

\proclaim{Theorem  \tnum{!!!}}Let $\pmb k$ be a field of characteristic $c$ and let $n$ and $N$ be positive integers. Then $\operatorname{pd}_{R_{\pmb k,n}} Q_{\pmb k,n,N}$ is finite if and only if at least one of the following conditions hold{\rm:}
\roster\item $n$ divides $N$, or
\item $c=2$ and $n\le N$, or
\item $c=p$ is an odd prime and there exist an odd integer $J$ and a power $q=p^e$ of $p$ with $e\ge 1$ and $|Jq-\frac Nn|< \{\frac q3\}$.\endroster \flushpar In particular, if $c=0$, then $\operatorname{pd}_{R_{\pmb k,n}} Q_{\pmb k,n,N}$ is finite if and only if $n$ divides $N$.
\endproclaim
\demo{Proof} According to Theorem \tref{???},  $\operatorname{pd} Q_{\pmb k,n,N}$ is finite if and only if 
$\lfloor\frac{N}n\rfloor\in T_c$ or $n|N$. Definition \tref{D5.1} shows that $T_0$ is empty, $T_2$ is the set of all positive integers, and if $p\ge 3$, then
$$\tsize \lfloor\frac{N}n\rfloor \in T_p\iff Jq-\{\frac q3\}\le \lfloor\frac{N}n\rfloor<Jq+\{\frac q3\},$$for some $J$ and $q$ as described in the statement of the result.  The case $n|N$ has been taken care of in (1). When we restrict our attention to $n\not |N$, then $\lfloor\frac{N}n\rfloor<\frac Nn$ and 
$$\tsize \lfloor\frac{N}n\rfloor \in T_p\iff Jq-\{\frac q3\}\le \lfloor\frac{N}n\rfloor<\frac Nn<Jq+\{\frac q3\}\iff |Jq-\frac{N}n|<\{\frac q3\}. 
 \qed $$\enddemo

\SectionNumber=7\tNumber=1
\heading Section \number\SectionNumber. \quad The proof of Lemma \tref{NandI}.
\endheading

Recall, from  Definition \tref{D21.1}, that if  $d$, $a$,  and $b$  are non-negative integers then
$$\operatorname{Poly}_{d,a,b}(A,B)=\sum\limits_{i=0}^{d}(-1)^i \binom{a+d-i}{a}\binom{b+i}{b}A^{d-i}B^i$$ in   the polynomial ring $\Bbb Z[A,B]$. 
In this section we prove  Lemma \tref{NandI}.
\proclaim{Lemma \tref{NandI}}For each positive integer $\delta$, the polynomials
$$ P_{2\delta-1}(A,B,C)=\left\{\matrix \format\l\\
(-1)^{\delta}A\operatorname{Poly}_{\delta-1,\delta,\delta-1}(A,B)\operatorname{Poly}_{\delta-1,\delta,\delta-1}(A,C)\\
+B\operatorname{Poly}_{\delta-1,\delta,\delta-1}(B,A)\operatorname{Poly}_{\delta-1,\delta,\delta-1}(A,C)\\
+C\operatorname{Poly}_{\delta-1,\delta,\delta-1}(A,B)\operatorname{Poly}_{\delta-1,\delta,\delta-1}(C,A)\\
+ (-1)^{\delta+1}\binom{2\delta}{\delta}\binom{3\delta-1}{\delta-1}A^{2\delta-1}
\endmatrix\right. $$
and $$P_{2\delta}(A,B,C)=\left\{\matrix\format\l\\
(-1)^{\delta+1}\operatorname{Poly}_{\delta,\delta,\delta}(A,B)\operatorname{Poly}_{\delta,\delta,\delta}(A,C)\\
+B\operatorname{Poly}_{\delta-1,\delta,\delta}(B,A)\operatorname{Poly}_{\delta,\delta,\delta}(A,C)\\
+C\operatorname{Poly}_{\delta,\delta,\delta}(A,B)\operatorname{Poly}_{\delta-1,\delta,\delta}(C,A)\\
+(-1)^{\delta}\binom{2\delta}\delta\binom{3\delta}\delta A^{2\delta}\endmatrix \right.$$ are in the ideal $(A+B+C)\Bbb Z[A,B,C]$.\endproclaim

The proof of Lemma \tref{NandI} is carried out in Lemmas \tref{T?} and \tref{l2b}. In Lemma \tref{T?} we exhibit the  polynomial $Q_{2\delta-1}(A,B,C)\in \Bbb Z[A,B,C]$ such that 
$$P_{2\delta-1}(A,B,C)= (A+B+C)Q_{2\delta-1}(A,B,C)\tag\tnum{eq}$$ in $\Bbb Z[A,B,C]$.
In Lemma \tref{l2b} we
express $P_{2\delta}(A,B,C)$ as an element of the ideal  of $\Bbb Z[A,B,C]$ generated by 
$P_{2\delta-1}(A,B,C)$ and $A+B+C$; hence, $P_{2\delta}$ is also an element of the ideal $(A+B+C)\Bbb Z[A,B,C]$. 

\definition{Definition \tnum{D35}} \roster \item For integers $d$ and $\delta$, with $0\le d$, define the polynomial 
 $$H_{d,\delta}(A,B)=\sum\limits_{i=0}^{d}(-1)^i \binom{2\delta-1-i}{d-i}\binom{d+i}{d}A^{d-i}B^i
$$ of $\Bbb Z[A,B,C]$. \item For each positive integer $\delta$, define the polynomial $Q_{2\delta-1}(A,B,C)$  to be 
  $$\eightpoint  (-1)^{\delta+1}\sum\limits_{d=0}^{\delta-1}\frac{\binom {3\delta-1}{2\delta +d}\binom {2\delta -d-1}{\delta} (2d+1)}{\delta \binom{\delta-1}{d}^2}H_{d,\delta}(A,B)H_{d,\delta}(A,C)A^{2(\delta-d-1)}$$
in $\Bbb Q[A,B,C]$. \endroster \enddefinition

\proclaim{Lemma \tnum{T?}}If $\delta$ is a positive integer, then  the polynomial $Q_{2\delta-1}(A,B,C)$ is in $\Bbb Z[A,B,C]$ and  equation {\rm(\tref{eq})} holds in  $\Bbb Z[A,B,C]$.
\endproclaim

\demo{Proof}   The polynomial $P_{2\delta-1}(A,B,C)$ is in $\Bbb Z[A,B,C]$,  the polynomial $A+B+C$ generates a prime ideal of $\Bbb Z[A,B,C]$,  $Q_{2\delta-1}(A,B,C)$ is in $\Bbb Q[A,B,C]$, and the ring $\Bbb Z[A,B,C]$ is a Unique Factorization Domain. Once we prove that equation (\tref{eq}) holds in $\Bbb Q[A,B,C]$, then it  follows that $Q_{2\delta-1}(A,B,C)$ is actually in $\Bbb Z[A,B,C]$ and the equality (\tref{eq})  takes place in $\Bbb Z[A,B,C]$.

If $\delta=1$, then $P_{2\delta-1}(A,B,C)=A+B+C$ and $Q_{2\delta-1}(A,B,C)=1$; thus, equation (\tref{eq}) holds in this case. Henceforth, we assume that $2\le \delta$. 

To show that (\tref{eq}) holds in $\Bbb Q[A,B,C]$, we must show that
$$\eightpoint (*)=\left\{\matrix \format\l\\
\phantom{+}A\operatorname{Poly}_{\delta-1,\delta,\delta-1}(A,B)\operatorname{Poly}_{\delta-1,\delta,\delta-1}(A,C)\\
+(-1)^{\delta}B\operatorname{Poly}_{\delta-1,\delta,\delta-1}(B,A)\operatorname{Poly}_{\delta-1,\delta,\delta-1}(A,C)\\
+(-1)^{\delta}C\operatorname{Poly}_{\delta-1,\delta,\delta-1}(A,B)\operatorname{Poly}_{\delta-1,\delta,\delta-1}(C,A)\\
+ (A+B+C)\sum\limits_{d=0}^{\delta-1}\frac{\binom {2\delta -d-1}{\delta} \binom {3\delta-1}{2\delta +d}(2d+1)}{\delta \binom{\delta-1}{d}^2}H_{d,\delta}(A,B)H_{d,\delta}(A,C)A^{2(\delta-d-1)}\endmatrix\right. $$ 
equals $\binom{2\delta}{\delta}\binom{3\delta-1}{\delta-1}A^{2\delta-1}$. Observe that if $A^aB^bC^c$ appears in $(*)$, then $0\le b,c\le \delta$, $0\le a\le 2\delta-1$, and $a+b+c=2\delta-1$. 
If $S$ is a statement, then define
$$\chi(S)=\cases 1&\text{if $S$ is true}\\0&\text{if $S$ is false}.\endcases$$
Fix integers $b$, $c$, and $\delta$, with $0\le b,c\le \delta$ and $b+c\le 2\delta-1$. We see that the coefficient of $A^{2\delta-1-b-c}B^bC^c$ in $(*)$ is $(-1)^{b+c}$ times $(**)_{b,c,\delta}$, where $(**)_{b,c,\delta}$ is the rational number
$$ (**)_{b,c,\delta}=D_{b,c,\delta} +\sum\limits_{d=0}^{\delta-1}\frac{\binom {3\delta-1}{2\delta +d}\binom {2\delta -d-1}{\delta} (2d+1)}{\delta \binom{\delta-1}{d}^2}E_{b,c,d,\delta},\tag\tnum{**}$$ for
$$D_{b,c,\delta}= \left\{\matrix \format\l\\
\phantom{+} \chi(b,c\le \delta-1) \binom{2\delta-1-b}{\delta-1-b}\binom{\delta-1+b}{\delta-1} 
  \binom{2\delta-1-c}{\delta-1-c}\binom{\delta-1+c}{\delta-1} 
\\
+ 
\chi(1\le b)\chi(c\le \delta-1) \binom{2\delta-1-b}{\delta-1}\binom{\delta-1+b}{b-1}
 \binom{2\delta-1-c}{\delta-1-c}\binom{\delta-1+c}{\delta-1}
\\
+ \chi(b\le \delta-1)\chi(1\le c) \binom{2\delta-1-b}{\delta-1-b}\binom{\delta-1+b}{\delta-1}
 \binom{2\delta-1-c}{\delta-1}\binom{\delta-1+c}{c-1}
\endmatrix\right. $$ 
and $$E_{b,c,d,\delta}=\left\{\matrix \format\l\\
 \phantom{-} \chi(b,c\le d) \binom{2\delta-1-b}{d-b}\binom{d+b}{d}
  \binom{2\delta-1-c}{d-c}\binom{d+c}{d}
\\
 - \chi(1\le b\le d+1)\chi(c\le d) \binom{2\delta-b}{d-b+1}\binom{d+b-1}{d}
 \binom{2\delta-1-c}{d-c}\binom{d+c}{d}
\\
 -  \chi(b\le d)\chi(1\le c\le d+1)
 \binom{2\delta-1-b}{d-b}\binom{d+b}{d}
 \binom{2\delta-c}{d-c+1}\binom{d+c-1}{d}.
\endmatrix\right. $$
We complete the proof by showing that $$(**)_{b,c,\delta}=\cases \binom{2\delta}{\delta}\binom{3\delta-1}{\delta-1}&\text{if $b=c=0$}\\0&\text{otherwise}.\endcases\tag\tnum{goal}$$
The data has been chosen to ensure that $\delta$ is not zero; so, $\frac 1\delta$ is a well defined rational number. Use $$\tsize \binom{2\delta-1-\gamma}{\delta-1-\gamma}=\frac{\delta-\gamma}{\delta}\binom{2\delta-1-\gamma}{\delta-1}\quad\text{and}\quad 
\binom{\delta-1+\gamma}{\gamma-1}=\frac \gamma\delta\binom{\delta-1+\gamma}{\delta-1},$$ for $\gamma\in\{b,c\}$, to transform $D_{b,c,\delta}$ into 
$$\left\{\matrix \format\l \\
  \chi(b,c\le \delta-1) \frac{(\delta-b)}{\delta}\binom{2\delta-1-b}{\delta-1}\binom{\delta-1+b}{\delta-1} 
 \frac{(\delta-c)}\delta \binom{2\delta-1-c}{\delta-1}\binom{\delta-1+c}{\delta-1}  
\\
 + 
\chi(1\le b)\chi(c\le \delta-1) \binom{2\delta-1-b}{\delta-1}\frac{b}{\delta}\binom{\delta-1+b}{\delta-1} 
 \frac{(\delta-c)}\delta \binom{2\delta-1-c}{\delta-1}\binom{\delta-1+c}{\delta-1} 
\\
 + \chi(b\le \delta-1)\chi(1\le c) \frac{(\delta-b)}{\delta}\binom{2\delta-1-b}{\delta-1}\binom{\delta-1+b}{\delta-1}
\binom{2\delta-1-c}{\delta-1} \frac c\delta\binom{\delta-1+c}{\delta-1}. 
\endmatrix\right. $$ 
Observe that if $\gamma\in\{b,c\}$, then
$$(\delta-\gamma)\chi(\gamma\le \delta-1)=(\delta-\gamma)\quad\text{and}\quad 
 \gamma \chi(1\le \gamma)=\gamma.$$ Indeed, the ambient hypothesis $0\le b,c\le \delta$ ensures that $\chi(\gamma\le \delta-1)$ takes the value $1$ unless $\gamma=\delta$, in which case the numbers $(\delta-\gamma)$ and $\chi(\gamma\le \delta-1)$ are both zero; and $\chi(1\le \gamma)$ takes the value $1$ unless $\gamma=0$, in which case the numbers $\gamma$ and $\chi(1\le \gamma)$ are both zero. Thus; $$\tsize D_{b,c,\delta}=\left (
  \frac{(\delta-b)}{\delta} 
 \frac{(\delta-c)}\delta  
 + 
  \frac{b}{\delta} 
 \frac{(\delta-c)}\delta 
 +   \frac{(\delta-b)}{\delta}
 \frac c\delta
\right )\binom{2\delta-1-b}{\delta-1}\binom{\delta-1+b}{\delta-1}\binom{2\delta-1-c}{\delta-1}\binom{\delta-1+c}{\delta-1}; $$ and therefore,
$$D_{b,c,\delta}=\frac{(\delta^2-bc)}{\delta^2}\binom{2\delta-1-b}{\delta-1}\binom{\delta-1+b}{\delta-1}\binom{2\delta-1-c}{\delta-1}\binom{\delta-1+c}{\delta-1}.\tag\tnum{!!}$$ 
The parameter $\delta$ is at least $2$; so we may use the identities
$$\tsize \binom{2\delta-1-\gamma}{\delta-1}= \frac{\delta+1-\gamma}{\delta-1}\binom{2\delta-1-\gamma}{\delta-2}\quad \text{and}\quad 
\binom{\delta-1+\gamma}{\delta-1}=\frac \delta{\delta+\gamma}\binom{\delta+\gamma}{\delta},$$for $\gamma\in\{b,c\}$, to write 
$$\tsize D_{b,c,\delta}=\frac{(\delta^2-bc)(\delta+1-b)(\delta+1-c)}{(\delta-1)^2(\delta+b)(\delta+c)}\binom{2\delta-1-b}{\delta-2}\binom{\delta+b}{\delta}\binom{2\delta-1-c}{\delta-2}\binom{\delta+c}{\delta}.\tag\tnum{!!.5}$$
In a similar manner, we simplify $E_{b,c,d,\delta}$, provided $\max\{b,c,1\}\le d\le \delta-1$. 
Use $$\tsize \binom{2\delta-1-\gamma}{d-\gamma}=\frac{d-\gamma+1}{2\delta-1-d}  \binom{2\delta-1-\gamma}{2\delta-2-d},\quad
\binom{2\delta-\gamma}{d-\gamma+1}=\frac{2\delta-\gamma}{2\delta-1-d}  \binom{2\delta-1-\gamma}{2\delta-2-d},$$ 
and $\binom{d+\gamma-1}{d}=\frac{\gamma}{d+\gamma}\binom{d+\gamma}{d}$, for $\gamma\in\{b,c\}$, to obtain 
 $$ 
E_{b,c,d,\delta}=\left\{\matrix \format\l \\
\phantom{-}  \frac{d-b+1}{2\delta-1-d}  \binom{2\delta-1-b}{2\delta-2-d}\binom{d+b}{d}
 \frac{d-c+1}{2\delta-1-d} \binom{2\delta-1-c}{2\delta-2-d} \binom{d+c}{d} 
\\
 - \chi(1\le b)\frac{(2\delta-b)}{2\delta-1-d}\binom{2\delta-b-1}{2\delta-2-d} \frac{b}{d+b}\binom{d+b}{d}
 \frac{d-c+1}{2\delta-1-d} \binom{2\delta-1-c}{2\delta-2-d}\binom{d+c}{d} 
\\
 -  \chi(1\le c)
 \frac{d-b+1}{2\delta-1-d}  \binom{2\delta-1-b}{2\delta-2-d}\binom{d+b}{d}
 \frac{(2\delta-c)}{2\delta-1-d}\binom{2\delta-c-1}{2\delta-2-d}\frac{c}{d+c}\binom{d+c}{d}. \endmatrix\right.  $$
The hypothesis $\max\{b,c,1\}\le d\le \delta-1$ ensures that the denominators $2\delta-1-d$, $d+b$, and $d+c$ all are non-zero. We still have $\gamma\chi(1\le \gamma)=\gamma$ for $\gamma\in\{b,c\}$. It follows that 
$$\tsize E_{b,c,d,\delta}=\frac{E_{b,c,d,\delta}'}{(2\delta-1-d)^2(d+b)(d+c)} \binom{2\delta-1-b}{2\delta-2-d}\binom{d+b}{d}\binom{2\delta-1-c}{2\delta-2-d}\binom{d+c}{d} \tag\tnum{3!}$$
for$$\eightpoint E_{b,c,d,\delta}'=(d-b+1)(d+b)(d+c)(d-c+1) - (2\delta-b)b(d+c)(d-c+1) -  (d-b+1)(2\delta-c) c(d+b),$$ provided
$\max\{b,c,1\}\le d\le \delta-1$.

If $0\le r\le \delta$, then let
$$W_{b,c,r,\delta}=D_{b,c,d}+\sum\limits_{d=r}^{\delta-1} \frac{ \binom {3\delta-1}{2\delta +d}\binom {2\delta -d-1}{\delta}(2d+1)}{\delta \binom{\delta-1}{d}^2}E_{b,c,d,\delta}.\tag\tnum{def}$$Our goal is given in (\tref{goal}); we must evaluate the rational number $(**)_{b,c,\delta}$, which is given in (\tref{**}). Write 
$$(**)_{b,c,\delta}= F_{b,c,\delta}+G_{b,c,\delta}\tag\tnum{B/B}$$
for $$F_{b,c,\delta}= W_{b,c,\max\{b,c,1\},\delta}$$ and
$$G_{b,c,\delta}=\sum\limits_{d=0}^{\max\{b,c,1\}-1} \frac{\binom {3\delta-1}{2\delta +d}\binom {2\delta -d-1}{\delta} (2d+1)}{\delta \binom{\delta-1}{d}^2}E_{b,c,d,\delta}.$$

We next show, by induction,  that if $\max\{b,c,1\}\le r\le \delta$, then $$  \tsize \eightpoint W_{b,c,r,\delta}= \frac{\binom {3\delta-1}{2\delta+r-1}\binom {2\delta-1-r}{\delta-1}}{(2\delta-1-r)^2\binom{\delta}{r}^2} 
\frac{(r+1-b)(r+1-c) ( r^2 -b c)}{(r+b)(r+c)}
\binom{2\delta-1-b}{2\delta-2-r}\binom{r+b}{r}\binom{2\delta-1-c}{2\delta-2-r}\binom{r+c}{r}. \tag\tnum{Wr} $$  
If $r=\delta$, then (\tref{def}) shows that $W_{b,c,\delta,\delta}=D_{b,c,\delta}$; thus, $W_{b,c,\delta,\delta}$ is given in (\tref{!!.5}). There is no difficulty in seeing that (\tref{!!.5}) is the same as the right side of (\tref{Wr}) when $r=\delta$. Assume that (\tref{Wr}) holds at $r+1$ and $\max\{b,c,1\}\le r<\delta$. Then, according to (\tref{def}), 
$$W_{b,c,r,\delta}=W_{b,c,r+1,\delta}+\frac{\binom{3\delta-1}{2\delta+r}\binom{2\delta-r-1}{\delta}(2r+1)}{\delta\binom{\delta-1}r^2}E_{b,c,r,\delta}.\tag\tnum{t9}$$We see, from (\tref{Wr}), that $W_{b,c,r+1,\delta}$
$$\tsize =\frac{\binom {3\delta-1}{2\delta+r}\binom {2\delta-2-r}{\delta-1}}{(2\delta-2-r)^2\binom{\delta}{r+1}^2} 
\frac{(r+2-b)(r+2-c) ( (r+1)^2 -b c)}{(r+1+b)(r+1+c)}
\binom{2\delta-1-b}{2\delta-3-r}\binom{r+1+b}{r+1}\binom{2\delta-1-c}{2\delta-3-r}\binom{r+1+c}{r+1}.$$
Use $\binom{2\delta-1-\gamma}{2\delta-3-r}= \frac{2\delta-2-r}{r+2-\gamma}\binom{2\delta-1-\gamma}{2\delta-2-r}$ and $\binom{r+1+\gamma}{r+1}=\frac{r+1+\gamma}{r+1}\binom{r+\gamma}{r}$, for $\gamma\in\{b,c\}$, to see that
   $$\tsize W_{b,c,r+1,\delta}=\frac{\binom {3\delta-1}{2\delta+r}\binom {2\delta-2-r}{\delta-1}}{\binom{\delta}{r+1}^2} 
\frac{  ((r+1)^2 -b c)}{(r+1)^2}
\binom{2\delta-1-b}{2\delta-2-r} \binom{r+b}{r}\binom{2\delta-1-c}{2\delta-2-r} \binom{r+c}{r}.$$
The identities $(r+1)\binom{\delta}{r+1}=\delta\binom{\delta-1}r$ and $\frac 1\delta\binom{2\delta-2-r}{\delta-1}=\frac 1{2\delta-1-r}\binom{2\delta-1-r}\delta$ allow us to transform  $W_{b,c,r+1,\delta}$ into
$$\tsize W_{b,c,r+1,\delta}=\frac{\binom {3\delta-1}{2\delta+r}\binom {2\delta-1-r}{\delta}}{\delta\binom{\delta-1}{r}^2} 
\frac{  ((r+1)^2 -b c)}{2\delta-1-r}
\binom{2\delta-1-b}{2\delta-2-r} \binom{r+b}{r}\binom{2\delta-1-c}{2\delta-2-r} \binom{r+c}{r}.$$ 
Apply (\tref{t9}) and   (\tref{3!}) to see that $W_{b,c,r,\delta}$
$$= \left\{\matrix \format\l\\
\frac{\binom {3\delta-1}{2\delta+r}\binom {2\delta-1-r}{\delta}}{\delta\binom{\delta-1}{r}^2} 
\frac{  ((r+1)^2 -b c)}{2\delta-1-r}
\binom{2\delta-1-b}{2\delta-2-r} \binom{r+b}{r}\binom{2\delta-1-c}{2\delta-2-r} \binom{r+c}{r}\\
+\frac{\binom{3\delta-1}{2\delta+r}\binom{2\delta-r-1}{\delta}(2r+1)}{\delta\binom{\delta-1}r^2}
\frac{E_{b,c,r,\delta}'}{(2\delta-1-r)^2(r+b)(r+c)} \binom{2\delta-1-b}{2\delta-2-r}\binom{r+b}{r}\binom{2\delta-1-c}{2\delta-2-r}\binom{r+c}{r}, \endmatrix\right.$$
for$$\eightpoint E_{b,c,r,\delta}'=(r-b+1)(r+b)(r+c)(r-c+1) - (2\delta-b)b(r+c)(r-c+1) -  (r-b+1)(2\delta-c) c(r+b).$$ Thus, 
$W_{b,c,r,\delta}$
$$\eightpoint \tsize = \frac{\binom {3\delta-1}{2\delta+r}\binom {2\delta-1-r}{\delta}\binom{2\delta-1-b}{2\delta-2-r}\binom{r+b}{r}\binom{2\delta-1-c}{2\delta-2-r}\binom{r+c}{r}}{(2\delta-1-r)^2\delta\binom{\delta-1}{r}^2(r+b)(r+c)} \left\{\matrix \format\l\\ 
  ((r+1)^2 -b c) (2\delta-1-r)(r+b)(r+c)
\\ 
+(2r+1)
E_{b,c,r,\delta}'  \endmatrix\right.$$
One may check that 
$$\eightpoint ((r+1)^2 -b c) (2\delta-1-r)(r+b)(r+c)
+(2r+1)
E_{b,c,r,\delta}'=(r -b + 1)(r - c + 1)  (r^2-bc) (2 \delta + r);$$
hence, $W_{b,c,r,\delta}$ equals
$$\eightpoint \tsize  \frac{\binom {3\delta-1}{2\delta+r}\binom {2\delta-1-r}{\delta}\binom{2\delta-1-b}{2\delta-2-r}\binom{r+b}{r}\binom{2\delta-1-c}{2\delta-2-r}\binom{r+c}{r}}{(2\delta-1-r)^2\delta\binom{\delta-1}{r}^2(r+b)(r+c)}  
    (r -b + 1)(r - c + 1)(r^2-bc) (2 \delta + r).  $$
Apply the identities
$$\eightpoint \tsize \delta\binom{\delta-1}r=(\delta-r)\binom \delta r,\quad \delta\binom{2\delta-1-r}\delta=(\delta-r)\binom{2\delta-1-r}{\delta-1},\quad\text{and}\quad
(2\delta+r)\binom{3\delta-1}{2\delta+r}=(\delta-r)\binom{3\delta-1}{2\delta+r-1},$$ in order to see that 
$$\eightpoint \tsize  \frac{\binom {3\delta-1}{2\delta+r}\binom {2\delta-1-r}{\delta}(2 \delta + r)}{\delta\binom{\delta-1}{r}^2}
=\frac{\binom {3\delta-1}{2\delta+r-1}\binom {2\delta-1-r}{\delta-1}}{\binom{\delta}{r}^2}$$ 
and 
$$  \tsize \eightpoint W_{b,c,r,\delta}= \frac{\binom {3\delta-1}{2\delta+r-1}\binom {2\delta-1-r}{\delta-1}}{(2\delta-1-r)^2\binom{\delta}{r}^2} 
\frac{(r+1-b) (r+1-c)  ( r^2 -b c) }{(r+b)(r+c)}
 \binom{2\delta-1-b}{2\delta-2-r}\binom{r+b}{r}\binom{2\delta-1-c}{2\delta-2-r}\binom{r+c}{r};   $$ thereby completing the proof of (\tref{Wr}).

We saw at (\tref{B/B}) that $(**)_{b,c,\delta}= F_{b,c,\delta}+G_{b,c,\delta}$. All of the formulas are symmetric in $b$ and $c$; so no harm is done if we only check (\tref{goal}) for $b\le c$. Apply (\tref{Wr}) at $$r=\max\{b,c,1\}=\cases 1&\text{if $b=c=0$}\\c&\text{$b,1\le c$}\endcases$$ 
to see that 
$$\tsize \eightpoint F_{b,c,\delta}= \cases \frac{4\binom {3\delta-1}{2\delta}\binom {2\delta-2}{\delta-1}}{(2\delta-2)^2\delta^2} 
\binom{2\delta-1}{2\delta-3}^2,&\text{if $b=c=0$,} \\ 
\frac{\binom {3\delta-1}{2\delta+c-1}\binom {2\delta-1-c}{\delta-1}}{(2\delta-1-c)^2\binom{\delta}{c}^2} 
\frac{(c+1-b)( c -b)}{2(c+b)}
\binom{2\delta-1-b}{2\delta-2-c}\binom{c+b}{c}\binom{2\delta-1-c}{2\delta-2-c}\binom{2c}{c},& \text{if $b,1\le c$.}\endcases$$ 
On the other hand, 
$$G_{b,c,\delta}=\sum\limits_{d=0}^{\max\{b,c,1\}-1} \frac{\binom {3\delta-1}{2\delta +d}\binom {2\delta -d-1}{\delta} (2d+1)}{\delta \binom{\delta-1}{d}^2}E_{b,c,d,\delta},$$
with $$E_{b,c,d,\delta}=\left\{\matrix \format\l\\
 \phantom{-} \chi(b,c\le d) \binom{2\delta-1-b}{d-b}\binom{d+b}{d}
  \binom{2\delta-1-c}{d-c}\binom{d+c}{d}
\\
 - \chi(1\le b\le d+1)\chi(c\le d) \binom{2\delta-b}{d-b+1}\binom{d+b-1}{d}
 \binom{2\delta-1-c}{d-c}\binom{d+c}{d}
\\
 -  \chi(b\le d)\chi(1\le c\le d+1)
 \binom{2\delta-1-b}{d-b}\binom{d+b}{d}
 \binom{2\delta-c}{d-c+1}\binom{d+c-1}{d}.
\endmatrix\right. $$ If $b=c=0$ and $0\le d\le \max\{b,c,1\}-1$, then $\chi(1\le b)=\chi(1\le c)=0$ and $d=0$. If $b,1\le c$ and $d\le \max\{b,c,1\}-1$, then $\chi(c\le d)=0$ and $\chi(c\le d+1)$ is non-zero only when $d=c-1$; so, 
$$G_{b,c,\delta}=\cases 
 \frac{\binom {3\delta-1}{2\delta }\binom {2\delta-1}{\delta}}{\delta }, 
 &
\text{if $b=c=0$},\\
 - \frac{\binom {3\delta-1}{2\delta +c-1}\binom {2\delta -c}{\delta} (2c-1)  \chi(b\le c-1)}{\delta \binom{\delta-1}{c-1}^2} 
 \binom{2\delta-1-b}{c-1-b}\binom{c-1+b}{c-1}
 \binom{2c-2}{c-1},
 &\text{if $b,1\le c$}.\endcases
$$ 
Use (\tref{B/B}) to see that if $b=c=0$, then 
$$(**)_{0,0,\delta}= \frac{4\binom {3\delta-1}{2\delta}\binom {2\delta-2}{\delta-1}}{(2\delta-2)^2\delta^2} 
\binom{2\delta-1}{2\delta-3}^2+ \frac{\binom {3\delta-1}{2\delta }\binom {2\delta-1}{\delta}}{\delta }.$$
Apply $$\tsize \frac {4\binom{2\delta-1}{2\delta-3}^2}{(2\delta-2)^2}=(2\delta-1)^2\quad\text{and}\quad \binom{2\delta-2}{\delta-1}(2\delta-1)=\binom{2\delta-1}\delta\delta$$ to see that 
$$(**)_{0,0,\delta}= \frac{\binom {3\delta-1}{2\delta }\binom {2\delta-1}{\delta}}{\delta }((2\delta-1)+1)=\binom {3\delta-1}{2\delta }\binom {2\delta}{\delta},$$as required in (\tref{goal}). 

If $b,1\le c$, then (\tref{B/B}) gives that
$$(**)_{b,c,\delta}=\left\{\matrix\format\l\\
\phantom{-}\frac{\binom {3\delta-1}{2\delta+c-1}\binom {2\delta-1-c}{\delta-1}}{(2\delta-1-c)^2\binom{\delta}{c}^2} 
\frac{(c+1-b)( c -b)}{2(c+b)}
\binom{2\delta-1-b}{2\delta-2-c}\binom{c+b}{c}\binom{2\delta-1-c}{2\delta-2-c}\binom{2c}{c}\\
 - \frac{\binom {3\delta-1}{2\delta +c-1}\binom {2\delta -c}{\delta} (2c-1)  \chi(b\le c-1)}{\delta \binom{\delta-1}{c-1}^2} 
 \binom{2\delta-1-b}{c-1-b}\binom{c-1+b}{c-1}
 \binom{2c-2}{c-1}\endmatrix \right.$$
In the present calculation, $b$ and $c$ are integers with $b\le c$; thus, $$\chi(b\le c-1)=\cases 0&\text{if $b=c$}\\1&\text{if $b<c$}\endcases$$ and $(c-b)=(c-b)\chi(b\le c-1)$. It follows that
$$\eightpoint \tsize (**)_{b,c,\delta}=\chi(b\le c-1)\binom {3\delta-1}{2\delta+c-1}\left\{\matrix\format\l\\
\phantom{-}\frac{\binom {2\delta-1-c}{\delta-1}}{(2\delta-1-c)^2\binom{\delta}{c}^2} 
\frac{(c+1-b)( c -b)}{2(c+b)}
\binom{2\delta-1-b}{2\delta-2-c}\binom{c+b}{c}\binom{2\delta-1-c}{2\delta-2-c}\binom{2c}{c}\\
 - \frac{\binom {2\delta -c}{\delta} (2c-1)  }{\delta \binom{\delta-1}{c-1}^2} 
 \binom{2\delta-1-b}{c-1-b}\binom{c-1+b}{c-1}
 \binom{2c-2}{c-1}.\endmatrix \right.$$Use 
$$\tsize (2c-1)\binom{2c-2}{c-1}=\binom{2c}c\frac c2,\quad \frac{\binom{2\delta-1-c}{2\delta-2-c}}{2\delta-1-c}=1,\quad
\frac{(c+1-b)(c-b)}{(2\delta-c-1)}\binom{2\delta-1-b}{2\delta-2-c}=\binom{2\delta-1-b}{c-b-1}(2\delta-c),$$
$$\tsize \frac{\binom{c+b}c}{c+b}= \binom{c+b-1}{c-1}\frac 1c,\quad (2\delta-c)\binom{2\delta-1-c}{\delta-1}=\binom{2\delta-c}\delta\delta,\quad \text{and}\quad \delta\binom{\delta-1}{c-1}=\binom \delta cc.$$ to see that $(**)_{b,c,\delta}=0$ for $1,b\le c$, as required by (\tref{goal}), and the proof is complete. \qed \enddemo

\proclaim{Lemma \tnum{l2b}}If $\delta$ is a positive integer, then the polynomial 
$$\eightpoint P_{2\delta}(A,B,C)+3AP_{2\delta-1}(A,B,C)+(-1)^{\delta}A(A+B+C)\operatorname{Poly}_{\delta-1,\delta,\delta}(A,B)\operatorname{Poly}_{\delta-1,\delta,\delta}(A,C)$$of $\Bbb Z[A,B,C]$ is the zero polynomial.
\endproclaim
\demo{Proof} Let $X$ be the polynomial which is recorded in the statement. Use the fact that $3\binom{3\delta-1}{\delta-1}=\binom{3\delta}{\delta}$ to see that 
$$X=\left\{\matrix\format\l\\
(-1)^{\delta+1}\operatorname{Poly}_{\delta,\delta,\delta}(A,B)\operatorname{Poly}_{\delta,\delta,\delta}(A,C)\\
+B\operatorname{Poly}_{\delta-1,\delta,\delta}(B,A)\operatorname{Poly}_{\delta,\delta,\delta}(A,C)\\
+C\operatorname{Poly}_{\delta,\delta,\delta}(A,B)\operatorname{Poly}_{\delta-1,\delta,\delta}(C,A)\\
\\
(-1)^{\delta}3A^2\operatorname{Poly}_{\delta-1,\delta,\delta-1}(A,B)\operatorname{Poly}_{\delta-1,\delta,\delta-1}(A,C)\\
+3AB\operatorname{Poly}_{\delta-1,\delta,\delta-1}(B,A)\operatorname{Poly}_{\delta-1,\delta,\delta-1}(A,C)\\
+3AC\operatorname{Poly}_{\delta-1,\delta,\delta-1}(A,B)\operatorname{Poly}_{\delta-1,\delta,\delta-1}(C,A)\\
\\
+(-1)^{\delta}A(A+B+C)\operatorname{Poly}_{\delta-1,\delta,\delta}(A,B)\operatorname{Poly}_{\delta-1,\delta,\delta}(A,C).\endmatrix\right.$$
The polynomial $X$ in $\Bbb Z[A,B,C]$ is homogeneous of degree $2\delta$. For non-negative integers $b$ and $c$, the coefficient of $A^{2\delta-b-c}B^bC^c$ in $X$ is
$$X_{b,c}=(-1)^{\delta+b+c}\left\{\matrix\format\l\\
-\binom{2\delta-b}{\delta}\binom{\delta+b}{\delta}\binom{2\delta-c}{\delta}\binom{\delta+c}{\delta}\\
+\binom{\delta+b-1}{\delta}\binom{2\delta-b}{\delta}\binom{2\delta-c}{\delta}\binom{\delta+c}{\delta}\\
+\binom{2\delta-b}{\delta}\binom{\delta+b}{\delta}\binom{\delta+c-1}{\delta}\binom{2\delta-c}{\delta}
\\
\\
+3\binom{2\delta-b-1}{\delta}\binom{\delta+b-1}{\delta-1}\binom{2\delta-c-1}{\delta}\binom{\delta+c-1}{\delta-1}\\
+3\binom{\delta+b-1}{\delta}\binom{2\delta-b-1}{\delta-1}\binom{2\delta-c-1}{\delta}\binom{\delta+c-1}{\delta-1}\\
+3\binom{2\delta-b-1}{\delta}\binom{\delta+b-1}{\delta-1}\binom{\delta+c-1}{\delta}\binom{2\delta-c-1}{\delta-1}\\
\\
+\binom{2\delta-b-1}{\delta}\binom{\delta+b}{\delta}\binom{2\delta-c-1}{\delta}\binom{\delta+c}{\delta}\\
-\binom{2\delta-b}{\delta}\binom{\delta+b-1}{\delta}\binom{2\delta-c-1}{\delta}\binom{\delta+c}{\delta}\\
-\binom{2\delta-b-1}{\delta}\binom{\delta+b}{\delta}\binom{2\delta-c}{\delta}\binom{\delta+c-1}{\delta}.\endmatrix\right.$$We use 
$$\matrix \format \c&\ \c\ &\c\\(\delta+b) \binom{\delta+b-1}{\delta-1}&=& \binom {\delta+b}\delta \delta\\
(\delta+b) \binom{\delta+b-1}{\delta}&=& \binom {\delta+b}\delta b\\
(\delta+c) \binom{\delta+c-1}{\delta-1}&=& \binom {\delta+c}\delta \delta\\
(\delta+c) \binom{\delta+c-1}{\delta}&=& \binom {\delta+c}\delta c\\
\endmatrix\qquad\qquad \matrix  \format \c&\ \c\ &\l\\
(2\delta-b)\binom{2\delta-1-b}{\delta}&=&\binom{2\delta-b}{\delta}(\delta-b)\\
(2\delta-b)\binom{2\delta-1-b}{\delta-1}&=&\binom{2\delta-b}{\delta}\delta \\
(2\delta-c)\binom{2\delta-1-c}{\delta}&=&\binom{2\delta-c}{\delta}(\delta-c)\\
(2\delta-c)\binom{2\delta-1-c}{\delta-1}&=&\binom{2\delta-c}{\delta}\delta \endmatrix$$ to write
$$(\delta+b)(\delta+c)(2\delta-b)(2\delta-c)X_{b,c}=(-1)^{\delta+b+c}X'_{b,c}\binom{2\delta-b}{\delta}\binom{\delta+b}\delta\binom{\delta+c}\delta\binom{2\delta-c}{\delta},$$ where $X'_{b,c}=X^{(1)}_{b,c}+X^{(2)}_{b,c}+X^{(3)}_{b,c}$, for
$$\eightpoint \split X^{(1)}_{b,c}&{}= - (2\delta-b)(\delta+b)(2\delta-c)(\delta+c) 
+ (2\delta-b)b(2\delta-c)(\delta+c) 
+  (2\delta-b)(\delta+b)(2\delta-c)c\\&{}=-(2\delta-b)(2\delta-c)(\delta^2-bc),\\
X^{(2)}_{b,c}&{}=3\delta(\delta-b)\delta(\delta-c)
+3 \delta b\delta(\delta-c)
+3 \delta(\delta-b)\delta c=3\delta^2(\delta^2-bc),\\
X^{(3)}_{b,c}&{}=(\delta-b)(\delta+b)(\delta-c)(\delta+c)-(2\delta-b) b(\delta-c)(\delta+c)
-(\delta-b)(\delta+b)(2\delta-c) c\\&{}=(\delta^2-2\delta b-2\delta c+bc)(\delta^2-bc).\endsplit $$Thus, $X'_{b,c}$ is equal to 
$$[-(2\delta-b)(2\delta-c)+3\delta^2+(\delta^2-2\delta b-2\delta c+bc)](\delta^2-bc)=0.$$ The coefficient $X_{b,c}$ is also zero and $X=0$. \qed \enddemo

\SectionNumber=8\tNumber=1
\heading Section \number\SectionNumber. \quad Applications to Frobenius powers.
\endheading

In this section the characteristic of the field $\pmb k$ is $p>0$.  Keep in mind that the $t^{\text{th}}$-Frobenius power of $Q_{\pmb k,n,N}$ is
$$F^t(Q_{\pmb k,n,N})=Q_{\pmb k,n,p^tN}.$$
Three questions are investigated in this section. When is there a module $Q_{\pmb k,n,N}$ of infinite projective dimension such that some Frobenius power of it has finite projective dimension? Is the tail of the resolution of $Q_{\pmb k,n,p^tN}$ (up to shift)  eventually a periodic function of $t$? Can one use socle degrees to predict that the tail of the resolution of $Q_{\pmb k,n,p^tN}$ is a shift of the tail of the resolution of  $Q_{\pmb k,n,N}$? The third question was answered in \cite{\rref{KV}} in the case that both modules have finite projective dimension (hence the infinite tail of both resolutions  is zero). It is shown in \cite{\rref{KU}} how the socle degrees can be used to predict that the tail of the resolution of $Q_{\pmb k,n,p^tN}$ is a shift of the tail of the resolution of  $Q_{\pmb k,n,N}$, {\bf as a graded module}. We show that the condition of \cite{\rref{KU}} actually produces an isomorphism of complexes {\bf with differential}. 

\heading Subsection \number\SectionNumber.1. \quad When is there a module $Q_{\pmb k,n,N}$ of infinite projective dimension such that some Frobenius power of it has finite projective dimension?
\endheading

At the Georgia State University -- University of South Carolina Commutative Algebra Seminar in October, 2006, Florian Enescu asked, ``When does the ring $R$ have an ideal $J$ so that $R/J$ has infinite projective dimension but $R/J^{[p]}$ has finite projective dimension?'' The participants at the seminar produced a partial answer to this result.    Recall that the ideal $I$ is {\it not Frobenius closed} if there is an element $r\in R$ with $r\notin I$, but $r^q\in I^{[q]}$ for some $q$. 
The ring $(R,\frak m)$ is {\it $F$-injective} if the Frobenius map induces an injective map on  local cohomology
  modules $\operatorname{H}_{\frak m}^i(R) \to \operatorname{H}_{\frak m}^i(R)$, for every $0 \le  i \le \dim(R)$. A   local Cohen-Macaulay ring is $F$-injective if and only if all of its   parameter ideals are Frobenius closed. 

\proclaim{Proposition \tnum{EKVY}} {\bf ([Enescu, Kustin, Vraciu, Yao])} Let $(R,\frak m)$ be a   Cohen-Macaulay ring. If $R$  is not $F$-injective, then there exists an $\frak m$-primary ideal $J$ in $R$ with $\operatorname{pd}_RR/J=\infty$, but $\operatorname{pd}_RR/J^{[p]}$ finite.  \endproclaim
 
\demo{Proof} The ring $R$ is not $F$-injective; so, there exists a system of parameters $\pmb x=(x_1,\dots x_d)$ for $R$ (with $d=\dim R$) and an element $u$ of $R$ with $u\notin (\pmb x)$ but $u^p\in (\pmb x^{[p]})$. Replace $u$ by some multiple of $u$, if necessary, in order to have $u$ in the ideal $(\pmb x)\:\! \frak m$. The ideal $I=(\pmb x,u)$ is linked to the maximal ideal of a local ring which is not regular; so, $I$ has infinite projective dimension. On the other hand, $I^{[p]}=(x_1^p,\dots x_d^p)$, because $u^p$ is already in the ideal on the right side; thus, $I^{[p]}$ is generated by a regular sequence and has finite projective dimension. \qed \enddemo

It is natural to ask the following question. 
 
\proclaim{Question \tnum{Q3.5}}Does  the converse  of Proposition \tref{EKVY} hold? \endproclaim

Recently, Hailong Dao showed us that the answer to Question \tref{Q3.5} would be ``No'' if the requirement that $J$ is $\frak m$-primary were removed. 
\example{Example \tnum {Dao}} {\bf([Dao])} Let $\pmb k$ be a field of characteristic $2$. If $J$ is the ideal $(x,z)$ of the ring $R=\pmb k[x,y,z]/(xy-z^2)$, then $I$ has infinite projective dimension; but $I^{[2]}=(x^2,z^2)=(x^2,xy)\cong (x,y)$ has finite projective dimension. The ring $R$ is  $F$-injective because $xy-z^2$ is not in $\frak M^{[2]}$ where $\frak M$ is the maximal ideal $(x,y,z)$ of the polynomial ring $\pmb k[x,y,z]$; see \cite{\rref{F}, Prop.~2.1 and Lemma~3.3}. \endexample

It is known, by work of Fedder, which of the rings $R_{\pmb k,n}$ are not $F$-injective. Theorem \tref{F-p} shows that the ideals $(x^N,y^N,z^N)$ in the rings $R_{\pmb k,n}$ can not provide counter examples to Question \tref{Q3.5}. 

\proclaim{Theorem \tnum{F-p}}Let $\pmb k$ be a field of positive characteristic $p$ and let $n\ge 2$ be  an integer not divisible by $p$.
The following statements are equivalent.
\roster
\item There exists a positive integer $N$ such that  $\operatorname{pd}_{R_{\pmb k, n}}Q_{\pmb k, n,N}$ is infinite  and $\operatorname{pd}_{R_{\pmb k, n}}Q_{\pmb k, n,pN}$ is finite.

\item The integer $n$ satisfies $n\ge 4$, or else $n=3$ and $p\equiv 2\mod 3$.

\item  The ring $R_{\pmb k, n}$ is not $F$-injective.
\endroster
\endproclaim
\remark{Remark} The implication $(3)\Rightarrow (1)$ says more than Proposition \tref{EKVY} because Proposition \tref{EKVY} guarantees that when $R$ is not $F$-injective, then there exists   an $\frak m$-primary ideal $J$ with infinite projective dimension such that $J^{[p]}$ has finite projective dimension and the implication $(3)\Rightarrow (1)$ says that if, in addition,  $R$ has the form $R_{\pmb k,n}$, then $J$ may be taken to have the form $(x^N,y^N,z^N)$ for some $N$.\endremark
\demo{Proof} The equivalence of (2) and (3) is known from \cite{\rref{F}, Prop.~2.1 and Lemma~3.3}. The implication $(1) \Rightarrow (2)$ is established in Propositions \tref{28.2} and \tref{28.4} where it is shown  that if $n=2$ and $p$ is odd or else if $n=3$ and $p\equiv 1\mod 3$, then $$\operatorname{pd}_{R_{\pmb k, n}}Q_{\pmb k, n,pN}<\infty\implies \operatorname{pd}_{R_{\pmb k, n}}Q_{\pmb k, n,N}<\infty.$$ We do not consider  $p=n=3$ or $p=n=2$.  

We now prove $(2) \Rightarrow (1)$. If $p=2$, $n\ge 3$, and $N=n-1$, then $N<n<2N$; and therefore Theorem \tref{!!!} shows that $Q_{\pmb k, n,N}$ has infinite projective dimension, but $Q_{\pmb k, n,2N}$ has  finite projective dimension. 

Suppose $4\le n$ and $p=3$. Let $N=n+1$. Thus, $N=\theta n+r$ with $\theta=r=1$. We see from Example \tref{!E18.2} that $1\in S_3$; so Theorem \tref{???}  shows that $Q_{\pmb k, n,N}$ has infinite projective dimension. On the other hand, $3N=\theta' n+r'$ with $\theta'=r'=3$. (We still have $1\le r'\le n-1$.) We see in Example \tref{E1.2} that $3\in T_3$; so $3\notin S_3$ and  $Q_{\pmb k, n,3N}$ has finite projective dimension.

Suppose $4\le n$ and $p\equiv 1\mod 3$. Let $N=n(p-\frac{p-1}3)-1$. We see that 
$N$ is also equal to $\theta n+r$ for $\theta=p-\frac{p-1}3-1$ and $r=n-1$. In the language of Definition \tref{!D18.1}, we have
$$\tsize \theta= \frac{2(p-1)}3=2\pi_p\in S_p.$$ It follows from Theorem \tref{???} that  $Q_{\pmb k, n,N}$ has infinite projective dimension. On the other hand, $\frac{pN}n=p^2-p\left(\frac{p-1}3\right)-\frac pn$, and
$$\tsize\left\vert p^2-\frac{pN}n\right\vert =\frac {p^2-p}3+\frac pn<\frac{p^2-1}3.$$ The inequality on the right holds because
$$\tsize-\frac p3+\frac pn<-\frac 13\iff \frac 13+\frac pn<\frac p3\iff 1<\left(\frac {n-3}n\right)p.$$The parameters $n$ and $p$ satisfy $4\le n$ (so $\frac 14\le \frac{n-3}n$) and $7\le p$ (so $1<\left(\frac {n-3}n\right)p$). Apply Theorem \tref{!!!} with $J=1$ and $q=p^2$ to conclude that $Q_{\pmb k, n,pN}$ has  finite projective dimension. 

Finally, suppose that $p$ is an odd prime with $p\equiv 2\mod 3$ and $n\ge 3$. Let $N=n(p-\frac {p+1}3)-1$.  We see that 
$N$ is also equal to $\theta n+r$ for $\theta=p-\frac{p+1}3-1$ and $r=n-1$. In the language of Definition \tref{!D18.1}, we have
$$\tsize \theta= \frac{2(p-2)}3=2\pi_p\in S_p.$$ It follows from Theorem \tref{???} that  $Q_{\pmb k, n,N}$ has infinite projective dimension. On the other hand, $\frac {pN}n=p(p-\frac{p+1}3)-\frac pn$. Take $J$ to be the odd integer $p-\frac{p+1}3$. We see that 
$$\tsize |Jp-\frac {pN}n|= \frac pn<\frac{p+1}3.$$ The last inequality holds because $3\le n$. It follows from Theorem \tref{!!!} that $Q_{\pmb k, n,pN}$ has  finite projective dimension. \qed \enddemo

\proclaim{Proposition \tnum{28.2}} Consider the data $(\pmb k,n,N)$ with $\pmb k$ a field of characteristic $p$. If $p\equiv 1\mod 3$ and $n=3$, then 
$$  \operatorname{pd} Q_{\pmb k,n,pN}<\infty\implies   \operatorname{pd} Q_{\pmb k,n,N}<\infty.$$ 
\endproclaim

\demo{Proof}  If $n$ divides $pN$, then $n$ divides $N$ (since $n=3$ and $p$ is not divisible by $3$) and $\operatorname{pd} Q_{\pmb k,n,N}$ is automatically finite.
Throughout the rest of the proof we assume that $n$ does not divide $pN$. 
According to Theorem \tref{!!!}, there exists an odd integer $J$ and a power $q=p^e$ of $p$ with $e\ge 1$ and 
$|Jq-\frac {pN}3|< \frac{q-1}3$.  Multiply by $3$ to see 
$|3Jq- pN| <q-1$. The integer $pN$ is divisible by $p$. There are no integers in the interval $(q-p,q-1)$ which are divisible by $p$; hence,
$|3Jq-pN| \le  q-p$. The integers $q-p$ and $3Jq$ are divisible by $3$. The integer $pN$ is not divisible by $3$; so
$|3Jq-pN| <  q-p$. Let $q'= \frac qp$ and divide by $3p$  to see that
$|Jq'-\frac N3|< \frac{q'-1}3$. This inequality shows that $1<q'$ and Theorem \tref{!!!} gives $ \operatorname{pd} Q_{\pmb k,n,N}<\infty$. \qed \enddemo

\proclaim{Proposition \tnum{28.4}} Consider the data $(\pmb k,n,N)$ with $\pmb k$ a field of characteristic $p$. If $p$ is an odd prime and $n=2$, then 
$$  \operatorname{pd} Q_{\pmb k,n,pN}<\infty\implies   \operatorname{pd} Q_{\pmb k,n,N}<\infty.$$ 
\endproclaim

\demo{Proof}  If $n$ divides $pN$, then $n$ divides $N$ (since $n=2$ and $p$ is not divisible by $2$) and $\operatorname{pd} Q_{\pmb k,n,N}$ is automatically finite.
Throughout the rest of the proof we assume that $n$ does not divide $pN$. 
According to Theorem \tref{!!!}, there exists an odd integer $J$ and a power $q=p^e$ of $p$ with $e\ge 1$ and 
$$\tsize  |Jq- \frac {pN}2|< \{\frac{q}3\}.\tag\tnum{zz}$$ 
We first treat the case $p=3$. In this case, $\{\frac{q}3\}=\frac{q}3$. Divide (\tref{zz}) by $3$ to obtain $|J \frac q3-\frac N2| < \frac{q}9$.  The hypotheses that $N$ is an  odd integer guarantees that $q>3$ and therefore Theorem \tref{!!!} shows that $\operatorname{pd} Q_{\pmb k,n,N}$ is finite.

Henceforth, we assume that $p\ge 5$. Multiply (\tref{zz}) by $6$ to see 
$$|6Jq- 3pN| < {\tsize 6\{\frac{q}3\}}.\tag\tnum{aa}$$ 
The integer $3pN$ is divisible by $3p$. 
The intervals 
$$\cases (2(q-p),6\{\frac{q}3\})&\text{if $\frac qp\equiv 1\mod 3$}\\
(6\{\frac{q}3\},2(q+p))&\text{if $\frac qp\equiv 2\mod 3$}\endcases$$ do not contain any integers which are divisible by $3p$. Hence, (\tref{aa}) is equivalent to
$$\cases
|6Jq-3pN| \le  2(q-p)&\text{if $\frac qp\equiv 1\mod 3$}\\
|6Jq- 3pN| <  2(q+p)&\text{if $\frac qp\equiv 2\mod 3$}.\endcases$$
The integers $6Jq$ and $2(q-p)$ are divisible by $2$. The integer $3pN$ is not divisible by $2$; so
$$\cases
|6Jq- 3pN| <  2(q-p)&\text{if $\frac qp\equiv 1\mod 3$}\\
|6Jq- 3pN| <  2(q+p)&\text{if $\frac qp\equiv 2\mod 3$}.\endcases$$ Let $q'={\tsize \frac qp}$ and  divide by $6p$ to see that
$$\cases |Jq'-\frac N2|< \frac{q'-1}3&\text{if $\frac qp\equiv 1\mod 3$}\\
|Jq'-\frac N2|< \frac{q'+1}3
&\text{if $\frac qp\equiv 2\mod 3$}.\endcases$$Thus, $1<q'$, 
$ |Jq'-\frac N2|< \{\frac {q'}{3}\}$, and   Theorem \tref{!!!} gives $ \operatorname{pd} Q_{\pmb k,n,N}<\infty$. \qed \enddemo

\heading Subsection \number\SectionNumber.2. \quad Is the tail of the resolution of $Q_{\pmb k,n,p^tN}$ (up to shift)  eventually a periodic function of $t$?
\endheading

\definition {Data \tnum{D8}}Fix data $(\pmb k,n)$, where $\pmb k$ is a field of   characteristic $p>0$ and  $n\ge 2$ is an integer which is relatively prime to $p$. Fix the ring $R=R_{\pmb k,n}$ and let $Q_N$ denote the $R$-module $Q_{\pmb k,n,N}$. Let $\pmb o$ be the  the order of $p$ in the group of multiplicative units in the ring $\Bbb Z/(n)$.\enddefinition

 In Theorem \tref{PY}, we show  that, for each positive integer $N$, the tail of the resolution of $Q_{p^tN}$, up to shift,  is eventually a periodic function of $t$ of period at most $2\pmb o$.  The ``tail of the resolution'' of $Q_N$ is equal to the resolution of the second syzygy of the $R$-module $Q_{N}$, which we denote by $\operatorname{syz}_2Q_N$ (see (\tref{trun})); therefore, the tail of the resolution of $Q_{qN}$ is isomorphic to a shift of the tail of the resolution of $Q_N$ if and only if  $\operatorname{syz}_2Q_{qN}$ is isomorphic to a shift of $\operatorname{syz}_2Q_N$.

\proclaim{Theorem \tnum{PY}} Adopt Data \tref{D8}. Then there exists a power $q=p^e$, of $p$,  for some integer $e\ge 1$, such that, for all positive integers $N$, there exists an integer $t_0(N)$, which depends on $N$, so that
$$\text{$\operatorname{syz}_2 Q_{q^sp^tN}$ is isomorphic to a shift of $\operatorname{syz}_2 Q_{p^tN}$}\tag\tnum{syz}$$ for all integers $s$ and $t$ with $s\ge 0$ and $t\ge t_0(N)$. Furthermore, $e$ may be taken to be at most $2 \pmb o$.\endproclaim

\demo{Proof} There are two cases. Either 
\medskip
\flushpar {\bf Case 1:} there exists an integer $t_0=t_0(N)$ (which depends on $N$) such that $Q_{p^{t_0}N}$ has finite projective dimension;  or else, 

\medskip \flushpar{\bf Case 2:} $\operatorname{pd} Q_{p^{t}N}$ is infinite for all $t$. 

In the first case, take $q=p$. It follows from the Theorem of Peskine and Szpiro \cite{\rref{PS}, Thm.~I.7} that $Q_{q^sp^tN}$ has finite projective dimension; and therefore  $\operatorname{syz}_2Q_{q^sp^tN}$ is a free $R$-module, for all $s\ge 0$ and $t\ge t_0$. Thus,    (\tref{syz}) holds.  

Recall from Theorem \tref{!!!} that if $p=2$, then $Q_N$ has finite projective dimension for all $N\ge n$. Thus, $p=2$ is covered in the first case.

Henceforth, we focus on the second case  with   $p\ge 3$.    Take $t_0(N)$ to be zero. 
We first identify an exponent $e$ so that $q=p^e$ has the form $bn+1$ for some positive even integer $b$. 
 Let $q'=p^{\pmb o}$ and write $q'=b'n+1$. If $b'$ is even, then take $e=\pmb o$, $q=q'$ and $b=b'$. If $b'$ is odd, then $n$ must be even since $b'n+1$ is the odd integer $q'$. Observe that
$(q')^2=((b')^2n+2b')n+1$, with $b=(b')^2n+2b'$ even. In this case, take $e=2\pmb o$ and $q=p^e$. 

Take any $t\ge 0$. Write $p^tN=an+r$, with $1\le r\le n-1$. (The remainder $r$ can not be zero because $\operatorname{pd} Q_{p^tN}=\infty$.) It follows that $qp^tN=(bp^tN+a)n+r$. Observe that $a$ and $bp^tN+a$ have the same parity and the remainder $r$ is the same in both numbers $p^tN$ and $qp^tN$.  The modules $Q_{p^tN}$ and $Q_{qp^tN}$ both have infinite projective dimension; so the integers $a$ and   $bp^tN+a$ are in $S_p$ according Theorem \tref{???}. Corollary \tref{I+} shows that $\operatorname{syz}_2 Q_{\pmb k,n,qp^tN}$ and  $\operatorname{syz}_2 Q_{\pmb k,n,p^tN}$ are both isomorphic to a shift of 
$$\cases \operatorname{coker}\varphi_{r,n-r}&\text{if $a$ is odd}\\
\operatorname{coker} \varphi_{n-r,r}&\text{if $a$ is even}. \endcases$$ One may iterate this procedure to see that (\tref{syz}) holds. \qed
\enddemo

Given Data  \tref{D8}, it is natural to ask for a bound $t_0$, which depends only on $p$ and $n$, such that, if $N\ge 1$ is an integer, then 
$$\text{$\operatorname{pd} Q_{p^tN}=\infty$ for some $t\ge t_0$} \implies \text{$\operatorname{pd} Q_{p^tN}=\infty$ for all $t\ge 0$}.\tag\tnum{*}$$
We have some partial results in the direction of describing such a number $t_0$. \item{$\bullet$} If $p=2$, then $\operatorname{pd} Q_{p^tN}<\infty$ for all $t\ge n$; so, any $t_0\ge \log_2 n$ has property (\tref{*}) for vacuous reasons. \item{$\bullet$} If $n=2$ and $p$ is odd, or $n=3$ and $p\equiv 1 \mod 3$, then Propositions \tref{28.2} and \tref{28.4} show that $t_0=0$ has property (\tref{*}). \item{$\bullet$} If
$n=3$ and $p\equiv 2\mod 3$ or $n\ge 4$ and $p$ is odd, then   Proposition \tref{P28.19'} shows that   $t_0=2\pmb o$ satisfies (\tref{*}) for all $N\equiv 1\mod n$. This special case of (\tref{*}) includes the most important case $N=1$, where $Q_N=\pmb k$: if $\operatorname{pd} Q_{p^{2\pmb o}} =\infty$, then  $\operatorname{pd} Q_{p^{t}} =\infty$ for all $t$.

\proclaim{Lemma \tnum{Oct13}} Let $p$ be an odd prime integer, and $n$ be an integer with either  $n\ge 4$ or else $n=3$ and $p\equiv 2\mod 3$. Suppose $q=p^e$ and $q=bn+1$ for some positive integers $b$ and $e$. If $a$ is a non-negative integer with $qa+b$  in $S_p$, then $a$ is even. \endproclaim 
\demo{Proof} Assume $a$ is odd. We prove that $qa+b$ is in $T_p$. We treat two cases. In the first case either $n\ge 4$ or $n=3$ and $q\equiv 2$. In the second case $n=3$ and $q\equiv 1$.
In the first case we have
$$\tsize |(qa+b)-qa|=b=\frac{q-1}n<\{\frac q3\}$$and Remark \tref{R1.0} shows that $qa+b$ is in $T_p$. We justify the inequality $\frac{q-1}n<\{\frac q3\}$. If $4\le n$, then $\frac{q-1}n<\frac{q-1}3\le \{\frac q3\}$. If $n=3$ and $q\equiv 2$, $\frac{q-1}n=\frac{q-1}3<\frac{q+1}3=\{\frac q3\}$.

Now we consider the case where $n=3$, $p\equiv 2\mod 3$, and $q\equiv 1\mod 3$. Let $q'=q/p$. Observe that  $q'\equiv 2 \mod 3$. We see that 
$$\tsize qa+b=qa+\frac{q-1}3= qa+\frac{pq'-1}3= qa+\frac{(p+1)q'-(q'+1)}3=q'(pa+\frac{p+1}3)-\frac {q'+1}3.$$ The integer $pa+\frac{p+1}3$ is odd; so,  $qa+b$ is in $T_p$ by Definition \tref{D5.1}. \qed \enddemo

\proclaim{Observation  \tnum{do28'}} Suppose $p\ge 5$ is a prime integer,  $n\ge 3$ is an integer, and $q=p^e$ for some integer $e\ge 1$. Suppose further that 
 $q=bn+1$ and  $b=\sum_{i=0}^E b_ip^i$  for integers $b_i$  with $|b_i|$ at most  $2\lfloor \frac p3\rfloor$ for all $i$.
Then   $b_e=\dots =b_E=0$.
\endproclaim

\demo{Proof}  Write $b=A+B$, with $A=\sum_{i=0}^{e-1} b_ip^i$ and $B=\sum_{i=e}^E b_ip^i$.  If $b_i\neq 0$ for some $i$ with $e\le i\le E$,   then $|B|$ is equal to $p^e$ times a positive integer; thus, $p^e\le |B|$. On the other hand, $|A|\le 2\lfloor \frac p3\rfloor\frac {p^e-1}{p-1}$; so,
$$\tsize q\le |B|=|(A+B)-A|\le |A+B|+|A|\le \frac{q-1}n+2\left\lfloor \frac p3\right\rfloor\frac {q-1}{p-1} \le (q-1)\left(\frac 1n+\frac 23\right)\le q-1.$$
We used the fact that $A+B=b=\frac {q-1}n$. We also used the fact that  $3$ does not divide $p$, so $\lfloor \frac p3 \rfloor\le \frac {p-1}3$.
This contradiction shows that $b_i=0$ for $e\le i\le E$. 
\qed\enddemo

\proclaim{Proposition \tnum{P28.19'}} Adopt Data \tref{D8} with $p\ge 3$.  Assume  that $n\ge 4$ or else $n=3$ and $p\equiv 2\mod 3$.  Let $q$ be the integer $p^{\pmb o}$. Write $q=bn+1$ and $N=an+1$ for some integers $b$ and $a$. The following statements are equivalent:
\roster
\item  $\operatorname{pd} Q_{p^tN}=\infty$  for all integers $t\ge 0$, 
\item $\operatorname{pd} Q_{q^2N}=\infty$,
\item the integers $a$ and $b$ are even elements of $S_p$.
\endroster

\flushpar Furthermore,  if the above equivalent conditions   are in effect, then $\operatorname{syz}_2 Q_{q^sN}$ is isomorphic to a shift of $\operatorname{syz}_2 Q_{N}$, for all non-negative integers $s$.
\endproclaim

 \demo{Proof}Observe,  by induction, that 
$$  \matrix\format\l\\ q^0N=an+1\\ q^1N=(qa+b)n+1\\q^2N=(q^2a+qb+b)n+1\\
 q^sN=\left[q^sa+\left(\sum_{i=0}^{s-1}q^i\right)b\right]n+1.\endmatrix $$
  $(2)\Rightarrow (3)$.  Theorem \tref{???} shows that the integers $a$, $qa+b$, and $q^2a+qb+b$ all are in $S_p$. 
Apply Lemma \tref{Oct13} to $a$ and to $qa+b$. We have $qa+b$ and $q(qa+b)+b$ are both in $S_p$; so we conclude that $a$ and $qa+b$ are both even. It follows that $b$ is even. 

We next show that $b$ is in $S_p$. 
We use the special $p$-adic expansion of $b$. We treat the cases $p=3$ and $p\ge 5$ separately. We first assume that $p=3$. Use Remark \tref{O-10-13}. The fact that $a$ and $qa+b$ are both even elements of $S_p$ ensures that $a$ and $qa+b$ are divisible by $4$. So, $b$ is also divisible by $4$. Write $$a=4\sum_{i=0}^r \epsilon_i3^i\quad\text{and}\quad b=4\sum_{i=0}^s \delta_i3^i\tag\tnum{exp}$$ for $\epsilon_i\in \{0,1\}$ and $\delta_i\in \{0,1,2\}$. Observe that
$$3^{s+1}\le 4(3^s)\le b=\frac {q-1}{3}<\frac q3=3^{{\pmb o}-1};$$so, $s\le {\pmb o}-2$ and $$qa+b=4(\epsilon_r 3^{r+{\pmb o}}+\dots+ \epsilon_03^{\pmb o}+\delta_s 3^s+\dots +\delta_0).\tag\tnum{qa+b}$$ The fact that $qa+b$ is in $S_p$ says that every coefficient in the expansion (\tref{qa+b}) is in $\{0,1\}$. Thus, every $\delta_i\in \{0,1\}$ and $b\in S_p$.

Now we assume $p\ge 5$. Expand $b$  and $a$ in the $p$-adic expansion as described in Notation \tref{N5.1}. Observation \tref{do28'} shows that
$$b=\sum_{i=0}^{{\pmb o}-1} b_ip^i.$$
The $p$-adic expansion of $qa+b$, in the sense of Notation \tref{N5.1}, is obtained by concatenating the $p$-adic expansion of $qa$ with the $p$-adic expansion of $b$.  
We have shown that $qa+b$ is an even element of $S_p$. Remark \tref{R1.1}  shows that every coefficient of the $p$-adic expansion of $qa+b$, in the sense of Notation \tref{N5.1}, is even. Therefore,  every coefficient of $b$ is even; and  therefore Remark \tref{R1.1} shows that  $b$ are even elements of $S_p$. 
The proof of $(2)\Rightarrow (3)$ is complete.

The implication $(3)\Rightarrow (1)$ now follows readily. If $p=3$, then $a$ and $b$ are given in (\tref{exp}) with all coefficients $\epsilon_i$ and $\delta_i$ from the set $\{0,1\}$. If $p\ge 5$,  every coefficient of the $p$-adic expansions of $a$ and $b$, in the sense of Notation \tref{N5.1}, is even.
For all characteristics $p\ge 3$,
 the $p$-adic expansion of $\lfloor\frac{q^sN}n\rfloor$ is obtained by concatenating the $p$-adic expansions of $q^sa,q^{s-1}b,\dots qb,b$. The form of the resulting 
$p$-adic expansion shows that 
$\lfloor\frac{q^sN}n\rfloor\in S_p$;  and therefore, $Q_{q^sN}$ has infinite projective dimension by Theorem \tref{I-MAIN}.  The theorem of Peskine and Szpiro \cite{\rref{PS}, Thm.~I.7} ensures that $Q_{p^tN}$ has infinite projective dimension for all non-negative integers $t$.

The final assertion  is an immediate consequence of Theorem 3.5 where it is shown that
$\operatorname{syz}_2 Q_{q^sN}$ is a shift of the cokernel of $\varphi_{n-1,1}$. 
\qed
\enddemo

\heading Subsection \number\SectionNumber.3. \quad Can one use socle degrees to predict that the tail of the resolution of $Q_{\pmb k,n,p^tN}$ is a shift of the tail of the resolution of  $Q_{\pmb k,n,N}$?
\endheading

We denote the socle of $Q$ by $\operatorname{soc} Q$, see Definition \tref{num?}; and, if $\Bbb F$ is a complex, then $\Bbb F_{\ge i}$ is the truncation $\dots \to F_{i+1}\to F_i$ of $\Bbb F$, see (\tref{trun}). 
\proclaim {Theorem \tnum{T29}} Let $\pmb k$ be a field and $n$, $N_1$, and $N_2$ be positive integers. Write $R$ for the ring $R_{\pmb k,n}$ and $Q_{N_i}$ for the $R$-module $Q_{\pmb k,n,N_i}$. Assume that $Q_{N_1}$ and $Q_{N_2}$ both have infinite projective dimension over $R$. Let $\Bbb F_{i,\bullet}$ be the minimal homogeneous resolution of $Q_{N_i}$ by free $R$-modules and let $N_i=a_in+r_i$ with $1\le r_i\le n-1$. The following statements are equivalent:
\roster
\item   $\operatorname{soc} Q_{N_2}$  is isomorphic to a shift of $\operatorname{soc} Q_{N_1}$ as a graded vector space,
\item   $\Bbb F_{2,\ge 2}$  is isomorphic to a shift of $\Bbb F_{1,\ge 2}$ as a graded $R$-{\bf module},
\item   $\Bbb F_{2,\ge 2}$ is isomorphic to a shift of $\Bbb F_{1,\ge 2}$ as a graded {\bf complex},
\item the $R$-module  $\operatorname{syz}_2Q_{N_2}$ is isomorphic to a shift of $\operatorname{syz}_2Q_{N_1}$, and 
\item either
{\rm (a)} $a_1+a_2$ is even and $r_1=r_2$; or else,
 {\rm (b)} $a_1+a_2$ is odd and $r_1+r_2=n$.
\endroster 
\flushpar Furthermore, if $N_2=qN_1$ for some positive integer $q$ 
 and conditions {\rm(1)--(5)} occur, then 
$$\operatorname{soc}   Q_{N_2}\cong \operatorname{soc}   Q_{N_1}(-w),\quad \operatorname{syz}_2  Q_{N_2}\cong \operatorname{syz}_2  Q_{N_1}(-w),\quad\text{and}\quad 
\Bbb F_{2,\ge 2}\cong \Bbb F_{1,\ge 2}(-w)$$
for $w=\frac 32N_1(q-1)$.
\endproclaim

\remark{Remark} The implication $(1)\Rightarrow(2)$ is \cite{\rref{KU}, Thm.~1.1}. We have an explicit formula for the socle degrees of $Q_{N_i}$ (see Theorem \tref{I-MAIN}); so, we may verify hypothesis (c) of \cite{\rref{KU}, Thm.~1.1} directly. The parameter $b+2a$ from \cite{\rref{KU}} is $3N_i+2n-6$ in the present notation.
The techniques of \cite{\rref{KU}} were not able to prove that $(1)\Rightarrow(3)$. One of the motivations for the present work was to establish this implication. The shift $\frac 32 N_1(q-1)$ is identified in \cite{\rref{KU}, Cor.~2.1}.
  \endremark

\demo{Proof} Assertions (3) and (4) are equivalent because $\Bbb F_{i,\ge 2}$ is the minimal resolution of $\operatorname{syz}_2Q_{N_i}$. It is clear that $(3)\Rightarrow (2)$ because $(3)$ is a statement about graded Betti numbers and differentials; whereas $(2)$ is only a statement about graded Betti numbers.

Apply Theorem \tref{???} to see that $a_i\in S_c$, where $c$ is the characteristic of $\pmb k$, and apply Theorem \tref{I-MAIN} to read the socle degrees of $Q_{N_i}$ and the resolution $\Bbb F_{i,\bullet}$ from the data $N_i=a_in+r_i$. The assertion $(2)\Rightarrow (1)$ may be read from Theorem \tref{I-MAIN} where it is shown that 
$$\operatorname{soc} Q_{N_i}\cong \frac{\Bbb F_{i,3}}{(x,y,z)\Bbb F_{i,3}}(+3)\tag\tnum{Keep}$$ as a graded vector space. 

\flushpar $(1) \Rightarrow (5)$ The socle degrees of $Q_{N_i}$ have the form $D_i\:\! 3, d_i\:\! 1$ for $(D_i,d_i)$ given in Theorem \tref{I-MAIN}. Condition (1) asserts that there exists an integer $w$ with $D_2=D_1+w$ and  $d_2=d_1+w$. It follows that $$D_2-d_2=D_1-d_1.\tag\tnum{taf}$$ We see from 
Theorem \tref{I-MAIN} that
$$D_i-d_i=\cases -n+2r_i &\text{if $a_i$ is odd}\\  n-2r_i &\text{if $a_i$ is even}.\endcases$$If $a_1$ and $a_2$ have the same parity, then (\tref{taf}) shows that $r_1=r_2$. If $a_1$ and $a_2$ have the different parity, then (\tref{taf}) shows that $n=r_1+r_2$.

\flushpar $(5)\Rightarrow (4)$ If hypothesis (a) holds, then Theorem \tref{I-MAIN} shows that $\operatorname{syz}_2Q_{N_1}$ and $\operatorname{syz}_2Q_{N_2}$ are both isomorphic to a shift of 
$$\cases \operatorname{coker} \varphi_{n-r_1,r_1}&\text{if $a_1$ is even}\\\operatorname{coker} \varphi_{r_1,n-r_1}&\text{if $a_1$ is odd}.\endcases$$ If hypothesis (b) holds with $a_1$ even, then $$N_2=(a_2-1)n+(n-r_2)=(a_2-1)n+r_1,$$ with $a_2-1$ even and Theorem \tref{I-MAIN} shows that $\operatorname{syz}_2Q_{N_1}$ and $\operatorname{syz}_2Q_{N_2}$ are both isomorphic to a shift of $\operatorname{coker} \varphi_{n-r_1,r_1}$.

Finally, we assume (1) -- (5) hold. These conditions guarantee that there exist  integers $w_1$, $w_2$, $w_3$  with
$$\operatorname{soc}   Q_{N_2}\cong \operatorname{soc}   Q_{N_1}(-w_1),\ \ \operatorname{syz}_2  Q_{N_2}\cong \operatorname{syz}_2  Q_{N_1}(-w_2),\ \ \text{and}\ \  
\Bbb F_{2,\ge 2}\cong \Bbb F_{1,\ge 2}(-w_3).$$ The complex $\Bbb F_{i,\ge 2}$ is a resolution of $\operatorname{syz}_2 Q_{N_i}$; so $w_2=w_3$, and the isomorphism (\tref{Keep}) shows that $w_1=w_3$.  We now identify the common value $w_1=w_2=w_3$, when $N_2=qN_1$. Condition (5) is in effect. If $a_1$ and $a_2$ are both even, then Theorem \tref{I-MAIN} gives 
$$\tsize w_1=\frac 32(a_2n-a_1n)= \frac 32([a_2n+r]-[a_1n+r])=\frac 32(qN_1-N_1).$$
If $a_1$ and $a_2$ are both odd, then  
$$\tsize w_1=\frac 32([a_2+1]n-[a_1+1]n)= \frac 32([a_2n+r]-[a_1n+r])=\frac 32(qN_1-N_1).$$
If $a_1$ is even and $a_2$ is odd, then  
$$\tsize w_1=\frac 32(a_2+1)n-(\frac 32a_1n+3r)= \frac 32([(a_2+1)n-r]-[a_1n+r])=\frac 32(qN_1-N_1). \qed$$
\enddemo

\SectionNumber=9\tNumber=1
\heading Section \number\SectionNumber. \quad Two variables.
\endheading

Consider $$\bar R_{\pmb k,n}=\pmb k[x,y]/(x^n+y^n)\quad\text{and}\quad \bar Q_{\pmb k,n,N}=\bar R_{\pmb k,n}/(x^N,y^N), \tag\tnum{bar}$$ where $\pmb k$ is a field. 
The results in this section were announced in the Introduction to \cite{\rref{KU}}. In particular, Example \tref{last} exhibits a situation where the tail of the resolution of $F^t(\bar Q_{\pmb k,n,N})$, after shifting, is periodic as a function of $t$, with an arbitrarily large period. We have used $F$ to represent the Frobenius functor.

Observation \tref{.O1.2} and Corollary \tref{C2.4} are the two variable analogues of Theorem \tref{I-MAIN}.   Observation \tref{.O1.22} is the two variable analogue of Theorem \tref{T29}. Two phenomenon  occur in 3 variables that do not occur in 2 variables. Let $N=an+r$, with $0\le r\le n-1$. In 3 variables, $a$ plays an important role; in 2 variables, $a$ is irrelevant. In 2 variables $\operatorname{pd}_{\bar R_{\pmb k,n}}\bar Q_{\pmb k,n,N}$ is infinite whenever $n$ does not divide $N$; in 3 variables $Q_{\pmb k,n,N}$ often has finite projective dimension over $R_{\pmb k, n}$. 
\proclaim{Observation \tnum{.O1.2}}Let $x,y$ be a regular sequence in a ring $P$. Let $a$, $n$, and $r$ be fixed integers with $0\le r\le n-1$ and $0\le a$, and  let  $f=x^n+y^n$,  and $N=an+r$.  Then the  ideal  $K=(x^N,y^N,f)$ of $P$ is perfect of grade two and $P/K$ is resolved by
$$0\to P^2@> d_2>> P^3@> d_1 >> P\to P/K\to 0,$$
with
$$d_1=\bmatrix x^N&y^N&f\endbmatrix\quad\text{and}\quad
d_2=\bmatrix x^{n-r}&(-1)^{a-1}y^r\\(-1)^ay^{n-r}&x^r\\  L&M\endbmatrix,$$for
$$L=-\sum\limits_{i=0}^a(x^n)^{a-i}(-y^n)^i\quad\text{and}\quad M= (-1)^{a}x^ry^r\sum\limits_{i=0}^{a-1}(x^n)^{a-1-i}(-y^n)^i.$$
\endproclaim
\demo{Proof}It suffices to show that the entries of $d_1$ are the signed maximal order minors of $d_2$. (See (\tref{smom}) if necessary.) Let $\Delta_i$ be $(-1)^{i+1}$ times the determinant of $d_2$ with row $i$ deleted.    We see that $\Delta_3=f$,
$$\split \Delta_2&{}=-x^{n-r}M+(-1)^{a-1}y^rL\\&{}=(-1)^{a+1}x^ny^r\sum\limits_{i=0}^{a-1}(x^n)^{a-1-i}(-y^n)^i+(-1)^ay^r\sum\limits_{i=0}^a(x^n)^{a-i}(-y^n)^i\\&{}=(-1)^{a}y^r\left(-\sum\limits_{i=0}^{a-1}(x^n)^{a-i}(-y^n)^i+\sum\limits_{i=0}^a(x^n)^{a-i}(-y^n)^i\right)=(-1)^{a}y^r\left(  (-y^n)^a\right)\\&{}=y^N,\endsplit $$ and
$$\split \Delta_1&{}= (-1)^ay^{n-r}M-x^rL\\&{}=x^ry^n\sum\limits_{i=0}^{a-1}(x^n)^{a-1-i}(-y^n)^i+x^r\sum\limits_{i=0}^a(x^n)^{a-i}(-y^n)^i\\&{}=x^r\left(-\sum\limits_{i=0}^{a-1}(x^n)^{a-1-i}(-y^n)^{i+1}+\sum\limits_{i=0}^a(x^n)^{a-i}(-y^n)^i\right)\\&{}=x^r\left(-\sum\limits_{i=1}^{a}(x^n)^{a-i}(-y^n)^{i}+\sum\limits_{i=0}^a(x^n)^{a-i}(-y^n)^i\right)\\&{}=x^N. \qed \endsplit$$\enddemo

\proclaim{Corollary \tnum{C2.4}}Keep the notation of Observation \tref{.O1.2}. Assume that $f$ is a regular element in $P$. Let $R=P/(f)$, $$D=\bmatrix x^{n-r}&(-1)^{a-1}y^r\\(-1)^ay^{n-r}&x^r\endbmatrix,\quad\text{and}\quad \check D=\bmatrix x^r&(-1)^{a}y^r\\(-1)^{a+1}y^{n-r}&x^{n-r}\endbmatrix.$$
Then the $R$-resolution of $R/(x^N,y^N)R$ is
$$\dots @>\check D>>R^2@> D>> R^2@> \check D>> R^2@> D>> R^2@> \bmatrix x^N&y^N\endbmatrix >> R\to R/(x^N,y^N)R\to 0.$$ \endproclaim
\demo{Proof}
Notice that $D$ is the top two rows of $d_2$ and that $\check D$ is the classical adjoint of $D$; so $D\check D=\check D D=fI_2$ in $P$. It is clear that the proposed resolution is a complex.  If $Dv=fw$ in $P$, then multiply by $\check D$ to learn that $fv=f\check D w$; but $f$ is a regular element in $P$; so $v=\check D w$. One may also reverse the roles of  $D$ and $\check D$. If $\bmatrix x^N&y^N\endbmatrix v=af$ in $P$, then $\left[\smallmatrix v\\-a\endsmallmatrix\right]$ is in $\operatorname{ker} d_1=\operatorname{im} d_2$; so there exists $u$ in $P^2$ with 
$$\bmatrix D\\ \vspace{.05in}\noalign{\hrule} \\ L\ M\endbmatrix u= \bmatrix v\\-a\endbmatrix;$$in particular, $Du=v$. \qed
\enddemo

\proclaim{Observation \tnum{.O1.22}}  Let $\pmb k$ be a field and $n$, $N_1$, and $N_2$ be positive integers. Write $ R$ for the ring $\bar R_{\pmb k,n}$ and $ Q_{N_i}$ for the $ R$-module $\bar Q_{\pmb k,n,N_i}$, as described in {\rm(\tref{bar})}.   Let $\Bbb F_{i,\bullet}$ be the minimal homogeneous resolution of $ Q_{N_i}$ by free $ R$-modules. Then the following statements are equivalent:
\roster
\item   $\operatorname{soc}  Q_{N_2}$  is isomorphic to a shift of $\operatorname{soc}  Q_{N_1}$ as a graded vector space,
\item   $\Bbb F_{2,\ge 1}$  is isomorphic to a shift of $\Bbb F_{1,\ge 1}$ as a graded $ R$-{\bf module},
\item   $\Bbb F_{2,\ge 1}$ is isomorphic to a shift of $\Bbb F_{1,\ge 1}$ as a graded {\bf complex},
\item the $ R$-module  $\operatorname{syz}_1 Q_{N_2}$ is isomorphic to a shift of $\operatorname{syz}_1 Q_{N_1}$, and 
\item $N_1\equiv \pm N_2 \mod n$.
\endroster 
\flushpar Furthermore, if $N_2=qN_1$ for some positive integer $q$ and conditions {\rm(1)--(5)} occur, then 
$$\operatorname{soc}  Q_{N_2}\cong \operatorname{soc}  Q_{N_1}(-w),\quad \operatorname{syz}_1 Q_{N_2}\cong \operatorname{syz}_1 Q_{N_1}(-w),\quad\text{and}\quad 
\Bbb F_{2,\ge 1}\cong \Bbb F_{1,\ge 1}(-w)$$
for $w=N_1(q-1)$.
\endproclaim

\demo{Proof}Let $N_i=a_in+r_i$ with ${1\le r_i\le n-1}$, $f$ be the polynomial $x^n+y^n$ in $P=\pmb k[x,y]$, and $K_i$ be the ideal $(x^{N_i},y^{N_i},f)$   of $P$. Observation \tref{.O1.2} shows  that the $P$-resolution of $P/K_i$ is
$$0\to P(-N_i-n+r_i)\oplus P(-N_i-r_i) \to P(-N_i)^2\oplus P(-n) \to P.$$   It follows from Remark \tref{R?} that the socle degrees of $P/K_i= Q_{N_i}$ are
$$\{N_i+n-r_i-2,N_i+r_i-2\}.\tag\tnum{.}$$  
We also know from Corollary \tref{C2.4} that the first syzygy module of the $R$-modules $ Q_{N_i}$ is presented by
$$ 
\matrix R(-N_i-n+r_i)\\\oplus\\R(-N_i-r_i)\endmatrix @>D_i >> R(-N_i)^2@> >> \operatorname{syz}_1^R( Q_{N_i})\to 0,\tag\tnum{..}$$ 
with $$D_i={\bmatrix x^{n-r_i}& -y^{r_i}\\y^{n-r_i}&x^{r_i}\endbmatrix}.$$
We read from (\tref{.}) and (\tref{..}) that
$$\operatorname{soc} Q_{N_i}\cong \frac{\Bbb F_{i,2}}{(x,y)\Bbb F_{i,2}}(+2)\tag\tnum{Keep.}$$ as a graded vector space. 

We have $(4)\iff (3)\Rightarrow (2) \Rightarrow (1)$ exactly as in the proof of Theorem \tref{T29} once (\tref{Keep}) is replaced with (\tref{Keep.}). The socle degrees of $Q_{N_i}$ are $\{d_i,d_i'\}$, as given in (\tref{.}). It follows that 
$$\split (1)&{}\implies |d_1-d_1'|=|d_2-d_2'|\iff |n-2r_1|=|n-2r_2|\\&{}\implies r_1=r_2 \text{ or } r_1+r_2=n\iff N_1\equiv \pm N_2 \mod n\iff (5).\endsplit$$

\flushpar$(5)\implies (4)$ The $R$-module $\operatorname{syz}_1(Q_{N_i})$ is presented by $D_i$ as shown in (\tref{..}). We see that
$$D_2=\cases D_1&\text{if $r_1=r_2$}\\
J^{\text{\rm T}}D_1J  
&\text{if $r_1+r_2=n$},\endcases$$
where $J$ is the matrix $J=\left[\smallmatrix 0&-1\\ 1&0\endsmallmatrix\right]$. 
Thus, if (5) holds, then 
$$\operatorname{syz}_1(Q_{N_2})=\operatorname{syz}_1(Q_{N_1})[N_1-N_2].$$ 
Finally, we assume (1) -- (5) hold. These conditions guarantee that there exist  integers $w_1$, $w_2$, $w_3$  with
$$\operatorname{soc}   Q_{N_2}\cong \operatorname{soc}   Q_{N_1}(-w_1),\ \ \operatorname{syz}_1  Q_{N_2}\cong \operatorname{syz}_1  Q_{N_1}(-w_2),\ \ \text{and}\ \  
\Bbb F_{2,\ge 1}\cong \Bbb F_{1,\ge 1}(-w_3).$$ The complex $\Bbb F_{i,\ge 1}$ is a resolution of $\operatorname{syz}_1  Q_{N_i}$; so $w_2=w_3$, and the isomorphism (\tref{Keep.}) shows that $w_1=w_3$.  We have already seen that when  condition (5) is in effect, then $w_2=N_2-N_1$. \qed 
\enddemo

\example{Example \tnum{last}} Fix a positive integer $e$ and a field $\pmb k$ of characteristic $p>0$. Take $n=p^e+1$ and $N$ to be any integer which is relatively prime to $n$. Let $R$ be the ring $\bar R_{\pmb k,n}$ and for each integer $M$, let $Q_M$ be the $R$-module $\bar Q_{k,n,M}$ as described in (\tref{bar}). Let $t_1\le t_2$ be positive integers. We observe that 
$$\matrix \format\l\\\text{$\operatorname{syz}_1 Q_{p^{t_2}N}$ is isomorphic to a shift of $\operatorname{syz}_1 Q_{p^{t_1}N}$ if and only if the integer}\\\text{$t_2-t_1$ is divisible by $e$.}\endmatrix \tag\tnum{div}$$ 
In other words, $\operatorname{syz}_1 Q_{p^{i}N}$ represent different isomorphism classes, even after shifting, for $0\le i\le e-1$; but $\operatorname{syz}_1Q_{p^{k+e}N}$ is isomorphic 
to a shift of  $\operatorname{syz}_1Q_{p^{k}N}$ for all integers $k\ge 0$. 

Assertion (\tref{div}) follows from Observation \tref{.O1.22}, where it is shown that $\operatorname{syz}_1 Q_{p^{t_2}N}$ is isomorphic to a shift of $\operatorname{syz}_1 Q_{p^{t_1}N}$ if and only if $p^{t_2}N\equiv \pm p^{t_1}N\mod n$. We have arranged the data so that $p^e\equiv -1\mod n$; but if $0\le t_1<t_2\le e-1$, then $p^{t_2}$ is not congruent to $\pm p^{t_1}$ mod $n$.
\endexample 

\enddocument